\documentclass[11pt,a4paper]{report}
\usepackage[latin1]{inputenc}
\usepackage[english]{babel}
\usepackage{amssymb,amsmath,amscd,amstext}
\usepackage{amsfonts}
\usepackage{amsthm}
\usepackage[dvips]{graphicx}
\usepackage{psfrag}
\usepackage{epsfig}
\usepackage{palatino}
\usepackage{xypic}
\input{xy}
\xyoption{all}
\newcounter{eno}
\newcounter{ch}
\newcounter{sec}
\newcounter{no}
\def \sec#1 {\stepcounter{sec} \bigskip\par{\underline{\bf \S \arabic{ch}.\arabic{sec} {\bf #1}  \setcounter{eno}{1} \setcounter{no}{1}}}}
\def \prop#1 {\bigskip\par{\bf Proposition \arabic{ch}.\arabic{sec}.\arabic{no}} \stepcounter{no} {\it #1}}
\def \cor#1 {\bigskip\par{\bf Corollary \arabic{ch}.\arabic{sec}.\arabic{no}} \stepcounter{no} {\it #1}}
\def \lem#1 {\bigskip\par{\bf Lemma \arabic{ch}.\arabic{sec}.\arabic{no}} \stepcounter{no} {\it #1}}
\def \thm#1 {\bigskip\par{\bf Theorem \arabic{ch}.\arabic{sec}.\arabic{no}} \stepcounter{no} {\it #1}}
\def \conj#1 {\bigskip\par{\bf Conjecture \arabic{ch}.\arabic{sec}.\arabic{no}} \stepcounter{no} {\it #1}}
\def \define#1 {\bigskip\par{\bf Definition
\arabic{ch}.\arabic{sec}.\arabic{no}} \stepcounter{no} {\it #1}\hspace{2mm}}
\def \eg#1 {\bigskip\par{\bf Example \arabic{ch}.\arabic{sec}.\arabic{no}} \stepcounter{no} {\it #1}}

\newcommand {\Hq} {\mathbb H}
\newcommand {\R} {\mathbb R}
\newcommand {\C} {\mathbb C}
\newcommand {\Z} {\mathbb Z}
\newcommand {\Q} {\mathbb Q}
\newcommand {\ZG} {\mathbb {Z}[G]}
\newcommand {\QG} {\mathbb {Q}[G]}

\newcommand {\RG} {\mathbb {R}[G]}
\newcommand {\ZH} {\mathbb {Z}[H]}
\newcommand {\rk} {{\rm rk}_{\Z}}
\newcommand {\ZD} {\mathbb {Z}[D_{4n}]}

\newcommand {\Zd} {\mathbb {Z}[D_{8}]}
\newcommand {\QD} {\mathbb {Q}[D_{4n}]}

\newcommand {\Zh} {\mathbb {Z}[\frac {1}{2}]}

\newcommand {\ZhD} {R}
\newcommand {\ZhDf} {T}

\newcommand {\D} {D_{8n+4}}

\def \proof#1 {\bigskip \hspace{2mm}Proof:$\,\,\,$  #1 \newline ${}$\hfill $\Box$
\,\, \newline \bigskip}

\begin{document}
\thispagestyle{empty}
\vspace*{2cm}
\begin{center}{\Huge
LOW DIMENSIONAL ALGEBRAIC COMPLEXES OVER INTEGRAL GROUP RINGS}\\[3cm]
{\LARGE Wajid Hassan Mannan\linebreak
D\,e\,p\,a\,r\,t\,m\,e\,n\,t \,o\,f
\,M\,a\,t\,h\,e\,m\,a\,t\,i\,c\,s\linebreak
U\,n\,i\,v\,e\,r\,s\,i\,t\,y \,C\,o\,l\,l\,e\,g\,e \,L\,o\,n\,d\,o\,n\linebreak
U\,n\,i\,v\,e\,r\,s\,i\,t\,y \,o\,f \,L\,o\,n\,d\,o\,n}\\[3cm]
{\Large A thesis submitted for the degree of\linebreak
\emph{Doctor of Philosophy.}\linebreak
Supervisor:\ \ Prof.\! F.E.A\! Johnson \linebreak
J\,U\,L\,Y \,2\,0\,0\,6}
\end{center}

\newpage

\begin{center}
I, Wajid Hassan Mannan, confirm that the work presented in this thesis is my own. Where information 
has been derived from other sources, I confirm that this has been indicated in the thesis.
\end{center}

\newpage
\thispagestyle{empty}
\vspace*{-2cm}
\begin{center}
{\bf Abstract}\\
\end{center}

{\small \hspace{2mm} The realization problem asks:  When does an
algebraic complex arise, up to homotopy, from a geometric complex?  In the case of $2$- dimensional algebraic complexes, this is equivalent to
the D2 problem, which asks when homological methods can distinguish between 2 and 3 dimensional complexes.  

\bigskip
\hspace{2mm}
We approach the realization problem (and hence the D2 problem) by
classifying all possible algebraic 2- complexes and showing that they are realized.  We show that if a dihedral group has order $2^n$, 
then the algebraic complexes over it are parametrized by their second homology groups, which we refer to as algebraic second homotopy
groups.  A cancellation theorem of Swan (\cite{Swan2}), then allows us to solve the realization problem for the group $D_8$.

\bigskip
\hspace{2mm}Let $X$ be a finite geometric $2$- complex.  Standard isomorphisms
give $\pi_2(X) \cong H_2(\tilde{X};\Z)$, as modules over $\pi_1(X)$.  Schanuel's 
lemma may then be used to show that the stable class 
of $\pi_2(X)$ is determined by $\pi_1(X)$.  We show how $\pi_3(X)$ may be calculated similarly. Specifically, we show that as a module 
over the fundamental group, 
$\pi_3(X)=S^2(\pi_2(X))$, where $S^2(\pi_2(X))$ denotes the symmetric part of the module $\pi_2(X) \otimes_\Z \pi_2(X)$.  As a 
consequence, we are able to show that when the order of
$\pi_1(X)$ is odd, the stable class of $\pi_3(X)$ is also determined by $\pi_1(X)$.

\bigskip
\hspace{2mm}Given a closed, connected, orientable $5$- dimensional manifold,
with finite fundamental group, we may represent it, up to homotopy equivalence, by an algebraic complex.  Poincare duality induces a homotopy
equivalence between this algebraic complex and its dual.  We consider how similar this homotopy equivalence may be made to the identity,
(through appropriate choice of algebraic complex).  We show that it can be taken to be the identity on $4$ of the $6$ terms of the chain 
complex.  However, by finding a homological obstruction, we show that in general the homotopy equivalence may not be written as the identity.
}

\newpage
\thispagestyle{empty}
\noindent{\huge {\bf Acknowledgments}}
\vspace{1.4cm}

\noindent As with anything, my first thanks are to God.  This thesis is largely born out of my late father's love and respect for
academia, as well as my mother's care and support whilst producing it.  I would also like to mention all my extended siblings, as well as my
nieces, who have been my inspiration: Amirah, Tasnim, Samirah, Layla, and Rohan. 

\bigskip
I have been very fortunate in having a supervisor who has been constantly available over the last few years, and has stayed
interested and involved in everything I have done.

\bigskip
I would like to thank my sponsors EPSRC and the John Hawkes Scholarship Fund. 
Their generosity has enabled me to pursue my long held ambition of researching
mathematics.

\bigskip
Finally, I would like to thank the administrators for the departmental computers, especially Jon Linkins.  Also in this context, I would like
to thank Asad Khair for building and maintaining a computer for me to work on at home.

\newpage
\thispagestyle{empty}
\vspace*{5cm}
\begin{center}
{\Large
\emph{In memory of my late father, who inspired all his children to value study and learning over other worldly endeavors,\\[1cm] }}
\end{center}

\tableofcontents



\newpage
\setcounter{sec}{0}
\setcounter{ch}{0}
\part{Introduction}

\hspace{2mm}  There are three main theorems in Chapter One.  Theorem 1.1.1 states that any two projective 
resolutions, of a module, having equal finite length, may be stabilized to the same homotopy type.  This is stronger than Schanuel's lemma,
which states that they have the same final homology group.

\hspace{2mm} Theorem 1.4.11 states that the only groups of period 2 are cyclic.  This result is well 
known though an explicit proof is hard to find in
the literature.   Swan makes a 
comment outlining an argument in \cite{Swan}.  We give a more direct proof which is elementary and avoids the technical difficulties of Swan's argument.

\hspace{2mm} Theorem 1.5.1 says that an initial segment of
a partial free resolution of a module may be changed without altering homotopy type, possibly at the cost of introducing a stably free module
to the resolution.  This will allow us to parametrize certain algebraic 
complexes by the last map in their sequence. 

\hspace{2mm}The question of when an $n+1$- dimensional CW- complex is homotopic to an $n$ dimensional one has been addressed by C.T.C. Wall 
\cite{Wall}.  He showed that for $n
\neq 2$, the vanishing of $n+1^{th}$ cohomology over all coefficient bundles is sufficient.

\hspace{2mm} Wall's methods are not effective for the case $n=2$, which is called the D(2) problem.  The problem is
parametrized by the fundamental group of the CW- complex in question.

\hspace{2mm}  Chapter Two is concerned with the D(2) property for dihedral groups.  In \cite{John1} (62.3) it is shown that the 
dihedral groups $D_{4n+2}$ satisfy the D(2) property.  The smallest dihedral group
not covered by this is $D_8$.  Theorem 2.3.4 states that the D(2) property does hold for $D_8$.  

\hspace{2mm}More generally, for dihedral groups of order $4n$, we show that a minimal 
element of $\Omega_3(\Z)$ is realized as the $\pi_2$ of a presentation (proposition
2.3.2).  In \S2.4 we parametrize all possible minimal elements of $\Omega_3(\Z)$ by a finite group.

\hspace{2mm}In the case of dihedral groups of order $2^n$, $n \in \Z$, we are further able to show that up to chain homotopy equivalence there is a unique algebraic 
2- complex with a "standard" $\pi_2$ (theorem 2.2.11).

\hspace{2mm}In Chapter Three, theorem 3A states that given a geometric $2$- complex, $X$,
with finite fundamental group $G$, we
have $\pi_3(X) \cong S^2(J)$, where $J=\pi_2(X)$.  We will define a 
module over
$\ZG$, $V_G$ and show that $\pi_3(X)$ is determined by $G$, up to stabilization by 
copies of $\ZG$ and copies of $V_G$ (theorem 3.5.5).  Rationally, we show that 
$\pi_3(X) \otimes \Q \cong \QG^a \oplus (V_G \otimes \Q)^b$ for integers $a,b$
(theorem 3.6.5).

\hspace{2mm}In the case where $G$ is a group of odd order, we have 
$V_G \cong 0$.  Hence in this case, the stable class of $\pi_3(X)$ is determined by $G$, and
$\pi_3(X)$ is rationally free (corollaries 3.5.6 and 3.6.6).

\hspace{2mm}Let $M$ be a closed, connected, orientable $5$- dimensional manifold,
with finite fundamental group $G$ (we assume manifolds to be without boundary).  In chapter four we consider algebraic
complexes $C_*(\tilde{M'})$, where $M'$ is a finite CW- complex, with $M \sim M'$.  $C_*(\tilde{M'})$ must satisfy 
Poincare duality.  We use this to show 
that up to chain homotopy equivalence, we may represent it by an algebraic $2$- complex, $\mathcal{A}$, connected to its dual via a 
$G$- invariant bilinear form, $\beta$,
on $(\pi_2(\mathcal{A}))^*$.  We denote the resulting algebraic $5$- complex $(\mathcal{A}, \beta)$.

\hspace{2mm}  We next consider the homotopy equivalence induced by Poincare Duality.  In particular we are interested in how similar it can
be made to the identity.  We show that it can be taken as the identity on $4$ of the $6$ terms of the chain complex.  However, we 
find a homological obstruction to this homotopy equivalence actually being the identity.  In particular, certain
manifolds described in \cite{Bard} do not satisfy the homological condition necessary, for being able to write the homotopy equivalence as
the identity.

\hspace{2mm}  We now make some notational points:

\hspace{2mm} All rings are assumed to contain a multiplicative identity. 
Modules are assumed to be right modules over the relevant ring.  A map between modules is assumed to be linear over the relevant ring, unless
we describe it as a map of sets.

\hspace{2mm} Suppose $A$ and $B$ are modules. If $\alpha:A\to B$ is a map, and $B$ is a summand of a third module $C$, then 
$\alpha:A \to C$ will denote $\alpha:A \to B$ 
composed with inclusion of the summand $B$, in $C$.

\hspace{2mm} Again, let $A$, $B$, $C$ be modules.  If $f:A \to C$ and $g:B \to C$ are maps, then $f \oplus g:A \oplus B \to C$ is 
defined by $(f\oplus g)(a\oplus b)= f(a)+g(b)$.

\hspace{2mm}  The abbreviation f.g. will be used to denote finitely generated, whether in the context of a module or a ring.

\hspace{2mm}  Contrary to some conventions, the subscript on a group will denote its order.  Hence $D_n$ will
denote the dihedral group of order $n$, $S_n$ will denote the symmetric group
of order $n$ and $A_n$ will denote the alternating group of order $n$.

\hspace{2mm}  Let $R$ be a ring and $G$ a group.  $R[G]$ denotes the free $R$- module with group elements as basis.  It is a ring, with 
product
structure given by group multiplication, for the basis elements, and extended linearly over $R$, for the remaining elements.  The "group ring" of a
group $G$, will refer to $\ZG$.  

\hspace{2mm}  In the context of a $\ZG$- module, $\Z$ will denote the $\ZG$- module whose underlying Abelian group is $\Z$,
and on which the action of $G$ is trivial.

\hspace{2mm}  The map sending $\sum{\lambda_i g_i}$ to $\sum{\lambda_i}$, for $\lambda_i \in \Z, \,\,\, g_i \in G$ will be referred to as
augmentation, and will usually be denoted $\epsilon$.  Its kernel, denoted $IG$ will be called "the augmentation ideal".

\hspace{2mm} The dual of a module, $M$, denoted $M^*$ is the set of $\Z$- linear maps
$M \to \Z$.  $M^*$ has the structure of an Abelian group with respect to point-wise addition.  In fact, it is a module over
$\ZG$, with $G$ action given by $\alpha g(m)= \alpha(m g^{-1})$ for all $m \in M$, $\alpha \in M^*$ and $g \in G$.  Similarly  we define 
the dual of a map, $f$, to be precomposition with $f$.  The dual of a map, $f$, will be denoted by $f^*$.  

\hspace{2mm} As we work over finite groups, this definition is consistent with the one where $M^*$ is defined as the set of $\ZG$- linear
maps $M \to \ZG$.

\hspace{2mm} If $\beta:M \times M \to \Z$ is a $G$- invariant bilinear form on a module $M$, it will also be regarded as
a map $M \to M^*$, which sends $x \in M$ to the element of $M^*$ defined by  
$\beta(x)(y)=\beta(x,y)$, for all $y \in M$.

\newpage
\part{}\label{}
\newpage

\chapter{Algebraic methods}\label{}
\stepcounter{ch}
\setcounter{sec}{0}

\hspace{2mm}  There are three main theorems in this chapter.  Theorem 1.1.1 states that any two projective 
resolutions of a module, having equal finite length, may be stabilized to the same homotopy type.  Theorem 1.4.11 states that the only groups of period 2 are cyclic. 
Although this result is well known, the proof is technically simpler than any in publication.  Theorem 1.5.1 says that an initial segment of
a partial free resolution of a module may be changed without altering homotopy type, possibly at the cost of introducing a stably free module
to the resolution.  This will allow us to parametrize certain algebraic 
complexes by the last map in their sequence.

\sec{Schanuel's Lemma}
 
\bigskip
\hspace{2mm}  Schanuel's lemma plays a key role in our algebraic study of
homotopy.  We adapt the proof to give us a stronger result; any two projective 
resolutions of a module may be stabilized to the same homotopy type.

\hspace{2mm}  Let $R$ be a ring.  An algebraic complex over $R$
consists of a sequence of modules $C_i$ and maps $\delta_i:C_i \to C_{i-1}$,
such that $\delta_i \delta_{i-1}=0$ for each $i$.  It may be denoted $(C_i,\,
\delta_i)$.  

\hspace{2mm}  A chain map $f:(C_i,\,\delta_i) \to (D_i,\,\partial_i)$, consists
of a sequence of maps $f_i:C_i \to D_i$ such that for each $i$, $f_{i-1}
\delta_i = \partial_i f_i$.

\hspace{2mm}  A chain homotopy, $I$, between two chain maps $f:(C_i,\,\delta_i)
\to (D_i,\,\partial_i)$ and 

$g:(C_i,\,\delta_i) \to (D_i,\,\partial_i)$ consists
of a sequence of maps $I_i:C_i \to D_{i+1}$ such that for each $i$, 
$\partial_{i+1} I_i +I_{i-1} \delta_i =f_i -g_i$.

\hspace{2mm} A chain homotopy equivalence,
$f:(C_i,\,\delta_i) \to (D_i,\,\partial_i)$, is a chain map for which there
exists a chain map $g:(D_i,\,\partial_i) \to (C_i,\,\delta_i)$ and a chain
homotopy, $I$, between the identity and $fg$, and a chain homotopy $J$, between
the identity and $gf$.  Two algebraic complexes, $(C_i,\,\delta_i)$ and $(D_i,\,\partial_i)$ are said to be homotopy equivalent precisely
when there exists a homotopy equivalence between them.  In this case, we may write $(C_i,\,\delta_i) \sim (D_i,\,\partial_i)$.

\hspace{2mm}  An important notational point is that whenever we write down an
algebraic chain complex, any maps denoted by dotted arrows or modules connected
to the complex by such arrows, are not part of the complex.

\bigskip
\hspace{2mm}  Let $A$ be a left module over $R$ and let $C_* = (C_i,\,\delta_i)$ be an algebraic complex (of right modules as usual).  Let 
$(C_i \otimes_R A,\,{\delta_i}_*)$ denote the tensor product of $(C_i,\,\delta_i)$ with $A$ over $R$.  Then $H_i(C_*; A)$ denotes the kernel of ${\delta_i}_*$ quotiented by 
the image
of ${\delta_{i+1}}_*$.

\hspace{2mm}  Let $A$ be a right module over $R$ and let $C_* = (C_i,\,\delta_i)$ be an algebraic complex (of right modules as usual).  Let 
$({\rm Hom}_R(C_i,  A),\,\delta_i^*)$ denote the functor ${\rm Hom}_R (\_, A)$ applied to $(C_i,\,\delta_i)$.  Then 
$H^i(C_*; A)$ denotes the kernel of $\delta_{i+1}^*$ quotiented by the image
of $\delta_{i}^*$.

\hspace{2mm} This notation for homology and cohomology is slightly non-standard.

\bigskip

\hspace{2mm}  We say a sequence of modules and maps is exact when the image of
each map is the kernel of the next.

\hspace{2mm} In this section, we work over a fixed ring $R$, with unit.  All maps, modules, algebraic chain complexes and chain homotopies  
will be assumed to be over $R$.  Modules denoted by ``$P_i$" and ``$Q_i$" may be assumed to be projective.  

\bigskip
\hspace{2mm}Suppose we have an algebraic complex 

\xymatrix{ \cdots \ar[r]^{\partial_{n+1}}&F_n \ar[r]^{\partial_{n}}&F_{n-1}
 \ar[r]^{\partial_{n-1}}&\cdots
}$\hfill (1)$

\bigskip
\hspace{2mm}  Replacing it with the complex 

\xymatrix{ \cdots \ar[r]^{\partial_{n+1}}&F_n \oplus F \ar[r]^{\partial_{n}
\oplus 1}&
F_{n-1} \oplus F
 \ar[r]^{\qquad\partial_{n-1}}&\cdots}$\hfill (2)$

is called a simple homotopy equivalence.  We may define a chain map $i$ from (1)
to (2) which is the identity on $F_r$, for $r \neq n, n-1$ and the natural
inclusion for $r=n, n-1$.  Similarly we may define a chain map $j$ from (2)
to (1) which is the identity on $F_r$, for $r \neq n, n-1$ and the natural
projection for $r=n, n-1$.  Then $ji$ is the identity and $ij$ differs from the
identity by $0$ on the $F_i$ and the identity on the copies of $F$.

\bigskip
\hspace{2mm}  Let $I_r:F_r \to F_{r+1}$ be $0$ for $r \neq n-1$ and let
$I_n:F_{n-1} \oplus F \to F_{n-1} \oplus F$ be $0$ on $F_{n-1}$ and the identity
on $F$.  Then we have $1-ij=\partial' I + I \partial'$ where $\partial'$ agrees
with $\partial_r$ for $r \neq n$, and is the map $\partial_n \oplus 1$ on $F_n
\oplus F$.  From this we see that a simple homotopy equivalence is indeed a
homotopy equivalence.

\bigskip
Let 

$$
P_n \stackrel{\partial_n}{\longrightarrow} \dots \stackrel{\partial_2}{\longrightarrow} P_1 \stackrel{\partial_1}
{\longrightarrow} P_0 \stackrel{\epsilon}{\dashrightarrow} M \dashrightarrow 0
$$

and 

$$
Q_n \stackrel{\partial_n'}{\longrightarrow} \dots \stackrel{\partial_2'}{\longrightarrow}Q_1 \stackrel{\partial_1'}
{\longrightarrow} Q_0  \stackrel{\epsilon'}{\dashrightarrow} M \dashrightarrow 0
$$

be exact sequences.

\bigskip
Let $R_0=P_0$ and $S_0=Q_0$.  Define $R_i, S_i$, for $i=1,\dots n$ by 

\bigskip
$R_i = S_{i-1} \oplus P_i$ 

$S_i = R_{i-1} \oplus Q_i$

\bigskip
\hspace{2mm} Note $R_i$ and $S_i$ are projective, as they are constructed as direct sums of projective modules.

\bigskip
\thm{The complexes}

{\it

$$
P_n \oplus S_n \stackrel{\partial_n \oplus 0}{\longrightarrow} P_{n-1} \stackrel{\partial_{n-1}}{\longrightarrow}
\dots \stackrel{\partial_2}{\longrightarrow} P_1 \stackrel{\partial_1}{\longrightarrow} P_0 
\eqno(1)
$$

and 

$$
Q_n \oplus R_n \stackrel{\partial_n'\oplus 0}{\longrightarrow} Q_{n-1} \stackrel{\partial_{n-1}'}{\longrightarrow}
\dots \stackrel{\partial_2'}{\longrightarrow}Q_1 \stackrel{\partial_1'}{\longrightarrow} Q_0 
\eqno(2)
$$

are chain homotopy equivalent.}

\bigskip
\hspace{2mm} Chain homotopies induce isomorphisms on homology groups.  In particular, a chain homotopy between (1) 
and (2) will induce an isomorphism between ${\rm Ker}(\partial_n \oplus 0)$ and ${\rm Ker}(\partial_n' \oplus 0)$, so 
we have the following corollary:
 
\bigskip
\cor{(Schanuel)$\rm{Ker} \partial_n \oplus S_n = \rm{Ker} \partial_n'
\oplus R_n$}

\bigskip
\hspace{2mm} Proof of theorem:  We perform a series of simple homotopy equivalences, $u_i$ on $(1)$, for
$i=1,\dots, n$.  Recall $R_0=P_0$.  Define $u_1$ to consist of replacing 

$$
\stackrel{\partial_{2}}{\longrightarrow} P_1 \stackrel{\partial_1}{\longrightarrow} R_{0}
\stackrel{\epsilon}{\dashrightarrow} M
$$

with

$$
\stackrel{\partial_{2}}{\longrightarrow} P_1 \oplus S_{0} \stackrel{\delta_1}{\longrightarrow} R_{0} \oplus S_{0} 
\stackrel{\epsilon \oplus 0}{\dashrightarrow} M
\eqno(3)
$$

where $\delta_1$ is defined by

$$\delta_1 = \left( \begin{array}{cc} \partial_1&0\\ 0&1 \end{array} \right)$$  

Note $R_1 \, = \, P_1 \oplus S_{0}$, so we can write $(3)$ as, 

$$
\stackrel{\partial_{2}}{\longrightarrow} R_1 \stackrel{\delta_1}{\longrightarrow} R_{0} \oplus S_{0} 
\stackrel{\epsilon \oplus 0}{\dashrightarrow}M
$$

Define $u_i$ to consist of replacing 

$$
\stackrel{\partial_{i+1}}{\longrightarrow} P_i \stackrel{\partial_i}{\longrightarrow} R_{i-1}
\stackrel{\delta_{i-1}}{\longrightarrow}
$$

with

$$
\stackrel{\partial_{i+1}}{\longrightarrow} P_i \oplus S_{i-1} \stackrel{\delta_i}{\longrightarrow} R_{i-1} \oplus S_{i-1} 
\stackrel{\delta_{i-1} \oplus 0}{\longrightarrow}
\eqno(4)
$$

where $\delta_i$ is defined by

$$\delta_i = \left( \begin{array}{cc} \partial_i&0\\ 0&1 \end{array} \right)$$  

Note $R_i \, = \, P_i \oplus S_{i-1}$, so we can write $(4)$ as, 

$$
\stackrel{\partial_{i+1}}{\longrightarrow} R_i \stackrel{\delta_i}{\longrightarrow} R_{i-1} \oplus S_{i-1} 
\stackrel{\delta_{i-1} \oplus 0}{\longrightarrow}
$$

Finally, note that $u_n$ replaces

$$
P_n \oplus S_n \stackrel{\partial_n \oplus 0}{\longrightarrow} R_{n-1} \stackrel{\delta_{n-1}}{\longrightarrow}
$$

with

$$
R_n \oplus S_n \stackrel{\delta_n \oplus 0}{\longrightarrow} R_{n-1} \oplus S_{n-1} \stackrel{\delta_{n-1} \oplus 0}
{\longrightarrow}
$$

Hence applying $u_1,\dots,u_n$ to $(1)$ gives

\bigskip
$$
R_n \oplus S_n \stackrel{\delta_n \oplus 0}{\longrightarrow} R_{n-1} \oplus S_{n-1} \stackrel{\delta_{n-1} \oplus 0}
{\longrightarrow} \dots \dots  \stackrel{\delta_2 \oplus 0}{\longrightarrow}
R_1 \oplus S_1 \stackrel{\delta_1 \oplus 0}{\longrightarrow} R_{0} \oplus S_0
\eqno(5)
$$

\bigskip
\bigskip
\hspace{2mm} Similarly we perform a series of simple homotopy equivalences $v_i$ on $(2)$, where $v_i$ consists of replacing   

$$
\stackrel{\partial_{i+1}'}{\longrightarrow} Q_i \stackrel{\partial_i'}{\longrightarrow} S_{i-1}
\stackrel{\delta_{i-1}'}{\longrightarrow}
$$

with

$$
\stackrel{\partial_{i+1}'}{\longrightarrow} Q_i \oplus R_{i-1} \stackrel{\delta_i'}{\longrightarrow} S_{i-1} \oplus R_{i-1} 
\stackrel{\delta_{i-1}' \oplus 0}{\longrightarrow}
\eqno(6)
$$

where $\delta_i'$ is defined by

$$\delta_i' = \left( \begin{array}{cc} \partial_i'&0\\ 0&1 \end{array} \right)$$  

Note $S_i \, = \, Q_i \oplus R_{i-1}$, so we can write $(6)$ as, 

$$
\stackrel{\partial_{i+1}'}{\longrightarrow} S_i \stackrel{\delta_i'}{\longrightarrow} S_{i-1} \oplus R_{i-1} 
\stackrel{\delta_{i-1}' \oplus 0}{\longrightarrow}
$$

Also $v_n$ replaces 

$$
Q_n \oplus R_n \stackrel{\partial_n' \oplus 0}{\longrightarrow} S_{n-1} \stackrel{\delta_{n-1}'}{\longrightarrow}
$$

with

$$
S_n \oplus R_n \stackrel{\delta_n' \oplus 0}{\longrightarrow} S_{n-1} \oplus R_{n-1} \stackrel{\delta_{n-1}' \oplus 0}
{\longrightarrow}
$$

\bigskip
\bigskip
Hence applying $v_1,\dots,v_n$ to $(2)$ gives

\bigskip
$$
S_n \oplus R_n \stackrel{\delta_n' \oplus 0}{\longrightarrow} S_{n-1} \oplus R_{n-1} \stackrel{\delta_{n-1}' \oplus 0}
{\longrightarrow} \dots \dots  \stackrel{\delta_2' \oplus 0}{\longrightarrow}
S_1 \oplus R_1 \stackrel{\delta_1' \oplus 0}{\longrightarrow} S_{0} \oplus R_0
\eqno(7)
$$

\bigskip
\bigskip

\hspace{2mm} Recall $R_i$, $S_i$ are projective.  As $(5)$
and $(7)$ are chain homotopy equivalent to $(1)$ and $(2)$ respectively, they are both exact.  Further, they extend to exact sequences

$$
R_n \oplus S_n \stackrel{\delta_n \oplus 0}{\longrightarrow} R_{n-1} \oplus S_{n-1} \stackrel{\delta_{n-1} \oplus 0}
{\longrightarrow} \dots \dots  \stackrel{\delta_2 \oplus 0}{\longrightarrow}
R_1 \oplus S_1 \stackrel{\delta_1 \oplus 0}{\longrightarrow} R_{0} \oplus S_0 \stackrel{\epsilon \oplus 0}
{\dashrightarrow} M
\dashrightarrow 0
\eqno(8)
$$

and 

$$
S_n \oplus R_n \stackrel{\delta_n' \oplus 0}{\longrightarrow} S_{n-1} \oplus R_{n-1} \stackrel{\delta_{n-1}' \oplus 0}
{\longrightarrow} \dots \dots  \stackrel{\delta_2' \oplus 0}{\longrightarrow}
S_1 \oplus R_1 \stackrel{\delta_1' \oplus 0}{\longrightarrow} S_{0} \oplus R_0 \stackrel{\epsilon' \oplus 0}
{\dashrightarrow} M
\dashrightarrow 0
\eqno(9)
$$

\hspace{2mm} We complete the proof of theorem 1.1.1 by constructing a pair of inverse chain isomorphisms, $h$, $k$, between $(8)$
and  $(9)$.

\bigskip
\hspace{2mm}  As $R_0$, $S_0$ are projective, we may pick $f_0$, $g_0$ so that the following diagrams commute:

\setcounter{equation}{9}
\begin{eqnarray*}
R_0 \stackrel{\epsilon}{\longrightarrow} M\,\,   \qquad \qquad       R_0 \stackrel{\epsilon}{\longrightarrow} M\,\, 
\nonumber\\ \downarrow f_0 \,\, \quad \downarrow 1             \qquad \qquad       \uparrow g_0 \,\, \quad \uparrow 1 
\nonumber\\ 
S_0 \stackrel{\epsilon'}{\longrightarrow} M\,\,  \qquad \qquad       S_0 \stackrel{\epsilon'}{\longrightarrow} M\,\,  
\end{eqnarray*}
\hfill (10)

Define $h_0:R_0 \oplus S_0 \to S_0 \oplus R_0$ and $k_0:S_0 \oplus R_0 \to R_0 \oplus S_0$ by

$$
h_0= \left(\begin{array}{cc} f_0&1-f_0g_0\\ 1& -g_0 \end{array}\right) \qquad \qquad 
k_0= \left(\begin{array}{cc} g_0&1-g_0f_0\\ 1& -f_0 \end{array}\right)
$$

Direct calculation shows that $h_0k_0 = 1$ and $k_0h_0 = 1$.

\bigskip
Also from commutativity of (10), we deduce

$$
(\epsilon' \quad 0)\left(\begin{array}{cc} f_0&1-f_0g_0\\ 1& -g_0 \end{array}\right) = (\epsilon'f_0 \quad
\epsilon'(1-f_0g_0))=(\epsilon \quad 0) 
$$

and

$$
(\epsilon \quad 0)\left(\begin{array}{cc} g_0&1-g_0f_0\\ 1& -f_0 \end{array}\right) = (\epsilon g_0 \quad
\epsilon(1-g_0f_0))=(\epsilon' \quad 0) 
$$

\bigskip
Hence the following diagrams commute:

\begin{eqnarray*}
R_0 \oplus S_0 \stackrel{\epsilon \oplus 0}{\longrightarrow} M\,\,   
\qquad \qquad       R_0 \oplus S_0 \stackrel{\epsilon \oplus 0}{\longrightarrow} M\,\,\\
\downarrow h_0 \,\,\quad\,\, \quad \downarrow 1         \quad \qquad \qquad       \uparrow k_0 \,\, \quad\,\, \quad 
\uparrow 1\\ 
S_0 \oplus R_0 \stackrel{\epsilon'\oplus 0}{\longrightarrow} M\,\,  \qquad \qquad       S_0 \oplus R_0 \stackrel{\epsilon' 
\oplus
0}{\longrightarrow} M\,\,\\   
\end{eqnarray*}

\hspace{2mm} Now suppose that for some $i \leq n$, we have defined $h_j:R_j \oplus S_j \to S_j \oplus R_j$ and $k_j:S_j \oplus R_j 
\to R_j \oplus S_j$ for $j=0,\dots,i-1$, so that for each $j$, we have $h_jk_j=1$ and $k_jh_j=1$.  We proceed by induction.

\bigskip
As before, pick $f_i$, $g_i$ so that the following diagrams commute:

\begin{eqnarray*}
R_i \stackrel{\delta_i}{\longrightarrow} R_{i-1} \oplus S_{i-1}\,\,   
\qquad \qquad    R_i \stackrel{\delta_i}{\longrightarrow} R_{i-1} \oplus S_{i-1}\,\, \nonumber\\  
\downarrow f_i \,\quad \quad \, \quad \downarrow h_{i-1} \,\,\,\qquad \qquad \uparrow g_i \,\, \quad \quad \quad 
\uparrow k_{i-1} \,\,\,\,\nonumber\\ 
S_i \stackrel{\delta_{i}'}{\longrightarrow} S_{i-1} \oplus R_{i-1} \,\,  
\qquad \qquad   S_i \stackrel{\delta_{i}'}{\longrightarrow} S_{i-1} \oplus R_{i-1} \,\,  
\end{eqnarray*}
\hfill(11)

Define $h_i:R_i \oplus S_i \to S_i \oplus R_i$ and $k_i:S_i \oplus R_i \to R_i \oplus S_i$ by

$$
h_i= \left(\begin{array}{cc} f_i&1-f_ig_i\\ 1& -g_i \end{array}\right) \qquad \qquad 
k_i= \left(\begin{array}{cc} g_i&1-g_if_i\\ 1& -f_i \end{array}\right)
$$

Direct calculation shows that $h_ik_i = 1$ and $k_ih_i = 1$.

\bigskip
Recall $h_{i-1}k_{i-1}=1$ and $k_{i-1}h_{i-1}=1$.  From commutativity of (11) we deduce

$$
(\delta_i' \quad 0)\left(\begin{array}{cc} f_i&1-f_ig_i\\ 1& -g_i \end{array}\right) = (\delta_i'f_i \quad
\delta_i'(1-f_ig_i))=h_{i-1}(\delta_i \quad 0) 
$$

and

$$
(\delta_i \quad 0)\left(\begin{array}{cc} g_i&1-g_if_i\\ 1& -f_i \end{array}\right) = (\delta_i g_i \quad
\delta_i(1-g_if_i))=k_{i-1}(\delta_i' \quad 0) 
$$

Hence the following diagrams commute:

\begin{eqnarray*}
R_i \oplus S_i \stackrel{\delta_i \oplus 0}{\longrightarrow} R_{i-1} \oplus S_{i-1}\,\,   
\qquad \qquad    R_i \oplus S_i \stackrel{\delta_i \oplus 0}{\longrightarrow} R_{i-1} \oplus S_{i-1}\,\,\\  
\downarrow h_i \,\quad \quad \, \quad \quad \downarrow h_{i-1} \,\,\,\qquad \quad \qquad \uparrow k_i \,\, \quad \quad \quad 
\quad 
\uparrow k_{i-1}\,\,\,\,\\ 
S_i \oplus R_i \stackrel{\delta_{i}' \oplus 0}{\longrightarrow} S_{i-1} \oplus R_{i-1} \,\,  
\qquad \qquad   S_i \oplus R_i\stackrel{\delta_{i}' \oplus 0}{\longrightarrow} S_{i-1} \oplus R_{i-1} \,\,  
\end{eqnarray*}

\hspace{2mm} Together with the identity on $M$, the $h_i$, $k_i$ are therefore a pair of mutually inverse chain maps, 
between $(8)$ and $(9)$.  Hence $(1)$ is chain chain homotopy equivalent to $(5)$, which is chain isomorphic to $(7)$ which in turn is 
chain chain homotopy equivalent to $(2)$.
\hfill $\Box$ \,\,

\bigskip
\hspace{2mm}If the modules $P_i$ and $Q_i$ are finitely generated and free, then
$R_n$ and $S_n$ are also finitely generated and free.  Let Free${}_n(M)$ denote 
the set of homotopy types of f.g. free $n$- term resolutions, of a module $M$.  We may 
give Free${}_n(M)$ the structure of a tree, by placing an edge between the homotopy
types of any pair of resolutions of the form

$$
F_n
\stackrel{\partial_{n}}{\longrightarrow} F_{n-1}
\stackrel{\partial_{2}}{\longrightarrow} \cdots
\stackrel{\partial_{2}}{\longrightarrow} F_1
\stackrel{\partial_1}{\longrightarrow} F_{0}
\stackrel{\epsilon}{\dashrightarrow} M
$$

and 

$$
F_n \oplus R
\stackrel{\partial_{n} \oplus 0}{\longrightarrow} F_{n-1}
\stackrel{\partial_{2}}{\longrightarrow} \cdots
\stackrel{\partial_{2}}{\longrightarrow} F_1
\stackrel{\partial_1}{\longrightarrow} F_{0}
\stackrel{\epsilon}{\dashrightarrow} M
$$

where $R$ is regarded as a module over itself.  From theorem 1.1.1, we may
conclude

\thm{Free${}_n(M)$ is a connected tree.}

\bigskip
\hspace{2mm}Another tree of interest is the stable tree of a module.

\define{Stable equivalence}\hspace{2mm}  Two modules over $R$, $L$ and $N$, are  
stably equivalent if there is an isomorphism $L\oplus R^a \to N \oplus R^b$, for
integers $a$ and $b$.

\hspace{2mm}The stable class of a module $K$ is the set of modules stably
equivalent to $K$.  We may give this the structure of a tree, by assigning an
edge to any pair of modules $L$, $N$, where $N$ is isomorphic to $L \oplus R$.

\hspace{2mm} From corollary 1.1.2, we know that if $(A_i, \delta_i)$ and $(B_i,
\partial_i)$ are elements of Free${}_n(M)$, then ker$(\delta_n)$ and
ker$(\partial_n)$ are in the same stable class.  We denote this class
$\Omega_{n+1}(M)$.  

\hspace{2mm} We may conclude that we have a map of trees Free${}_n(M)$ $\to$
$\Omega_{n+1}(M)$, which sends an element of Free${}_n(M)$, $(A_i, \delta_i)$,
to ker$(\delta_n)$.

\sec{The Derived Category}


\hspace{2mm}The previous section provides some motivation for considering a
category where stably equivalent modules are isomorphic objects.  We follow
Johnson[4], in constructing such a category.  All the results of this section
are explained in greater detail in [4], \S19.  We provide a summary, for
narrative purposes.

\hspace{2mm}  Modules over the group ring, $\ZG$, of a group form the objects of a
category, whose morphisms are $\ZG$- linear maps.  

\bigskip

\hspace{2mm}
 A map $f: X \to Y$ is said to factor through a projective if there exists a 
 projective module $P$ and maps $a$, $b$ such the following diagram commutes:

\xymatrix{
X\ar[d]_a \ar[r]^f &Y \\
P \ar[ur]^b}

\lem {A map factors through a projective module if and only if it factors 
through a free  module.}

\proof{Suppose $f =  ba$ as before.  We have some projective module $Q$ 
such that $P \oplus Q$ is free.  Let $i:P \hookrightarrow P \oplus Q$ be the
natural inclusion and $p: P \oplus Q \twoheadrightarrow  P$ be the natural
projection.  Then $pi=1$.  Hence $f=ba =(bp)(ia)$}

\hspace{2mm}  For any pair of modules, $M$ and $N$, define an equivalence 
relation $\sim$ on ${\rm Hom}_{\ZG} (M, N)$ given by $f \sim g$ if and only if $f -
g$
 factors through a projective module.   To see this is transitive, note that if 
$f - g$ factors through $P$ and $g - h$ factors through $Q$ then both $f-g$ and
$g-h$ factor through $P \oplus Q$ so
$f - h = (f - g) + (g - h)$ factors through $P\oplus Q$.

\bigskip
\hspace{2mm}
If $f$ factors through a projective then clearly any map composed 
with $f$ also does.  By linearity of 
composition we therefore have that composition of maps under $\sim$ is 
well defined.

\bigskip
\hspace{2mm}Hence we can define the derived category of the 
category of $\ZG$- modules, Der$(\ZG)$, to be the category whose objects 
are (left) $\ZG$- modules and whose morphisms are $\ZG$-linear homomorphisms 
under the
relation $\sim$.  

\lem {In the derived category, $0 \cong P$ for all $P$ projective.}

\proof{
We have unique maps $i:0 \to P$ and $p:P \to 0$. $pi$ is the identity on $0$ and
$ip-1$ factors through $1:P \to P$.}





\lem{In the derived category $M \cong M\oplus P$ for $P$ projective.}

\proof{Let $i:M \to M \oplus P$, $p:M \oplus P \to P$ be the natural 
inclusion, and natural projection
 between $M$ and  $M \oplus P$ respectively.  Clearly $pi = 1$. 
 Also $1-ip$ restricts to $0$ on the summand $M$ and $1:P \to P$, so $1-ip$ 
 factors through projection onto $P$.}

\prop{Suppose we have an exact sequence of modules over $\ZG$, $G$ a finite
group:
\newline
\bigskip
\xymatrix{
0\ar[r]&A\ar[r] &F_n \ar[r]&F_{n-1} \ar[r]&\cdots\ar[r]&F_0 \ar[r]&M\ar[r]&0 
}
where the $F_i$ are free.  Then ${\rm End}_{\rm Der}(A) \cong {\rm End}_{\rm
Der}(M)$. }

\proof{See \cite{John1}.  The isomorphism is given by taking taking any endomorphism
$f:M \to M$ and using projectivity of the $F_i$ to extend to a commutative
diagram such as 
\newline
\bigskip
\xymatrix{
0\ar[r]&A\ar[d]^{D(f)}\ar[r] &F_n\ar[d] \ar[r]&F_{n-1}\ar[d] 
\ar[r]&\cdots\ar[d]\ar[d]\ar[r]&F_0\ar[d] \ar[r]&M\ar[d]^f\ar[r]&0 \\
0\ar[r]&A\ar[r] &F_n \ar[r]&F_{n-1} \ar[r]&\cdots\ar[r]&F_0 \ar[r]&M\ar[r]&0
}
\bigskip
\newline
${}$\hspace{2mm} Then if $f'$ is equivalent in the derived category to $f$, any choice of $D(f)$
will be equivalent in the derived category to any choice of $D(f')$.
\newline
${}$\hspace{2mm} The inverse is constructed dually, using the relative injectivity of the $F_i$. 
See \cite{John1} for details of the diagram chases.}

\define{Algebraic $n$- complex} An algebraic $n$- complex is an algebraic
complex $(F_i, \partial_i)$, $i=1,\cdots,n$, over $\ZG$, satisfying:

\bigskip
i)$F_i$ is free and  finitely generated, $i=1,\cdots,n$.

\bigskip
ii)coker$(\partial_1)=\Z$.

\bigskip
iii)ker$(\partial_1)= $im$(\partial_2)$.

\bigskip \hspace{2mm}
Let $G$ be a finite group of order $n$ and let $J=H_2(X;\ZG)$ for some algebraic 2-complex, $X$, over 
$\ZG$.

\define{k- invariant}  We define a map $k:{\rm End}(J) \to \Z_n$.  Given any map $\alpha:J \to J$, one may extend it to a chain map, 
$X \rightarrow X$.  
The chain map will induce a map on the last cokernel, $\Z$.  This map will be multiplication by $m$, for 
some integer $m$.  A map $\Z \to \Z$ factors through a projective module if and only if it factors through a free module, hence if and 
only if it is some
sum of maps which factor through $\ZG$.  As any such map will be multiplication by a number divisible by $n$, the number $m$ is 
determined up to
congruence modulo $n$.  The k- invariant of $\alpha$, $k(\alpha)$, is then defined to be the congruence of $m$ in $\Z_n$.

\bigskip
\hspace{2mm}If $\alpha$ is an automorphism, $J \to J$, then $k(\alpha)$ will be a unit, modulo $n$. 

\define{Swan Map} The Swan map, ${\rm Aut}(J) \to \Z_n^*$, \hspace{1mm} sends an automorphism, $\alpha$, to its k- invariant, $k(\alpha)$.

\hspace{2mm}Note that $\Z_n$ is being identified here with End${}_{\rm Der}(J)$.  Hence if the Swan map was defined with respect to a
different algebraic $2$- complex, it would be the same map, as the only ring isomorphism $\Z_n \to \Z_n$ is the identity.

\thm{If the Swan map, ${\rm Aut}(J) \to \Z_n^*$, is surjective, then any algebraic $2$- complex, $Y$, with $H_2(Y;\ZG)=J$ satisfies $X \sim
Y$.}

\hspace{2mm} Proof:  Any element of Aut$(J)$ will induce a chain map 
$X$ to $Y$.  Let $k$ be the multiplication induced on $\Z$.  Now suppose that the Swan map is surjective, restricted to a map from 
${\rm Aut}_{\ZG}(J) \to {\rm Aut}_{\rm Der}(J)$.  Then 
we have a chain map
$X\to X$ with k-invariant the inverse of $k$ modulo $n$, and which induces an isomorphism on $J$. Composing this chain 
map with the chain map $X \to Y$ gives a chain map which induces
an automorphism on $J$ and has k-invariant 1.  We can replace this with a chain map that induces an automorphism on $J$ and actually induces
multiplication by 1 on $\Z$.  This will necessarily be a homotopy equivalence. 
Hence $X$ and $Y$ will be chain homotopy equivalent.  See \cite{John} for details
of the diagram chases.
\hfill $\Box$

\bigskip
\hspace{2mm} If the Swan map is surjective for a module $J$, then it is also surjective for $J \oplus \ZG$, as if $\alpha \in {\rm Aut} (J)$ then we have an
automorphism $J \oplus \ZG \to J \oplus \ZG$ with the same $k$- invariant, given
by the matrix:

$$
\left(\begin{array}{cc}  
\alpha&0 \\
0 &1 
\end{array} \right)
$$

\define{Minimal module} A module is minimal in its stable tree if it does not contain a summand isomorphic to $\ZG$.

\hspace{2mm}  We may conclude that it is sufficient to check that this property holds for all minimal modules in the stable class of
$J$, in order to deduce that it holds for all modules in that stable class. 
Again see \cite{John} for details of this method.

\bigskip
\bigskip
\bigskip
\bigskip
\bigskip
\sec{Group Cohomology}


\hspace{2mm}  A resolution for a module $M$ over a ring is an exact sequence of 
modules, $E_i, \partial_i$, such that the cokernel of $\partial_1$ is $M$.

\hspace{2mm}  When $\Z$ is regarded as a module over $\ZG$ for some group $G$,
we will assume the trivial action.  Let $G$ be a finite group.  There exists 
a resolution of finitely generated free modules, $F_i$, over $\ZG$, for $\Z$
(see \cite{Macl}):

\bigskip
\xymatrix{\cdots\ar[r]^{\partial_3}&F_2\ar[r]^{\partial_2}&F_1\ar[r]^{\partial_1}&F_0 \ar@{.>>}[r]|-{} &\Z}

\define{ For a left $\ZG$ module $A$, $H_n(G; A)$ is given by tensoring the 
resolution for $\Z$ with $A$ and
taking the kernel of $\partial_n$ quotiented out by the image of
$\partial_{n+1}$}

\bigskip
\define{ For a right $\ZG$ module $A$, $H^n(G; A)$ is given by applying 
${\rm Hom}_{\ZG}(\bullet,A)$ to the resolution for $\Z$ with $A$ and
taking the kernel of $\partial_{n+1}^*$ quotiented out by the image of
$\partial_{n}^*$}

\hspace{2mm} We note the following:

\prop{ Let $A$ be a module over $\ZG$.  The elements of $H^2(G;A)$
parametrize short exact sequences of groups of the form 
\newline
$1 \to A \to E \to G \to 1$, 
\newline where
the conjugation action of $G$ on $A$ is  
given by the $ZG$- action on $A$.}

(See \cite{Thom})

\bigskip
\hspace{2mm}  We may generalize cohomology groups to general rings:

\define{Given modules $M$ and $A$, ${\rm Ext}^n(M, A)$ is the $n$'th cohomology
group of a resolution of $M$ with coefficient module $A$.}

\prop{${\rm Ext}^1(M, A)$ parametrizes short exact sequences of modules of the
form $A \to X \to M$, up to chain isomorphism with identities at both ends.}

(See \cite{Macl})

\hspace{2mm}  Combining the previous two propositions we see that there is a
correspondence between short exact sequences of the form

$$
1 \to A \to ? \to G \to 1$$

and short exact sequences of modules of the form

$$
0 \to A \to ? \to IG \to 0$$

where $IG$ is the kernel of the map augmentation map $\epsilon:\ZG \to \Z$,
sending $1 \in ZG$ to $1 \in \Z$.  The kernel of the augmentation map is
called the augmentation ideal.

\bigskip
\hspace{2mm}As a result of this, any module, $A$ which occurs in an exact
sequence of the form 

$$0 \to A \to F \to \ZG \stackrel{\epsilon}{\to} \Z \to 0$$

corresponds to some surjection $E \to G$, for some group $E$.

\hspace{2mm}  Specifically, given a short exact sequence of groups

\xymatrix{1 \ar[r]& A \ar[r]^i &E\ar[r]^j&G \ar[r]& 1}

we have a $\ZG$- module $(IE)_G$, which is the augmentation ideal of $E$,
quotiented out by the ideal generated by elements of the form $x(e_1-e_2)$ where
$x \in IE$ and $e_1, e_2 \in E$ both map to the same element of $G$, via $j$. 
The action of an element $g \in G$ on $(IE)_G$ is the action of any preimage of
$g$.  The choice of preimage does not effect the action.

\bigskip
\hspace{2mm}  The map $j$ induces a map from $(IE)_G \to IG$.  The kernel of
this map is isomorphic to $A$ with $G$- action given by conjugation, as before.

\bigskip
\hspace{2mm} Conversely, we may construct an inverse to this operation:  Given
any short exact sequence of modules

\xymatrix{0 \ar[r]& A \ar[r]^i &M\ar[r]^j&G \ar[r] &0}

we let $E_M$ denote the subset of $M$ which maps to an element of the form $g-1, g
\,\,\in G$, via $j$.  We define a product on $E_M$ by setting for each $e_1, e_2
\in E_M$, the product $e_1 \circ e_2$ is equal to $e_1g+e_2$, where $e_2$ maps
to $g-1$.  We have a map $E_M \to G$ given by sending $e_1$ to $g$ where $j$
sends $e_1$ to $g-1$.  The kernel of this map is $A$, so we have 
a short exact sequence of groups

\xymatrix{1 \ar[r] &A \ar[r]^i &E_M\ar[r]^j&G \ar[r]& 1}

\bigskip
\hspace{2mm}  For any group $G$, we may denote by $F_G$ the free group generated
by the underlying set of $G$.  There is a natural surjection $F_G
\twoheadrightarrow G$.  Let $K_G$ denote the kernel of this surjection, and
let $K_G'$ denote the commutator subgroup of $K_G$.  We have a short exact
sequence:

\xymatrix{1 \ar[r]& K_G/K_G' \ar[r] &F_G/K_G'\ar[r]&G \ar[r]& 1} $${}\eqno{ (1)}$$

\prop{This sequence splits on the right, if and only if $H^2(G;A)$ vanishes for
all coefficient modules $A$.}

\proof{$K_G/K_G'$ is a $\ZG$- module, with $G$ action given by conjugation as
before.  If $H^2(G;K_G/K_G')=0$ then by proposition 1.3.3 (1) must split on
the right.
\bigskip
\newline
${}$\hspace{2mm}  Conversely, suppose (1) splits on the right, with $k:G \to
F_G/K_G'$ as the splitting map and we have any short
exact sequence
\newline
\xymatrix{1 \ar[r]& A \ar[r]^i &E\ar[r]^j&G \ar[r]& 1} $${}\eqno{ (2)}$$
\newline
with $A$ abelian.  We may construct maps $f_2:F_G/K_G' \to E$ and 
$f_1:K_G/K_G' \to A$
to make the following diagram commute:
\newline
\xymatrix{K_G/K_G' \ar[r]\ar[d]^{f_1} &F_G/K_G'\ar[r]\ar[d]^{f_2}&G\ar[d]^1\\
 A \ar[r]^i &E\ar[r]^j&G}
\bigskip\newline${}$\hspace{2mm}
 This is done by sending each generator of $F_G$ to some element of $E$ which
 lies in the preimage of the corresponding element of $G$.  This is well defined
 as $K_G$ maps into $A$ and $A$ is abelian.
 \bigskip\newline${}$\hspace{2mm}
 We have $j \circ (f_2 \circ k)=(j \circ f_2) \circ k=1$.  Hence (2) splits. 
 But $A$, together with the conjugation action of $G$ on it,was chosen 
 arbitrarily, hence $H^2(G;A)=0$, for all coefficient modules $A$.
 } 


\sec{Periodicity}

\hspace{2mm} The main result of this section is well known to those working in the field, though an explicit proof is hard to find in
the literature.  The result states that the only groups which are (homologically) of period $2$ are cyclic.  Swan makes a 
comment outlining an argument in \cite{Swan}.  We give a more direct proof which is elementary and avoids the technical difficulties of Swan's argument.  
This result will be used in the next section.

\bigskip
\hspace{2mm}We say a finite group, $G$, has period 2 if and only if there exists an exact sequence of $\ZG$ modules 

$$
0 \rightarrow \Z \rightarrow P_1 \rightarrow P_0  \rightarrow \Z \rightarrow 0
$$ 

with $P_0$, $P_1$ finitely generated projective.  (We will assume the $G$ - action on $\Z$ to be trivial throughout).

\bigskip 
$\epsilon$ will denote the augmentation map $\ZG \to \Z$ which takes $1 \in \ZG$ to $1 \in \Z$.  $\epsilon^*$ will denote its dual, which sends $1 \in \Z$ to $\sum_{g \in G}g \in \ZG$.

\bigskip
\prop{Let $G$ be a finite group which has period 2.  Then there exists an exact sequence
\newline
$$
0 \longrightarrow \Z \stackrel {}{\longrightarrow} S \stackrel {}{\longrightarrow} \ZG  \stackrel {\epsilon}{\longrightarrow} \Z \longrightarrow 0
$$ 
\newline
where S is projective.}

\bigskip

\hspace{2mm} Proof:  If $G$ has period 2 then there exists an exact sequence

$$
0 \rightarrow \Z \rightarrow P_1 \rightarrow P_0  \rightarrow \Z \rightarrow 0
$$ 

with $P_0$, $P_1$ finitely generated and projective.  Choosing $Q$ such that $Q \oplus P_0 \cong \ZG^n$ for some $n$, and performing a simple
congruence, we get the exact sequence

$$
0 \rightarrow \Z \stackrel {a}{\rightarrow} P \stackrel {b}{\rightarrow} \ZG^n  \stackrel {c}{\rightarrow} \Z \rightarrow 0
\eqno (1)
$$ 

where $P=P_1 \oplus Q$ and is finitely generated, projective.

\bigskip
\hspace{2mm}  Let a $\ZG$- basis be chosen, $<e_1,...,e_n>$, for $\ZG^n$.  As $c$ is a surjection, and $c(e_ig)=c(e_i)$ for all $g$ 
and $i$, some $\Z$- linear combination of the $c(e_i), \quad i=1,..,n,$ must equal $1$.  Let $E$ denote the $\Z$- linear
 span of the $e_i, \quad i=1,...,n$.  It follows that $c$ restricts to a surjection of Abelian groups $E \twoheadrightarrow \Z$.
 
\bigskip
\hspace {2mm} Clearly $\Z$ is a free Abelian group, so this surjection splits, as a map of Abelian groups.  Therefore a $\Z$- basis of $E$,
$f_i, \quad i=1,...,n$ may be chosen, such that $c(f_1)=1$ and $c(f_i)=0, \quad i=2,...,n$.

\bigskip
\lem {The $f_i, \quad i=1,...,n$ are a $\ZG$- basis for $\ZG^n$.}

\hspace{2mm}Proof of lemma:  Any element $x \in \ZG^n$ can be uniquely written as $\sum_{g \in G} x_gg$
with $x_g \in E$.  Each $x_g$ may be uniquely written as a $\Z$- linear combination of the $f_i$.
Therefore $x$ may be written uniquely as a $\ZG$- linear combination of the $f_i$.

\hfill $\Box$ $\,\,\,$

(Proof of Proposition continued)

\hspace{2mm} The module, $\ZG^n$ splits as the direct sum of the $\ZG$- linear span 
of $f_1$ and the $\ZG$- linear span of $f_i, \quad i=2,...,n$.  Hence (1) may be
rewritten as

$$
0 \longrightarrow \Z \stackrel {a}{\longrightarrow} P \stackrel {b_1 \oplus b_2}{\longrightarrow} \ZG^{n-1} \oplus \ZG  \stackrel {0, \,\, \epsilon}{\longrightarrow} \Z \longrightarrow 0
\eqno (2)
$$ 

\hspace{2mm} By exactness, $b_1$ must surject onto $\ZG^{n-1}$, hence $b_1$ splits.  (Note $\ZG^{n-1}$ projective).  So $P=S \oplus \ZG^{n-1}$ where $S$ is the kernel of $b_1$.  $S$ is a summand of $P$, hence projective.  

\bigskip
\hspace{2mm} By exactness, $a$ maps $\Z$ into $S$.  Also $b_2$ restricts to a map $\beta$ from $S$ into the kernel of $\epsilon$ in $\ZG$. 
The following sequence is obtained:

$$
0 \longrightarrow \Z \stackrel {a}{\longrightarrow} S \stackrel {\beta}{\longrightarrow} \ZG  \stackrel {\epsilon}{\longrightarrow} \Z \longrightarrow 0
\eqno (3)
$$ 

\lem{This sequence is exact.}

\hspace{2mm} Proof of lemma:  The image of $\beta$ must equal the kernel of $\epsilon$, as for any $x \in
ker (\epsilon)$, $b_1 \oplus b_2$ must map some element of $P$ to $(0,x) \in \ZG^{n-1} \oplus \ZG$.  This element
of $P$ must lie in $S$.  (Note $S$ was defined as the kernel of $b_1$).  

\hspace{2mm} The kernel of $\beta$ in $S$ is by definition the intersection of the kernels of $b_1$ and $b_2$ in
$P$. This equals the kernel of $b_1 \oplus b_2$ in P, and hence is equal to the image of $a$ in $S$.  Therefore (3) is exact.

\hfill $\Box$ $\,\,\,$

This completes the proof of the proposition.

\hfill $\Box$ $\,\,\,$

\bigskip
\lem{All modules and cokernels in (3) are torsion free and of finite $\Z$- rank}. 

\hspace{2mm} Proof:  $\Z$ and $\ZG$ are torsion free and of finite $\Z$- rank.  $S$ is finitely generated and
projective, hence torsion free and of finite $\Z$- rank.  By exactness, the cokernels are submodules of torsion
free modules of finite $\Z$- rank, hence themselves torsion free and of finite $\Z$- rank.  

\hfill $\Box$ $\,\,\,$

\cor{The sequence (3) may be dualized to get the exact sequence
\newline 
$$
0 \longrightarrow \Z \stackrel {\epsilon^*}{\longrightarrow} \ZG \stackrel {\beta^*}{\longrightarrow} T  \stackrel {a^*}{\longrightarrow} \Z \longrightarrow 0 
\eqno (4) 
$$ 
\newline
where $T$ is the dual of $S$.}

\bigskip
T is projective, as dualizing commutes with taking direct sums.

\bigskip
\lem{There is no surjective homomorphism from $T$ to $\Z \oplus \Z$.}

\hspace{2mm} Proof: Tensoring (4) with $\Q$, the rationals, yields 

$$
0 \longrightarrow \Q \stackrel {}{\longrightarrow} \QG \stackrel {}{\longrightarrow} T \otimes \Q \stackrel {}{\longrightarrow}  \Q \longrightarrow 0 
$$ 

\hspace{2mm} This is an exact sequence of finitely generated modules over $\QG$. As $\QG$ is semi-simple, the
'Whitehead trick' may be performed: $T \otimes \Q \oplus \Q \cong \QG \oplus \Q$.  

\bigskip
\hspace{2mm} Cancellation gives:  $T \otimes \Q \cong \QG$.  

\bigskip
\hspace{2mm} Suppose $f$ was a surjection $T \twoheadrightarrow \Z \oplus \Z$. Tensoring $f$ with $\Q$ would yield a
$\QG$- linear surjection $\QG \twoheadrightarrow \Q \oplus \Q$.  This is impossible because surjections
between finitely generated $\QG$ modules split and $\QG$ does not contain a copy of $\Q \oplus \Q$.

\hfill $\Box$ $\,\,\,$

\bigskip
\hspace{2mm} For any $\ZG$ module, $M$, let $M_G$ denote the module resulting from quotienting $M$ by
the submodule generated by elements of the form 

$m(1-g), \,\, m \in M, \,\, g \in G$.

\bigskip
\hspace{2mm}  Let $p: T \twoheadrightarrow T_G$ denote the natural surjection.  For some $\ZG$ module $W$, $T \oplus W \cong \ZG^l$, for some
integer $l$.  So $T_G \oplus W_G \cong (T \oplus W)_G \cong (\ZG^l)_G \cong \Z^l$.  

\hspace{2mm}  Therefore, $T_G \cong \Z^j$ for some $j \leq l$.  But as $p$ surjects onto it, $T_G$ must equal
$\Z$ or $0$, by lemma 1.4.6.

\bigskip
\hspace{2mm}  Any map from $T$ to a $\ZG$ module with trivial $G$-action must factor through $p$.  So $a^*$ in (4) factors through $p$.

\begin {eqnarray*}
T_G \quad \quad \,\,\,&\\
\stackrel {p} {\nearrow} \quad \downarrow q\quad \quad&
\\
0 \longrightarrow \Z \stackrel {\epsilon^*}{\longrightarrow} \ZG \stackrel {\beta^*}{\longrightarrow} T  \stackrel {a^*}{\longrightarrow} \Z \longrightarrow 0
\end {eqnarray*}

$a^*$ is a surjection, so $T_G \neq 0$, hence $T_G \cong \Z$.  $q$ is a surjection $\Z \to \Z$.  Hence
$q$ must be an isomorphism and $ker(\alpha^*)=ker(p)$.  Consequently, the following sequence is exact

$$
0 \longrightarrow \Z \stackrel {\epsilon^*}{\longrightarrow} \ZG \stackrel {\beta^*}{\longrightarrow} T  \stackrel {p}{\longrightarrow} \Z \longrightarrow 0
\eqno (5)
$$

\bigskip
\lem{Let $G$ be a finite group, of order k, having period 2.  For any $\ZG$- module A, with trivial G-action, $H^1(G;A)= ker(\times k:A \rightarrow A)$, the kernel of multiplication by $k$.}

\bigskip
\hspace{2mm} Proof:  Combining (3) and (5) over $\Z$ gives the first few terms of a resolution for $\Z$ over $\ZG$:

\begin {eqnarray*}
\dots \dots \dots \longrightarrow S \stackrel {\beta} {\longrightarrow} \ZG  \stackrel {\epsilon^* \circ \epsilon} {\longrightarrow} \ZG \stackrel {\beta^*} {\longrightarrow} T \stackrel {p} {\twoheadrightarrow} \Z \\ 
\epsilon \searrow \quad \nearrow \epsilon^* \quad \qquad \qquad \,\,\,&\\
\Z \qquad \qquad \qquad \quad \,\,\,\,&
\end {eqnarray*}
$$
\eqno(6)
$$

\hspace{2mm} A map from $\ZG$ to $A$ is determined by the element of $A$ to which the identity is sent.  Hence ${\rm Hom}_{ZG}(\ZG,A)$ may be
identified with $A$.

\hspace{2mm} $\epsilon^* \circ \epsilon$ sends $1 \in \ZG$ to $\sum_{g \in G}g \in \ZG$.  So if $f:\ZG \to A$
sends $1 \in \ZG$ to $a \in A$, then  
$$
f\epsilon^* \circ \epsilon (1)=f(\sum_{g \in G}g)= f(1) \sum_{g \in G}g  = ak
$$ 

Suppose $f$ is a map from $T$ to A.  As the $G$- action on $A$ is trivial, $f$ must factor through $p$, hence $f \beta^* = 0$.

\begin {eqnarray*}
\dots\dots\dots\longrightarrow S \stackrel {\beta} {\longrightarrow} \ZG  \stackrel {\epsilon^* \circ \epsilon} {\longrightarrow} \ZG \stackrel {\beta^*} {\longrightarrow} T \stackrel {p} {\twoheadrightarrow} \Z \\ 
f \downarrow \quad \swarrow \,\,\,\,  &\\
A \qquad \,\,\,\,&
\end {eqnarray*}

Hence the cochain obtained by applying ${\rm {\rm Hom}}(\bullet , A)$ to (6) begins

$$
\dots\dots\dots \stackrel {}{\longleftarrow} A \stackrel {\times k}{\longleftarrow} A \stackrel {0} {\longleftarrow} {\rm Hom}_{ZG}(T,A)
$$

where $k$ is the order of $G$.  From this we see that $H^1 (G; A)$ is the kernel of multiplication by $k$ on $A$.

\hfill $\Box$ $\,\,\,$

\hspace{2mm}  The following is a standard result:

\prop{Let $G$ be a finite group and $A$ an Abelian group.  Regarding $A$ as a $\ZG$-module with
trivial $G$-action $H^1 (G; A) \cong {\rm Hom}(G,A)$}.

\prop{Let $G$ be a finite group and $A$ an Abelian group.  Let $G'$ denote the commutator
subgroup of $G$.  Then ${\rm Hom} (G, A) \cong {\rm Hom} (G/G', A)$.} 

\hspace{2mm}  Proof:   Given a homomorphism from $G$ to $A$, $G'$ must lie in its kernel.  Hence any homomorphism
from $G$ to $A$ factors uniquely through a homomorphism from $G/G'$ to $A$.

\hfill $\Box$ $\,\,\,$

\bigskip
\hspace{2mm}  Let $\bf{C}^*$ denote the Abelian group, consisting of unit complex numbers, with addition given by complex multiplication.  

\bigskip
\prop{If $B$ is a finite Abelian group then ${\rm Hom} (B, \C^*) \cong B$}

\hspace {2mm}  Proof:  Let $\Z_r$ denote the cyclic group of order $r$ and let $b$ be a generator of it. 
A map $\Z_r \to \bf{C}^*$ is determined by which $r^{th}$ root of unity $b$ is mapped to.  So ${\rm Hom} (\Z_r,
\bf{C}^*)$ is generated by any map sending $b$ to a primitive $r^{th}$ root of unity.  Hence ${\rm Hom} (\Z_r,
\bf{C}^*) \cong \Z_r$ 

\hspace {2mm} Now any finite Abelian group $B$ is the product of finite cyclic groups:
$
B \cong \bigoplus_{i=1}^m \Z_{r_{i}}
$

So 

$$
{\rm Hom} (B, \C^*) \cong {\rm Hom} (\bigoplus_{i=1}^m \Z_{r_{i}},\C ^*) \cong \bigoplus_{i=1}^m {\rm Hom} (\Z_{r_{i}}, \C^*) \cong \bigoplus_{i=1}^m \Z_{r_{i}}
\cong B
$$ 
\hfill $\Box$ $\,\,\,$

\bigskip
\hspace{2mm} $\C^*$ can be regarded as a $\ZG$- module, assuming the trivial action of $G$ on $\C^*$.  Note that the kernel of multiplication by an integer, $k$ on $\C^*$, is the cyclic group of $k^{th}$ roots of unity.

\bigskip 
\thm{Let $G$ be a finite group having period 2.  Then $G$ is cyclic.}

\hspace {2mm}  Proof: Given a finite group, $G$, which has period 2, let $k$ denote the order of $G$ and let $G'$ denote its commutator subgroup.

$$
G/G' \cong {\rm Hom}(G/G', \C ^*)  \cong {\rm Hom}(G,\C ^*) \cong H^1(G, \C ^*) \cong \Z_k
$$ \hfill (by propositions 1.4.10, 1.4.9, 1.4.8, and lemma 1.4.7)

\bigskip
As $\Z_k$ has the same order as $G$, 

$$
G \cong G/G' \cong \Z_k
$$

\hfill $\Box$ $\,\,\,$

The converse is also true:

\bigskip
\prop{If a finite group G is cyclic, then it has period 2.}

\hspace{2mm}  Proof: Let $t$ be a generator of $G$ and let $n$ be the order of $t$.  $IG$ is the submodule of $\ZG$ consisting of elements of the form $\sum_{i=0}^{n-1} t^i \lambda_i$ with $\lambda_i \in \Z, \sum_{i=0}^{n-1} \lambda_i=0$. It is generated by the element $t-1$.  

\bigskip
Let $\sigma$ denote $\sum_{i=0}^{n-1} t^i$. 

\bigskip
$(t-1) \sigma = 0$.  If $(t-1) \mu =0$ where $\mu= \sum_{i=0}^{n-1} \mu_i t^i$, then equating coefficients
gives $\mu_i =\mu_{i+1}$.  Hence $\sigma$ divides $\mu$.  Therefore $IG \cong \ZG / \sigma \ZG$.  

\bigskip
There exist exact sequences

$$
0 \longrightarrow IG \longrightarrow \ZG \stackrel {\epsilon}{\longrightarrow} \Z \longrightarrow 0
\eqno (7)
$$

and

$$
0 \longrightarrow \Z \stackrel{\epsilon^*}{\longrightarrow} \ZG \longrightarrow IG^* \longrightarrow 0 
\eqno (8)
$$

where (8) is the dual of (7).

\hspace{2mm} As $\epsilon^*(1)=\sigma$, from (8) it is observed that $IG^*$ is isomorphic to $\ZG / \sigma \ZG$, which is isomorphic to $IG$.
 Combining (8) and (7) over the isomorphism $IG \equiv IG^*$ gives

\begin {eqnarray*}
0 \longrightarrow \Z \stackrel{\epsilon^*} {\longrightarrow}  \!\!\!    &\ZG \,\,\, \longrightarrow   \,\,\, \ZG & \!\!\!   \stackrel {\epsilon}{\longrightarrow} \Z \longrightarrow 0\\
&\searrow \quad \quad \quad \nearrow& \\
&IG^*  \cong  IG&
\end{eqnarray*}

This is exact

\hfill $\Box$ $\,\,\,$

So a finite group is cyclic if and only if it has period 2.

\bigskip
\hspace{2mm}  If $G$ is a finite group which is not cyclic, there can be no isomorphism between $IG$ and $IG^*$. In fact they cannot be stably equivalent, as if they were, we would have 
$IG \oplus \ZG^n \cong IG^* \oplus \ZG^m$, for some $n, m$.  It would then be possible to construct the following exact sequence:

\begin {eqnarray*}
0 \longrightarrow \Z \stackrel{\epsilon^* \oplus 0} {\longrightarrow}     &\ZG \bigoplus \ZG^m \longrightarrow  \ZG \bigoplus \ZG^n &    \stackrel {\epsilon,\,\,\, 0}{\longrightarrow} \Z \longrightarrow 0\\
&\searrow \quad \quad \quad \quad \quad \quad \nearrow \,\,& \\
&IG^* \bigoplus \ZG^m  \cong  IG \bigoplus \ZG^n&
\end{eqnarray*}

contradicting theorem 3.1.11.

\bigskip
\hspace{2mm} Note that no infinite group can have period 2 in the sense that

\bigskip
\prop{ Let $H$ be any infinite group.  Let $f:\Z \to P$ be a map for some projective $\ZH$- module, $P$.  Then $f$=0.}

\hspace{2mm} Proof:  Suppose we have a non- zero map $f:\Z \to P$.  A 
module $Q$ may be chosen, such that $P \oplus Q$ is free.  Then $f \oplus 0:\Z \to P \oplus Q$ may be composed with projection onto 
a copy of $\ZH$, to yield a non-zero map $g:\Z \to \ZH$.  Let $x = g(1)$.  By $\ZH$- linearity, we know that $xh=x$ for all $h \in H$. However, as
$x$ is non-zero, it is expressible as a finite sum $\sum_i g_i \lambda_i$, with $\lambda_i$ integers and $\lambda_j$ non-zero for some $j$.  
Then, choosing $h$ such that $g_j h \neq g_k$ for any $k$, we have $xh \neq x$, giving the desired contradiction.   
\hfill $\Box$

\prop{If a finite group $G$ acts freely, on a
circle, then it is cyclic.}

\bigskip
\hspace{2mm} Proof:  By lifting the cells of the quotient manifold to the
circle, we would get a cellular resolution of the circle:

$$
\xymatrix{
\Z \ar[r] & F_1 \ar[r]&F_0\ar[r]& \Z}
$$             

\bigskip
\hspace{2mm} Any fixed point free map, $f$, on a circle is homotopic to the antipodal
map, as for each $x \in S^1$, we can choose a unique shortest arc between $f(x)$
and $-x$.  $f$ must therefore preserve orientation.  

\bigskip
\hspace{2mm} The group action on the left hand copy of the integers is therefore
trivial.  $F_1$ and $F_0$ are clearly finitely generated and free, so $G$ has
period 2 and is cyclic.

\hfill $\Box$


\bigskip
\hspace{2mm}Finally, in this section, we apply theorem 1.4.11 to obtain a lemma which we use in the following section.

\bigskip
\hspace{2mm} Let $G$ be a finite group. 

\define{Quaternionic} We say a real representation of $G$ is quaternionic if its endomorphism ring is $\Hq$, the
quaternions.

\hspace{2mm} Let $V_1, \cdots,V_k$ be the irreducible real representation of $G$.  Given a $\ZG$- module, $M$, of finite rank, we have the 
decomposition: 

$$M \otimes \R=\bigoplus_{i=1}^k V_i^{n_i}$$

\define{Eichler} $M$ satisfies the Eichler condition precisely when $V_i$ quaternionic implies that $n_i \neq 1$, for $i=1, \cdots,k$.

\thm{(Swan-Jacobinski)  Let $M$ be a torsion free $\ZG$- module of finite $\Z$- rank.  Suppose $M \oplus \ZG$ satisfies the Eichler
 condition.  Let $L$ be a $\ZG$- module with $\rk(L) \geq \rk(M \oplus \ZG)$, and $L$ stably equivalent to $M$. Then we have 
 $L \cong M \oplus \ZG^r$ for some $r \geq 1$. (See \cite{John1} \S 15)}

\hspace{2mm} In particular, note that for any finite group, $\ZG^2 \otimes\R$ will contain more than one copy of any irreducible
module.  Hence $\ZG^2$ satisfies the Eichler condition and has the form $\ZG \oplus \ZG$.  Therefore any stably free module of $\ZG$ -rank greater than 1, must be free.
 
\hspace{2mm} Note also that if $G$ is cyclic, $G$ does not have any irreducible real representations with endomorphism ring $\Hq$
so $0 \oplus \ZG$ satisfies the Eichler condition and all finitely generated stably free modules are free.

\lem{Let $G$ be a finite group.  Suppose there exists an exact sequence over $\ZG$ 
\newline
$$
SF\stackrel{\partial_{1}} {\longrightarrow} \ZG^a \stackrel{} {\longrightarrow} \Z \longrightarrow 0
$$
\newline
where $SF$ is a finitely generated stably free module.  Then $SF$ is free.}

\hspace{2mm} Proof:  If the $\ZG$- rank of $SF$ is greater than 1, the result would follow from the Swan-Jacobinski Theorem.  Suppose 
the $\ZG$- rank of $SF$ is 1.   It is sufficient to prove that $G$ is cyclic as then 
all finitely generated stably free modules over it 
would be free.  We will show that $G$ has period 2, as then it must be cyclic, by theorem 1.4.11.

\hspace{2mm}  Exactness and consideration of $\Z$- rank imply that the 
kernel of $\partial_1$, $K$, must have
$\Z$- rank congruent to 1 modulo the order of the group.  The $\Z$- rank of $K$ must be less than that of $SF$.  Hence if $SF$
has $\ZG$- rank 1, then $K$, must have $\Z$- rank 1. By 
the 'Whitehead Trick', $K \otimes \Q \cong \Q$  as
$\QG$- modules. Hence the $G$- action on $K$ is trivial and we have an exact sequence:

\bigskip
$$
0 \longrightarrow \Z \longrightarrow SF\stackrel{\partial_{1}} {\longrightarrow} \ZG^a \stackrel{} {\longrightarrow} \Z \longrightarrow 0
$$

\hfill $\Box$

\lem{Let $G$ be a finite group.  Suppose there exists an exact sequence over $\ZG$ 
\newline
$$
SF\stackrel{\partial_{2}} {\longrightarrow} \ZG^b \stackrel{\partial_{1}} {\longrightarrow} \ZG^a \stackrel{} {\longrightarrow}
\Z \longrightarrow 0
$$
\newline
where $SF$ is a finitely generated stably free module.  Then $SF$ is free.}

\hspace{2mm} Proof:  If the $\ZG$- rank of $SF$ is greater than 1, the result would follow from the Swan-Jacobinski Theorem.  Suppose 
the $\ZG$- rank of $SF$ is 1.   It is sufficient to prove that $G$ is cyclic as then 
all finitely generated stably free modules over it 
would be free.  We will show that $G$ has period 2, as then it must be cyclic, by theorem 1.4.11.

\hspace{2mm}  Exactness and consideration of $\Z$- rank imply that the 
kernel of $\partial_1$, $K'$, must have
$\Z$- rank congruent to 1 modulo the order of the group.  The $\Z$- rank of $K'$ must be less than that of $SF$.  Hence if $SF$
has $\ZG$- rank 1, then $K'$, must have $\Z$- rank 1. By 
the 'Whitehead Trick', $K \otimes \Q \cong \Q$  as
$\QG$- modules. Hence the $G$- action on $K'$ is trivial and we have an exact sequence:

$$
0\longrightarrow \Z {\longrightarrow} \ZG^b \stackrel{\partial_{1}} {\longrightarrow} \ZG^a \stackrel{} {\longrightarrow}
\Z \longrightarrow 0
$$

\hfill $\Box$

\bigskip
\sec{Free Resolutions}

\hspace{2mm}  The main theorem of this section, theorem 1.5.1, essentially shows that the homotopy type of any resolution may 
be represented by a resolution with a prespecified  initial segment:

\newpage
 
\thm{ Let 
\newline
$$
F_m \stackrel{\partial_{m}} {\longrightarrow} F_{m-1} \stackrel{\partial_{m-1}} {\longrightarrow} \dots\dots\dots\dots\dots\dots\dots F_1 \stackrel{\partial_{1}}
{\longrightarrow} F_{0} {\dashrightarrow} M \dashrightarrow 0      
\eqno[1]
$$
\newline
and
\newline
$$
\quad\quad \quad G_n \stackrel{\partial_{n}'} {\longrightarrow} G_{n-1} \stackrel{\partial_{n-1}'} {\longrightarrow}\dots\dots G_1 \stackrel{\partial_{1}'}
{\longrightarrow} G_{0} {\dashrightarrow} M \dashrightarrow 0      
$$
\newline
be exact sequences, over a ring $R$ with $n<m$ and the $F_i$, $G_i$ finitely generated free modules.  Then there exists an exact
sequence over $R$:
\newline
$$
F_m \stackrel{\partial_{m}} {\longrightarrow} F_{m-1} \stackrel{\partial_{m-1}} {\longrightarrow} \dots\dots \stackrel{\partial_{n+3}}
{\longrightarrow} F_{n+2} {\longrightarrow} SF \longrightarrow G_n \stackrel{\partial_{n}'} {\longrightarrow}
\dots\dots\dots G_1 \stackrel{\partial_{1}'} {\longrightarrow} G_{0} {\dashrightarrow} M \dashrightarrow 0      
$$
\newline
which is chain homotopy equivalent to [1] with $SF$ finitely generated, stably free.}

\hspace{2mm}Proof:  By Schanuel's Lemma, there exist finitely generated free modules $K$, $L$, such that 

$$
F_n \oplus L \stackrel{\partial_{n} \oplus 0} {\longrightarrow} F_{n-1} \stackrel{\partial_{n-1}} {\longrightarrow} \dots\dots F_1 \stackrel{\partial_{1}}
{\longrightarrow} F_{0}     
\eqno[2]
$$

is chain homotopy equivalent to 

$$
G_n \oplus K \stackrel{\partial_{n}' \oplus 0} {\longrightarrow} G_{n-1} \stackrel{\partial_{n-1}'} {\longrightarrow} 
\dots\dots G_1 \stackrel{\partial_{1}'} {\longrightarrow} G_{0}      
\eqno[3]
$$

 Let $f$ denote a chain homotopy
equivalence from $[2]$ to $[3]$ and let $g$ denote a chain homotopy inverse to $f$, from $[3]$ to $[2]$.  Let $I$ denote a chain homotopy from
$1_{[2]}$ to $gf$ and let $J$ denote a chain homotopy from $1_{[3]}$ to $fg$.  So $I_{r-1}\partial_r + \partial_{r+1}I_r=g_rf_r-1_{[2]}$  and 
$J_{r-1}\partial_r' + \partial_{r+1}'J_r=f_r g_r-1_{[3]}$.

\bigskip
\hspace{2mm} $[1]$ is chain homotopy equivalent to

$$
F_m \stackrel{\partial_{m}} {\longrightarrow} F_{m-1} \stackrel{\partial_{m-1}} {\longrightarrow} \dots
\stackrel{\partial_{n+2}}{\longrightarrow} F_{n+1} \oplus L \stackrel{\partial_{n+1} \oplus 1_L}{\longrightarrow} F_n \oplus L
\stackrel{\partial_{n} \oplus 0} {\longrightarrow} F_{n-1} \dots 
\stackrel{\partial_3}{\longrightarrow} F_{2} \stackrel{\partial_{2}} {\longrightarrow}F_1 \stackrel{\partial_{1}}
{\longrightarrow} F_{0}       
\eqno[4]
$$

\bigskip
\prop{The following sequence is exact:}

$$
F_m \stackrel{\partial_{m}} {\longrightarrow} F_{m-1} \stackrel{\partial_{m-1}} {\longrightarrow} \dots 
\stackrel{\partial_{n+2}}{\longrightarrow} F_{n+1} \oplus L \stackrel{f_n\circ(\partial_{n+1} \oplus 1_L)}{\longrightarrow} G_n \oplus K
\stackrel{\partial_{n}' \oplus 0} {\longrightarrow} G_{n-1} \stackrel{\partial_{n-1}'} {\longrightarrow}  \dots G_1 \stackrel{\partial_{1}'} {\longrightarrow} G_{0}      
\eqno[5]
$$

\hspace{2mm} Proof of Proposition: The image of $\partial_{n+1} \oplus 1_L$ is the kernel of $\partial_{n} \oplus 0$ and $f_n$ restricts to an isomorphism from the
kernel of $\partial_{n} \oplus 0$  to the kernel of $\partial_{n}' \oplus 0$, because it is the last term in a chain homotopy.  Hence the
image of $f_N \circ(\partial_{n+1} \oplus 1_L)$ is the kernel of $\partial_{n}' \oplus 0$ and the kernel of $f_N \circ(\partial_{n+1} \oplus
1_L)$ is the kernel of $\partial_{n+1} \oplus 1_L$.  
\hfill $\Box$

\bigskip
\hspace{2mm}  Let $\hat{f}_r = f_r$ for $r \leq n$ and $\hat{f}_r = 1$ for $r>n$.  Let $\hat{g}_r = g_r$ for 
$r \leq n$ and $\hat{g}_r = 1$ for $r>n$.

\bigskip
\prop {$\hat{f}$ is a chain map $[4]$ to $[5]$ and $\hat{g}$ is a chain map $[5]$ to $[4]$}.

\hspace{2mm} Proof of Proposition:  For commutativity, it is sufficient to check that $g_n f_n \partial_n = \partial_n$:   

\bigskip
$g_n f_n \partial_n = \partial_n +I_{n-1} \partial_{n-1} \partial_{n} = \partial_n$
\hfill $\Box$

\bigskip
\prop {$\hat{f}$ and $\hat{g}$ are homotopy inverse to each other}.

\hspace{2mm} Proof of Proposition:  Let $\hat{I}_r = I_r$ for $r \leq n-1$ and $\hat{I}_r = 0$ for $r > n-1$.  Let $\hat{J}_r = J_r$ for 
$r\leq n-1$ and $\hat{J}_r = 0$ for $r > n-1$.  Then $\hat{I}$ is a homotopy from $\hat{g}\hat{f}$ to $1$ and $\hat{J}$ is a homotopy 
from $\hat{f}\hat{g}$ to $1$.
\hfill $\Box$

\bigskip
So [1] is chain homotopy equivalent to [4] which is chain homotopy equivalent to [5].

\hspace{2mm}  $K$ is in the kernel of $\partial_n' \oplus 0$ so it is in the image of ${f_n\circ(\partial_{n+1} \oplus 1_L)}$.  Therefore composing
${f_n\circ(\partial_{n+1}\oplus 1_L)}$ with projection onto $K$ gives a surjection.  Let $SF$ denote the kernel of this surjection. As $K$ is free, this 
surjection splits.  Hence [5] can be written as

$$
F_m \stackrel{\partial_{m}} {\longrightarrow} F_{m-1} \stackrel{\partial_{m-1}} {\longrightarrow} \dots\dots
\stackrel{\partial_{n+2}}{\longrightarrow} SF \oplus K \stackrel{D}{\longrightarrow} G_n \oplus K \stackrel{\partial_{n}'} 
{\longrightarrow} G_{n-1} \stackrel{\partial_{n-1}'} {\longrightarrow} \dots\dots G_1 \stackrel{\partial_{1}'} {\longrightarrow} G_{0}      
\eqno[6]
$$
 
Note $SF \oplus K = F_{n+1} \oplus L $, so $SF$ is finitely generated, stably free.  With respect to the above decompositions, let $D$ be 
represented by 

\bigskip
$
D =
\left[\begin{array}{cc}  
 d&\phi\\
0&1_K
\end{array} \right]
$

\bigskip
The image of $\partial_{n+2}$ is contained in the kernel of $D$, hence in $SF$.  We have a sequence

$$
F_m \stackrel{\partial_{m}} {\longrightarrow} F_{m-1} \stackrel{\partial_{m-1}} {\longrightarrow} \dots\dots
\stackrel{\partial_{n+2}}{\longrightarrow} SF \stackrel{d}{\longrightarrow} G_n \stackrel{\partial_{n}'} 
{\longrightarrow} G_{n-1} \stackrel{\partial_{n-1}'} {\longrightarrow} \dots\dots G_1 \stackrel{\partial_{1}'} {\longrightarrow} G_{0}      
\eqno[7]
$$

\prop{$[7]$ is exact.}

\hspace{2mm} Proof of Proposition:  If an element of $SF$ is in the kernel of $d$ it is in the kernel of $D$, hence the image of
$\partial_{n+2}$.  If an element of $G_n$ is in the kernel of $\partial_{n}'$, then it is in the image of $D$.  However if $Dx \in G_n$ then
$x \in SF$ and $Dx=dx$.  
\hfill $\Box$
 
[7] is chain homotopy equivalent to 
 
$$
F_m \stackrel{\partial_{m}} {\longrightarrow} F_{m-1} \stackrel{\partial_{m-1}} {\longrightarrow} \dots\dots\dots
\stackrel{\partial_{n+2}}{\longrightarrow} SF \oplus K \stackrel{D'}{\longrightarrow} G_n \oplus K\stackrel{\partial_{n}'} 
{\longrightarrow} G_{n-1} \stackrel{\partial_{n-1}'} {\longrightarrow} \dots\dots\dots G_1 \stackrel{\partial_{1}'} {\longrightarrow} G_{0}      
\eqno[8]
$$

where $D'$ is represented by 
$
D =
\left[\begin{array}{cc}  
 d&0\\
0&1_K
\end{array} \right]
$ 

\bigskip
So to prove the theorem, it is sufficient to show that [6] is chain isomorphic to [8].

\hspace{2mm}  The image of $\phi$ is contained in the image of $d$ and $K$ is free, so there exists a map $\psi$ which makes the following
diagram commute:

$$
\begin{array}{ccc}  
K&&\\
\downarrow \psi &\searrow &\phi\\
SF &\stackrel{d}{\rightarrow}& G_n
\end{array} 
$$

Define a chain map $h$ from [6] to [8] by $h$ is the identity on all terms except $SF \oplus K$, where it is represented by the matrix

\bigskip
$
\left[\begin{array}{cc}  
1_{SF}&\psi\\
0&1_K
\end{array} \right]
$ 

\bigskip
Define a chain map $k$ from [8] to [6] by $k$ is the identity on all terms except $SF \oplus K$, where it is represented by the matrix

\bigskip
$
\left[\begin{array}{cc}  
1_{SF}&-\psi\\
0&1_K
\end{array} \right]
$

\bigskip

Then $hk=1_{[8]}$ and $kh=1_{[6]}$.  

\bigskip
\hspace{2mm}This completes the proof of Theorem 1.5.1.
\hfill $\Box$

\hspace{2mm} Let $G$ be a finite group.  We present two corollaries to theorem
1.5.1, lemma 1.4.18, and lemma 1.4.19.

\cor{ \it Suppose we have an exact sequence 
$$ 
0 \longrightarrow K\stackrel{} {\longrightarrow} \ZG^b \stackrel{\partial_{1}} {\longrightarrow} \ZG^a \stackrel{} {\longrightarrow}
\Z \longrightarrow 0
$$
Then we have a short exact sequence of the form
$$ 
0 \longrightarrow K\stackrel{} {\longrightarrow} \ZG^{b'} \stackrel{\partial_{1}'} {\longrightarrow} \ZG \stackrel{\epsilon} {\longrightarrow}
\Z \longrightarrow 0
$$
such that 
$$
ZG^b \stackrel{\partial_{1}} {\longrightarrow} \ZG^a
$$
is chain homotopy equivalent to 
$$
\ZG^{b'} \stackrel{\partial_{1}'} {\longrightarrow} \ZG
$$
}

\cor{Suppose we have exact sequences  
$$ 
0 \longrightarrow J\stackrel{} {\longrightarrow} \ZG^{c} \stackrel{\partial_2}{\longrightarrow}\ZG^{b} \stackrel{\partial_{1}} {\longrightarrow} \ZG^a \stackrel{} {\longrightarrow}
\Z \longrightarrow 0
$$
and
$$ 
0 \longrightarrow K\stackrel{} {\longrightarrow} \ZG^{b'} \stackrel{\partial_{1}'} {\longrightarrow} \ZG \stackrel{\epsilon} {\longrightarrow}
\Z \longrightarrow 0
$$
Then we have a short exact sequence of the form
$$ 
0 \longrightarrow J\stackrel{} {\longrightarrow} \ZG^{c'} \stackrel{\partial_2'}{\longrightarrow}\ZG^{b'} \stackrel{\partial_{1}'}
{\longrightarrow} \ZG \stackrel{\epsilon} {\longrightarrow}
\Z \longrightarrow 0
$$
such that
$$ 
\ZG^{c} \stackrel{\partial_2}{\longrightarrow}\ZG^{b} \stackrel{\partial_{1}} {\longrightarrow} \ZG^a 
$$
is chain homotopy equivalent to 
$$ 
\ZG^{c'} \stackrel{\partial_2'}{\longrightarrow}\ZG^{b'} \stackrel{\partial_{1}'}{\longrightarrow} \ZG 
$$
}

In particular, one may fix an exact sequence

$$ 
0 \longrightarrow K\stackrel{} {\longrightarrow} \ZG^{b} \stackrel{\partial_{1}} {\longrightarrow} \ZG \stackrel{\epsilon} {\longrightarrow}
\Z \longrightarrow 0
$$

with $K$ a minimal element of its stable class, so that any algebraic 2-complex is chain homotopy equivalent to one of the form 

$$ 
\ZG^{c} \stackrel{\partial_2}{\longrightarrow}\ZG^{b} \stackrel{\partial_{1}} {\longrightarrow} \ZG 
$$

\hspace{2mm}So we may parametrize algebraic 2-complexes, by maps from finitely generated
free modules to $K$.

\newpage

\chapter{The D2 problem for Dihedral groups}\label{}
\stepcounter{ch}
\setcounter{sec}{0}

\hspace{2mm}  This chapter will be concerned with the D(2)- problem for
dihedral groups of order $4n$.   In the first section we state and briefly discuss 
the general D(n)- problem, before focusing on the D(2) problem.  

\hspace{2mm}In \S2.2 we show that for dihedral groups whose order is a power of $2$, the Swan map (definition 1.2.7) is
surjective.  In \S2.3 we use this result to show that $D_8$ satisfies the D2
property.  

\hspace{2mm}In \S2.4 we give an exhaustive list of candidates for minimal elements of $\Omega_3(Z)$, over any 
dihedral group of order $4n$.  Finally, in \S2.5 we note that whilst not periodic over $\Z$, the dihedral groups of
order $8n+4$ are periodic over $\Zh$.

\sec{The D(n) problem}

\hspace{2mm}In this section we introduce the D(n) problem by outlining some of the content of \cite{John1}.

\bigskip
\hspace{2mm}We work in the category of CW- complexes and continuous maps. (See
\cite{Whit}, chapter II).  In particular, geometric complex will refer to a connected
finite dimensional, finite, CW complex, and geometric $n$- complex will refer to an $n$ dimensional geometric complex.  If geometric 
complexes $X$ and $Y$ are homotopy equivalent, we write $X \sim Y$.

\bigskip
\hspace{2mm}   Following the standard convention, we regard elements 
of $\tilde{X}$ as equivalence classes of paths in $X$, based at $b$.  Let $p:\tilde{X} \to X$ denote the covering 
map, which sends a path to its endpoint.  
.

\bigskip 
\hspace{2mm} 
If $X$ is a geometric complex, with base point $b$, then we may give $\tilde{X}$ the structure
of a geometric complex in the following way.
For each cell, $D$, of $X$, select a point $x_D \in {\rm Int}(D)$ and let $\phi_D:D \to X$ denote the natural insertion.  

\bigskip
\hspace{2mm}
Then 
for each $y
\in \tilde{X}$ such that
$p(y) = \phi_D(x_D)$, we may define a cell of $\tilde{X}$ which we denote $\tilde{D_y}$.  As a closed ball we identify its
points with those of $D$.  We define a map $\phi_y:\tilde{D_y} \to \tilde{X}$ which sends $t \in \tilde{D_y}$ to the 
path $y \circ \phi_D(z)$, where
$z$ is a path in  $\tilde{D_y}$ connecting $x_D$ to $t$.   We say $\tilde{D_y}$
lies above $D$.

\bigskip
\hspace{2mm} To reconstruct $\tilde{X}$ as a geometric complex, we need the cells $\tilde{D_y}$ together with attaching maps.  If 
the $n-1$ skeleton has been constructed, then
$\phi_y|_{\partial\tilde{D_y}}$ may be regarded as the attaching map for $\tilde{D_y}$ for each $n$- cell $D$, once we have 
shown $\phi_y({\partial\tilde{D_y}})$ is contained in the $n-1$ skeleton of $\tilde{X}$.  To see this, note that it is contained in the
universal cover of the $n-1$ skeleton of $X$.  Let $q$ be an element of this.  Then $p(q)= \phi_D(s)$ for some $n-1$ cell, $D$.  Let $z$ be a
path in $D$, from $s$ to $x_D$.  Then let $y=q\circ  \phi_D(z) $.  Then $q$ is in the $n-1$ cell $\tilde{D_y}$.

\bigskip
\hspace{2mm}We denote the corresponding complex of abelian groups $C_*(\tilde{X})$.
Let $G=\pi_1(X)$.  Then $G$ acts transitively on the cells
which lie above $D$, for each cell $D$ of $X$.  Hence the abelian groups in $C_*
(\tilde{X})$ have the structure of $\ZG$ modules.  These modules are free,
as no non-trivial element of $\pi_1(X)$ fixes any cell.  The boundary maps
are linear with respect to $\ZG$.  Hence $C_*(\tilde{X})$ is an algebraic complex over $\ZG$, which we will refer to as 
the associated algebraic complex of $X$.

\bigskip
\hspace{2mm}  If $B^n$ denotes an $n$- dimensional ball, then $C_*(\widetilde{X \vee B^n})$ is given by applying a simple homotopy equivalence
to  $C_*(\tilde{X})$, whereby a copy of $\ZG$ is
added to $C_n(\tilde{X})$ and $C_{n-1}(\tilde{X})$.

\bigskip
\hspace{2mm}  From our definition, any geometric $n$- complex  $X$ is connected.  Hence $\tilde{X}$ is necessarily  connected and simply 
connected.  As a result $C_*(\tilde{X})$ is exact at $C_1$ and the cokernel of the last boundary map $\partial_1$ is $\Z$.  Hence
$C_*(\tilde{X})$ is an algebraic $n$- complex.

\bigskip
\hspace{2mm}Given a geometric complex $X$ and a module, $A$, over $\pi_1(X)$,
 we define $H^n(X;A)=H^n((C_*\tilde{X},\partial_*);A)$
and $H_n(X;A)=H_n((C_*\tilde{X},\partial_*);A)$.  When a module is used in this
context, we refer to it as a coefficient bundle.  Note that to avoid ambiguity,
we must specify whether $A$ is an abelian group (in which case we are referring to the (co)homology of $C_*(X)$),
or a module over the fundamental
group.  Unless otherwise stated, in the context of coefficients for a (co)homology group,$\Z$ will denote the abelian group, rather than
the module with trivial group action.

\bigskip
\hspace{2mm}  The D(n) problem asks when a 
geometric $n+1$- complex, $X$, is homotopy
equivalent to a geometric $n$ complex.  C.T.C.Wall has shown that for $n\geq3$, a
necessary and sufficient condition is that $H^{n+1}(X;A)=0$ for all coefficient
bundles (see \cite{Wall}).  Note that this condition is equivalent to saying that 
$(C_*\tilde{X},\partial_*)$ is chain homotopy equivalent to an algebraic $n$-
complex, as if $H^{n+1}(X;C_{n+1}(\tilde{X}))=0$ then the following diagram commutes:

\bigskip
\xymatrix{
C_{n+1}(\tilde{X})\ar[d]_1 \ar[r]^{\partial_{n+1}} 
& C_{n}(\tilde{X})\ar@{.>}[dl]^f\\
C_{n+1}(\tilde{X})} 

\bigskip

for some map $f$.  Hence $\partial_{n+1}$ splits and $C_{n+1}(\tilde{X})$ has
some stably free complement $S$ in $C_{n}(\tilde{X})$.  Finally note a simple
homotopy equivalence connects the complexes

\bigskip
\xymatrix{
C_{n+1}(\tilde{X}) \ar[r]^{\partial_{n+1}} 
& C_{n+1}(\tilde{X}) \oplus S
\ar[r]^{\quad \partial_{n}} 
& \cdots \ar[r]^{\partial_{1}} 
&C_{0}(\tilde{X})}

\bigskip

and

\xymatrix{\qquad\qquad
0 \qquad \ar[r] 
& \qquad S
\ar[r]^{\partial_{n}} 
& \cdots \ar[r]^{\partial_{1}} 
&C_{0}(\tilde{X})}

\bigskip
\hspace{2mm}  If $n>0$, then in fact this complex will be chain homotopy
equivalent to 

\bigskip
\xymatrix{& S \oplus C_{n+1}(\tilde{X})
\ar[r]^{\partial_{n} \oplus 1} 
&C_{n-1}(\tilde{X}) \oplus C_{n+1}(\tilde{X})\ar[r]
& \cdots \ar[r]^{\partial_{1}} 
&C_{0}(\tilde{X})}

This is a complex containing only free modules.

\bigskip
\hspace{2mm}  An example of why the weaker condition, $H_2(X;\Z)=0$ is not
sufficient is the Klein bottle, $K$.  It is constructed by identifying the ends of an
oriented cylinder, so that the induced orientations on the ends agree.  Whilst
$H_2(K;\Z)=0$, it is not homotopy equivalent to a geometric $1$- complex, because
$H^2(K;\Z)= \Z_2$.

\lem{Let $X$ be a geometric $1$- complex and suppose $H^1(X;A)=0$ for 
all coefficient bundles $A$.  Then $X$ is homotopy equivalent to a set of
points.}

\proof{We know that the associated algebraic complex of $X$ is chain 
homotopy equivalent to a $0$- algebraic complex.  Hence we have a short exact
sequence: \newline  $$0 \to F \to \Z$$
\newline
where $F$ is free over $\pi_1(X)$.  Hence the augmentation ideal of  $\pi_1(X)$
must be $0$.  We may conclude that $\pi_1(X)$ is trivial.  Hence $X$ must be a
geometric $1$- complex with no cycles.  It is therefore a set of trees, each of
which is contractible to a point.}





\bigskip
\hspace{2mm}  Now suppose that $X$ is a geometric $2$- complex and 
suppose $H^2(X;A)=0$ for all coefficient bundles $A$.  Then the associated
algebraic complex of $X$ is an algebraic $2$- complex.  In
particular $H_2(X;\Z)=0$, so the associated algebraic complex of $X$ is exact. 
Consequently $X$ is homotopy equivalent to $k(\pi_1(X),1)$.  If 
$\pi_1(X)$ is free, then $k(\pi_1(X),1)$ must be homotopy equivalent 
to a wedge of circles which is a geometric $1$- complex.

\bigskip
\hspace{2mm}In fact $\pi_1(X)$ must be free as the following theorem implies:

\thm{(Stallings-Swan)  Suppose that we have an exact sequence over $\ZG$, for
some group $G$, of the form:
\newline
\newline
$$0 \to P_1 \to P_0 \to \Z \to 0$$
\newline
where $P_0$ and $P_1$ are projective.  Then $G$ is free}

\bigskip
\hspace{2mm}  We know that $\pi_1(X)$ satisfies the hypothesis of this theorem,
as the associated algebraic complex of $X$ is chain homotopy equivalent to an algebraic
$1$- complex.

\bigskip
\hspace{2mm}  Note also, that the theorem may be stated group theoretically.  
The existence of the exact sequence 

$$0 \to P_1 \to P_0 \to \Z \to 0$$

is equivalent to $H^2(G;A)$ vanishing for all coefficient bundles $A$.  By 
proposition 1.3.3, this in
turn is equivalent to saying that every group surjection, $j:E \to G$, with
${\rm ker}(j)$ abelian, will split.  So the theorem could be stated 

\thm{If every group surjection, with abelian kernel, to a group $G$, splits,
then $G$ is free.}

\hspace{2mm}  We may conclude that the D(0), D(1), and 
D(n), $n \geq 3$ problems all have the same solution.  In each case, a geometric
complex is homotopy equivalent to one of dimension 1 less than it, precisely if
its top cohomology group vanishes with respect to all coefficient bundles.  We now
discuss the extent to which this is known to hold for D(2).




\hspace{2mm}Due to the relationship between $2$- complexes and presentations of
fundamental groups, it is natural to
parametrize the D(2) problem over the fundamental group.  So for a finite group $G$, the
D(2) problem is:

\bigskip
\hspace{2mm} 
{\bf Let $X$ be a finite geometric $3$- complex, with fundamental group
$G$, and with vanishing third cohomology group for every coefficient bundle. 
Must $X$ be homotopy equivalent to a finite geometric $2$- complex?}

\bigskip
\hspace{2mm}If the answer is yes, then we say $G$ satisfies the  
D(2) property.

\bigskip
\hspace{2mm}We now outline some of the arguments in \cite{John1} which show that for each finite group $G$, the D(2) 
property for $G$ is equivalent
to every algebraic $2$- complex  over $\ZG$ being realized by a geometric $2$- complex.

\bigskip
\hspace{2mm}  If $X$ satisfies the hypothesis of the D(2) problem, then 
$C_*(\tilde{X}) \sim Y$  for some algebraic
$2$- complex, $Y$.  Suppose every such algebraic $2$- complex, was 
chain homotopy equivalent to $C_*(\tilde{Z})$, for some geometric $2$- complex, $Z$.   Then we would have $C_*(\tilde{X}) \sim 
C_*(\tilde{Z})$.  From \cite{John1}, theorem 59.4, 
we would then have $X \sim Y$ would then have.

\bigskip
\hspace{2mm}Conversely, let $(F_*, d_*)$ be an algebraic $2$- complex 
of free modules, over $\ZG$, which is not
realizable up to homotopy as a geometric $2$- complex. Take any finite
presentation of $G$ and let $Y$ denote its Cayley complex. By theorem 1.1.1
there exists a free finitely generated module $F$ and a number $n$, such that $Y \vee \bigvee_{1=1}^n S^2$, 
is chain homotopy equivalent to 

$$
F_2 \oplus F \stackrel{d_2 \oplus 0}{\longrightarrow}
F_1 \stackrel{d_1 }{\to} F_0
$$

\hspace{2mm}Let $f$ denote the homotopy equivalence between this complex and
$C_*(\tilde{Y'})$.  Let the $e_i$ be a basis for $F$.  As $F$ lies in the kernel
of $d_2 \oplus 0$ and $f_2$ restricts to a map between second homology
groups, for each $i$,
$f_2(e_i)$ is an element of $\pi_2(Y')$.  Starting with $Y'$, attach 
one $3$- cell, $B_i$ for each $i$, where
the attaching map from $\partial B_i$, is given by $f_2(e_i)$, the element of
$\pi_2(Y')$.  Let $Z$ denote the resulting geometric $3$- complex.  Note that
$C_3(\tilde{Z}) \cong F $ and $C_i(\tilde{Z}) \cong C_i(\tilde{Y'})$ for
$i=0,1,2$.  

\bigskip
\hspace{2mm}
The algebraic complex $(F_*, d_*)$ is related, via a simple homotopy
equivalence to 

$$
F\stackrel{(0,1)}{\longrightarrow}F_2 
\oplus F \stackrel{d_2 \oplus 0}{\longrightarrow}
F_1 \stackrel{d_1 }{\to} F_0
$$

\hspace{2mm}
This is in turn chain homotopy equivalent to $C_*(\tilde{Z})$ via the following homotopy
equivalence:

\xymatrix{F \ar[d]_1 \ar[r]^{(0,1)\quad\quad}  &F_2 \oplus F\ar[r]^{ \quad
d_2 \oplus
0}\ar[d]_{f_2}
&F_1 \ar[r]^{d_1} \ar[d]_{f_1}&F_0 \ar[d]_{f_0} \\
C_3(\tilde{Z})\ar[r]^{f_2}&C_2(\tilde{Z})\ar[r]^{\partial_2}
&C_1(\tilde{Z})\ar[r]^{\partial_2}&C_0(\tilde{Z})}

\bigskip
where the $\partial_i$ are the boundary maps for $Y$.

\bigskip
\hspace{2mm}  $Z$  cannot be homotopy equivalent to a geometric $2$- complex,
as the associated algebraic complex of such a complex would be chain homotopy
equivalent to $(F_*, d_*)$, contradicting the hypotheses that $(F_*, d_*)$ is 
not realizable up to homotopy as a geometric $2$- complex.  However, $Z$ clearly
satisfies the hypotheses of the general D(2) problem as its associated algebraic
complex is chain homotopy equivalent to  $(F_*, d_*)$.

\bigskip
\hspace{2mm}  Hence a finite group $G$ satisfies the D(2) property if and 
only if every
algebraic $2$- complex is realizable up to homotopy, as a
geometric $2$- complex.

\bigskip
\hspace{2mm}  If $G$ is a finite group then we also have the following
characterization of homologically $2$ dimensional complexes.

\lem{Over $\ZG$ an algebraic $3$- complex, $(F_i, d_i)$ is chain
homotopy equivalent to an algebraic $2$- complex if and only if $H_2((F_i, d_i),\Z)$
is torsion free and  $H_3((F_i, d_i),\Z)=0$.}

\proof{As $F_3^*$ is projective, $d_3^*$ has a right inverse.  
$d_3=d_3^{**}$ because the cokernel of $d_3$ is torsion free.  Hence $d_3$ is
the inclusion of a free summand into $F_2$.} 

\bigskip
\hspace{2mm}A finite geometric $2$- complex will always be homotopy equivalent
to the Cayley complex of a finite presentation.  Hence our problem reduces to a
purely algebraic one.  Following Fox, we may define an algorithm for constructing the
associated algebraic complex of a finite Cayley complex directly from the finite
presentation,
and our question becomes:

\bigskip
\hspace{2mm}  Is every algebraic $2$- complex 
chain homotopy equivalent to one constructed from a finite presentation?

\bigskip
\hspace{2mm}Let $F_2 \stackrel{\partial_2}{\rightarrow}
F_1 \stackrel{\partial_1}{\rightarrow}
F_0$ be an algebraic  $2$- complex. The kernel of $\partial_2$ will be referred to as $\pi_2$
of the algebraic $2$- complex.  
As the kernel of $\partial_1$ is equal to the
image of $\partial_2$ and the cokernel of  $\partial_1$ is isomorphic to the
$\ZG$- module $\Z$, we have that $\pi_2(F_*,\partial_*) \in \Omega_3(\Z)$.
  
\prop{Let $X$ be a geometric $2$- complex, with $\pi_1(X)=G$, for some finite group $G$.  Then $\pi_2(C_*(\tilde{X})) \cong \pi_2(X)$ as 
modules over $\ZG$.} 

\proof{We first define an isomorphism of abelian groups, and then show that it is $\ZG$- linear.  From the definition, we have a natural
identification of $H_2(\tilde{X}; \Z)$ with $\pi_2(C_*(\tilde{X}))$.  As $\tilde{X}$ is simply connected,
the Hurewicz isomorphism theorem tells us that the Hurewicz homomorphism induces an isomorphism
$h:\pi_2(\tilde{X}) \to H_2(\tilde{X}; \Z)$.  The covering 
map $p:\tilde{X} \to X$ induces an isomorphism $p_*:\pi_2(\tilde{X})\to \pi_2(X)$.
\newline 
${}$\hspace{2mm} Hence 
we have an isomorphism of Abelian groups $hp_*^{-1}:\pi_2(X) \to \pi_2(C_*(\tilde{X}))$, which lifts any element of $\pi_2(X)$ to
$\tilde{X}$ and then applies the resulting map $H_2(S^2;Z)\to H_2(\tilde{X};\Z)$ to 
the positive generator of $H_2(S^2;\Z)$.  We must check that this isomorphism
respects the action of $G$.  
\newline
${}$\hspace{2mm} Given $g \in G$, and $\alpha \in \pi_2(X)$, by construction, we
have that 
$\alpha g$ lifts to $(p_*^{-1}\alpha)g$, attached to the base point of $\tilde{X}$ via a lift of $g$.  Hence the 
element of $H_2(\tilde{X};\Z)$ which it sends the positive generator of $H_2(S^2;Z)$
to, is $(hp_*^{-1}\alpha)g$.  So $hp_*^{-1}(\alpha g)= (hp_*^{-1}\alpha)g$.
}

\define{}{The Euler characteristic of an algebraic complex of free modules is 
the alternating sum of the dimensions of its modules.}

\bigskip
\hspace{2mm}Our principal method of showing that a finite group, $G$, satisfies the D(2)
property will be to show that any algebraic $2$- complex over $G$, is determined 
up to chain homotopy equivalence by its Euler characteristic.  Once
that is done, it is only necessary to show that the algebraic complex with
minimal Euler characteristic is chain homotopy equivalent to one arising from a presentation. 
Then all the others will necessarily be chain homotopy equivalent to the
algebraic complex arising from the presentation, together with the appropriate
number of trivial relations, "$e=e$" added.

\bigskip
\hspace{2mm}  A necessary precursor to implementing this method is to show that
$\pi_2(X)$ for any algebraic two-complex, $X$, over a finite group $G$, is determined by its $\Z$- rank. 
Clearly, the Swan - Jacobinski theorem is a useful tool in this endeavor.  

\bigskip
\hspace{2mm}  Given an algebraic $2$- complex $X$, we may tensor it with $\R$. Then the "Whitehead trick" together with semi-simplicity of
$\RG$ tell us that $\R \otimes_{\Z}\pi_2(X)=\R \otimes_{\Z}IG \oplus \RG^k$ for some integer $k$.  Hence given a 
minimal element, $J \in\Omega_3(\Z)$, we know that $J \oplus \ZG$ satisfies the Eichler condition.  Any non-minimal 
element of $\Omega_3(\Z)$ is therefore equal to $J \oplus \ZG^r$ for some integer $r$.  The remaining problem is to prove that $J$ is the
unique minimal element of $\Omega_3(\Z)$.

\bigskip
\hspace{2mm} If for some group it is shown that there is a unique potential
$\pi_2$ for each Euler characteristic, it remains to show that the $\pi_2$ determines the
algebraic $2$- complex.  From theorem 1.2.8, we know that the surjectivity of the Swan map is sufficient to show this.

\bigskip
\hspace{2mm}We end this section with an example of a group satisfying the D(2) property.

\lem{The group $\{e\}$ satisfies the D(2) property.}

\hspace{2mm}Proof:  Over $\Z$, any finitely generated torsion free module is isomorphic to
$\Z^k$ for some $k$.  Also $\Z$ is a Noetherian ring.  Consequently,
any torsion free map $\Z^a \to \Z^b$, is the composition of projection onto a
free summand of $\Z^a$, with inclusion into a free summand of $\Z^b$.
\newline
${}$ \hspace{2mm} Hence any algebraic $2$- complex over $\Z$, is related through
simple homotopy equivalences to 
\newline
$\Z^k \to 0 \to \Z$
\newline
for some $k$.  This is realized geometrically as the wedge of $k$ copies of
$S^2$. 

\hfill $\Box$

\sec{Surjectivity of the Swan map.}

\hspace{2mm}The D(2) property has been verified for dihedral groups of order $4n+2$ in \cite{John1}.  The remainder of this chapter is concerned
with the D(2) property for dihedral groups of order $4n$.  In the next section we prove that $D_8$ satisfies the D(2) property.  We
begin with more general considerations.

\bigskip
\hspace{2mm} Let $D_{4n}$ be the group given by the presentation, $<a,b \,\, |a^{2n}=b^2=e,\,\,aba =b>$.  $\Sigma$ will denote 
$\sum_{i=0}^{2n-1} a^i$.  This presentation has a Cayley complex, which in turn has an associated algebraic complex.  This is an exact
sequence over $\ZD$:

$$
J\hookrightarrow \ZD^3 \stackrel{\partial_2}{\longrightarrow} \ZD^2 \stackrel{\partial_1}{\longrightarrow} \ZD \stackrel{\epsilon}{\twoheadrightarrow} \Z
\eqno (1)
$$

$\epsilon$ is determined by mapping $1 \in \ZD$ to $1 \in \Z$.  $J$ is the kernel of $\partial_2$. Let $e_1$, $e_2$ denote basis elements of $\ZD^2$.  Then $\partial_1 e_1 = a-1, \, \, \partial_1 e_2 = b-1$.  

\hspace{2mm} Let $E_1, \, E_2, \, E_3$ be basis elements of $\ZD^3$, which correspond to the relations in the presentation so that:

\bigskip
$
\partial_2 E_1 = e_1 \Sigma \newline
\partial_2 E_2 = e_2 (1+b) \newline
\partial_2 E_1 = e_1 + e_2a + e_1ba - e_2= e_1(1+ba)+e_2(a-1)
$

\bigskip
\hspace{2mm}  With respect to the basis $\{E_1,\, E_2,\,E_3\}$ and the basis $\{e1, \, e2\}$, $\partial_2$ is given by ;

$$
\left[ \begin{array}{ccc} \Sigma & 0 & 1+ba \\ 0 & 1+b & a-1 \end{array} \right]
$$

Let 

\bigskip
$\alpha_0 = 1+a+b$

\bigskip
$\alpha_1 = \left[ \begin{array}{cc} 1+a-ba & b-1 \\ 0 & 1 \end{array} \right]$

\bigskip
$\alpha_2 = \left[ \begin{array}{ccc}  1+a-ba & 0&0 \\ 0 & 1&0 \\ 0&0&1\end{array} \right]$

\bigskip
\hspace{2mm}It is easily verified that:

\prop{The following diagram commutes: 
\newline
\begin {eqnarray*}
J\hookrightarrow \ZD^3 \stackrel{\partial_2}{\longrightarrow} \ZD^2 \stackrel{\partial_1}{\longrightarrow}
\ZD \stackrel{\epsilon}{\twoheadrightarrow} \Z \,\,\,\\
\downarrow \theta  \,\,\,\,\, \quad \downarrow \alpha_2 \qquad \qquad \downarrow \alpha_1 \qquad \downarrow \alpha_0 \qquad\downarrow 3
\\
J\hookrightarrow \ZD^3 \stackrel{\partial_2}{\longrightarrow} \ZD^2 \stackrel{\partial_1}{\longrightarrow}
\ZD \stackrel{\epsilon}{\twoheadrightarrow} \Z \,\,\,
\end {eqnarray*}
\newline
where $\theta$ is the restriction of $\alpha_2$.}

\bigskip
\hspace{2mm} For the remainder we will assume $3$ coprime to $n$.  Note that if we regard the above diagram as a diagram of commutative $\Z$-
modules and $\Z$- linear maps,  there are well defined integer determinants for all the maps in the chain map.  A map is an 
isomorphism if and only if it has
determinant $\pm1$.  (As the property of being an isomorphism is dependent only on surjectivity and injectivity, it does not depend on
whether we are regarding modules as being over $\ZD$, or $\Z$). 

\hspace{2mm} Note also that over $\Z$, all the maps in the exact sequences above are quotienting of a summand, followed by inclusion
of a summand.  Consequently, the following proposition holds;   

\bigskip
\prop{$3{\rm Det}(\theta) {\rm Det}(\alpha_1)= {\rm Det}(\alpha_2){\rm Det}(\alpha_0)$}

\hspace{2mm} Proof:  Let $u$ be the restriction of $\alpha_1$ to the kernel of $\partial_1$ and let $v$ be the restriction of
$\alpha_0$ to the kernel of $\epsilon$.    Then by the previous discussion, we have 

$$
3{\rm Det}(\theta){\rm Det}(\alpha_1)= {\rm Det}(\theta){\rm Det}(u) {\rm Det}(v)3 = {\rm Det}(\alpha_2){\rm Det}(\alpha_0)
$$
\hfill $\Box$

\bigskip
\prop{${\rm Det}(1+a+b)=-3$}

\hspace{2mm} Proof:  Let $A$ be the matrix for left multiplication by $1+a+b$ in the regular representation, with basis
$\{a^{2n-1},a^{2n-2},\cdots, a,1,ba^{2n-1},ba^{2n-2},\cdots, ab,b\}$.  Then the upper right quadrant of $A$ and the lower left quadrant
of $A$ are copies of the identity matrix.  The upper left quadrant has $1$'s along the diagonal and immediately
above as well as a $1$ in the bottom left corner.  The lower right quadrant has $1$'s along the diagonal and
immediately below, as well as a $1$ in the top right corner.  All the other entries in $A$ are $0$.  

\newpage
\hspace{2mm} For example, if $n$ were equal to $4$, the matrix $A$ would be

$$ \tiny
\left[\begin{array}{ccccccccccccccccc}  
1&1&0&0&0&0&0&0&&1&0&0&0&0&0&0&0\\
0&1&1&0&0&0&0&0&&0&1&0&0&0&0&0&0 \\
0&0&1&1&0&0&0&0&&0&0&1&0&0&0&0&0\\
0&0&0&1&1&0&0&0&&0&0&0&1&0&0&0&0 \\
0&0&0&0&1&1&0&0&&0&0&0&0&1&0&0&0\\
0&0&0&0&0&1&1&0&&0&0&0&0&0&1&0&0 \\
0&0&0&0&0&0&1&1&&0&0&0&0&0&0&1&0\\
1&0&0&0&0&0&0&1&&0&0&0&0&0&0&0&1 \\
 & & & & & & & && & & & & & & & \\
1&0&0&0&0&0&0&0&&1&0&0&0&0&0&0&1\\
0&1&0&0&0&0&0&0&&1&1&0&0&0&0&0&0 \\
0&0&1&0&0&0&0&0&&0&1&1&0&0&0&0&0\\
0&0&0&1&0&0&0&0&&0&0&1&1&0&0&0&0 \\
0&0&0&0&1&0&0&0&&0&0&0&1&1&0&0&0\\
0&0&0&0&0&1&0&0&&0&0&0&0&1&1&0&0 \\
0&0&0&0&0&0&1&0&&0&0&0&0&0&1&1&0\\
0&0&0&0&0&0&0&1&&0&0&0&0&0&0&1&1 \\
\end{array} \right]
$$

\hspace{2mm} Label the rows of $A$, $v_1, v_2.....v_{4n}$. We will perform row
operations.

\hspace{2mm}  First let $v_{2n}'= v_{2n}-v_1+v_2-v_3....-v_{2n-1}$.  Now let $v_{2n}''=v_{4n}$ and $v_{4n}''=v_{2n}'$.  Let the
remaining $v_{i}''=v_{i}$.  This swap causes a change of sign in the determinant, so the matrix with rows $v_i''$ has determinant -Det$A$. 
In the case $n=4$, this matrix is

$$ \tiny
\left[\begin{array}{ccccccccccccccccc}  
1&1&0&0&0&0&0&0&&1&0&0&0&0&0&0&0\\
0&1&1&0&0&0&0&0&&0&1&0&0&0&0&0&0 \\
0&0&1&1&0&0&0&0&&0&0&1&0&0&0&0&0\\
0&0&0&1&1&0&0&0&&0&0&0&1&0&0&0&0 \\
0&0&0&0&1&1&0&0&&0&0&0&0&1&0&0&0\\
0&0&0&0&0&1&1&0&&0&0&0&0&0&1&0&0 \\
0&0&0&0&0&0&1&1&&0&0&0&0&0&0&1&0\\
0&0&0&0&0&0&0&1&&0&0&0&0&0&0&1&1 \\
 & & & & & & & && & & & & & & & \\
1&0&0&0&0&0&0&0&&1&0&0&0&0&0&0&1\\
0&1&0&0&0&0&0&0&&1&1&0&0&0&0&0&0 \\
0&0&1&0&0&0&0&0&&0&1&1&0&0&0&0&0\\
0&0&0&1&0&0&0&0&&0&0&1&1&0&0&0&0 \\
0&0&0&0&1&0&0&0&&0&0&0&1&1&0&0&0\\
0&0&0&0&0&1&0&0&&0&0&0&0&1&1&0&0 \\
0&0&0&0&0&0&1&0&&0&0&0&0&0&1&1&0\\
0&0&0&0&0&0&0&0&&-1&1&-1&1&-1&1&-1&1 \\
\end{array} \right]
$$

\bigskip
\bigskip
\hspace{2mm}  For each $2n+1 \leq i \leq 4n-2$, let $v_i'''=v_i''+v_{i+1}''-v_{i-2n}''$.  

\hspace{2mm} Let $v_{4n-1}'''=
v_{4n-1}''+v_{2n}''-v_{2n-1}$''.  

\hspace{2mm}For $i \leq 2n$ let $v_i'''=v_i''$.  

\newpage
\hspace{2mm}When $n=4$, the matrix with rows $v_i'''$ is

$$ \tiny
\left[\begin{array}{ccccccccccccccccc}  
1&1&0&0&0&0&0&0&&1&0&0&0&0&0&0&0\\
0&1&1&0&0&0&0&0&&0&1&0&0&0&0&0&0 \\
0&0&1&1&0&0&0&0&&0&0&1&0&0&0&0&0\\
0&0&0&1&1&0&0&0&&0&0&0&1&0&0&0&0 \\
0&0&0&0&1&1&0&0&&0&0&0&0&1&0&0&0\\
0&0&0&0&0&1&1&0&&0&0&0&0&0&1&0&0 \\
0&0&0&0&0&0&1&1&&0&0&0&0&0&0&1&0\\
0&0&0&0&0&0&0&1&&0&0&0&0&0&0&1&1 \\
 & & & & & & & && & & & & & & & \\
0&0&0&0&0&0&0&0&&1&1&0&0&0&0&0&1\\
0&0&0&0&0&0&0&0&&1&1&1&0&0&0&0&0 \\
0&0&0&0&0&0&0&0&&0&1&1&1&0&0&0&0\\
0&0&0&0&0&0&0&0&&0&0&1&1&1&0&0&0 \\
0&0&0&0&0&0&0&0&&0&0&0&1&1&1&0&0\\
0&0&0&0&0&0&0&0&&0&0&0&0&1&1&1&0 \\
0&0&0&0&0&0&0&0&&0&0&0&0&0&1&1&1\\
0&0&0&0&0&0&0&0&&-1&1&-1&1&-1&1&-1&1 \\
\end{array} \right]
$$

\bigskip
\bigskip
\hspace{2mm}
In general, the matrix with rows $v_i'''$ has an upper
triangular top left quadrant, with $1$'s along the diagonal and a lower left quadrant with no non-zero entries.  Let $B$ denote the
lower right quadrant.  Then ${\rm Det}(1+a+b)=-{\rm Det}(B)$.     

\hspace{2mm} Cycle the top $2n-1$ rows of $B$ upwards to get the matrix $B'$.  As this is a cycle of odd length, ${\rm Det}(B')={\rm
Det}(B)$.  When $n=4$, $B'$ is

$$ \tiny
\left[\begin{array}{ccccccccc}  
1&1&1&0&0&0&0&0 \\
0&1&1&1&0&0&0&0\\
0&0&1&1&1&0&0&0 \\
0&0&0&1&1&1&0&0\\
0&0&0&0&1&1&1&0 \\
0&0&0&0&0&1&1&1\\
1&1&0&0&0&0&0&1\\
-1&1&-1&1&-1&1&-1&1 \\
\end{array} \right]
$$

\bigskip
\hspace{2mm} Label the rows of $B'$ as $w_1,...,w_{2n}$.  Set $u_{i}=w_i-w_{i+1}$ for $i=1,2, \cdots,2n-3$.  Let $B''$ denote the matrix with
rows $u_i$.  After these row operations, we have ${\rm Det}(1+a+b)=-{\rm Det}(B'')$

\newpage
\hspace{2mm} When $n=4$, $B''$ is

$$ \tiny
\left[\begin{array}{ccccccccc}  
1&0&0&-1&0&0&0&0 \\
0&1&0&0&-1&0&0&0\\
0&0&1&0&0&-1&0&0 \\
0&0&0&1&0&0&-1&0\\
0&0&0&0&1&0&0&-1 \\
0&0&0&0&0&1&1&1\\
1&1&0&0&0&0&0&1\\
-1&1&-1&1&-1&1&-1&1 \\
\end{array} \right]
$$

\hspace{2mm}  We must consider two cases: $n$ congruent to $1$ modulo $3$ and $n$ congruent to $2$ modulo $3$.

\hspace{2mm}  If $n=1$ modulo $3$ then replace $u_{2n-1}$ with 

$u_{2n-1}-u_1-u_2-u_4-u_5-u_7-u_8\cdots - u_{2n-3}$.

\hspace{2mm}Also, replace $u_{2n}$ with 

$u_{2n} + (u_1-u_2+u_3)+(u_7-u_8+u_9)+(u_{13}-u_{14}+u_{15})...+(u_{2n-7}-u_{2n-6}+u_{2n-5})$.  

We are left with a
matrix with $1$'s along the diagonal and $0$'s below, except in the last four columns.  The 4 by 4 matrix in the bottom right corner is 

$$\left[\begin{array}{cccc} 1&0&0&-1 \\ 0&1&1&1 \\ 0&0&1&2 \\ 0&0&-1&1 \end{array} \right] $$

$${\rm Det}(1+a+b)=-{\rm Det}\left[\begin{array}{cccc} 1&0&0&-1 \\ 0&1&1&1 \\ 0&0&1&2 \\ 0&0&-1&1 \end{array} \right] 
=-{{\rm Det}}\left[\begin{array}{cccc} 1&0&0&-1 \\ 0&1&1&1 \\ 0&0&1&2 \\ 0&0&0&3 \end{array} \right]=-3$$

\bigskip
\hspace{2mm}  If $n=2$ modulo $3$ then replace $u_{2n-1}$ with 

$u_{2n-1}-u_1-u_2-u_4-u_5-u_7-u_8....-u_{2n-5}$.  

\hspace{2mm}Also, replace $u_{2n}$ with 

$u_{2n} + (u_1-u_2+u_3)+(u_7-u_8+u_9)+(u_{13}-u_{14}+u_{15})...+(u_{2n-9}-u_{2n-8}+u_{2n-7})$.  

We are left with a
matrix with $1$'s along the diagonal and $0$'s below, except in the last four columns.  The 4 by 4 matrix in the bottom right corner is 

$$\left[\begin{array}{cccc} 1&0&0&-1 \\ 0&1&1&1 \\ 1&1&0&1 \\ -1&1&-1&1 \end{array} \right] $$

\bigskip
\bigskip
$${\rm Det}(1+a+b)=-{\rm Det}\left[\begin{array}{cccc} 1&0&0&-1 \\ 0&1&1&1 \\ 1&1&0&1 \\ -1&1&-1&1 \end{array} \right] 
=-{\rm Det}\left[\begin{array}{cccc} 1&0&0&-1 \\ 0&1&1&1 \\ 0&1&0&2 \\ 0&1&-1&0 \end{array} \right]$$

\bigskip
$$=-{\rm Det}\left[\begin{array}{cccc} 1&0&0&-1 \\ 0&1&1&1 \\ 0&0&-1&1 \\ 0&0&-2&-1 \end{array} \right] = 
-{\rm Det}\left[\begin{array}{cccc} 1&0&0&-1 \\ 0&1&1&1 \\ 0&0&-1&1 \\ 0&0&0&-3 \end{array} \right] = -3$$

\hfill $\Box$

\bigskip
\prop{${\rm Det}(2-b)=3^{2n}$}

\hspace{2mm} Proof:  Let $A$ be the matrix for $2-b$ in the regular representation, with basis
$\{1, b, a,ba, a^2, ba^2......a^{2n-1},ba^{2n-1} \}$ .  Then $A$ consists of $2n$ two by two blocks of the form

$$\left[\begin{array}{cc}  2 & -1 \\ -1 & 2 \end{array} \right]$$

along the diagonal.  Hence ${\rm Det} (A) = 3^{2n}$

\hfill $\Box$

\bigskip
\prop{${\rm Det}(1+a-ba) \neq 0$}

\bigskip
\hspace{2mm} Proof:  Let $\alpha_2'=\left[\begin{array}{ccc}  2-b & 0&0 \\ 0 & 1&0 \\ 0&0&1\end{array}
\right]$

\bigskip
The following diagram commutes:

\begin {eqnarray*}
J\hookrightarrow \ZD^3 \stackrel{\partial_2}{\longrightarrow} \ZD^2 \stackrel{\partial_1}{\longrightarrow}
\ZD \stackrel{\epsilon}{\twoheadrightarrow} \Z \,\,\,\\
\downarrow \eta  \,\,\,\, \quad \downarrow \alpha_2' \qquad \qquad \downarrow \alpha_1 \qquad \downarrow \alpha_0 \qquad\downarrow 3
\\
J\hookrightarrow \ZD^3 \stackrel{\partial_2}{\longrightarrow} \ZD^2 \stackrel{\partial_1}{\longrightarrow}
\ZD \stackrel{\epsilon}{\twoheadrightarrow} \Z \,\,\,
\end {eqnarray*}

where $\eta$ is the restriction of $\alpha_2'$

\bigskip
Therefore $3{\rm Det}(\eta) {\rm Det}(\alpha_1)= {\rm Det}(\alpha_2'){\rm Det}(\alpha_0)$

\bigskip
So $3*{\rm Det}(\eta){\rm Det}(1+a-ba) = -3*3^{2n}$.  Hence ${\rm Det}(1+a-ba)$ cannot be $0$.
\hfill $\Box$ 

\bigskip
\prop{$\theta$ is an isomorphism}

\hspace{2mm} Proof: $3{\rm Det}(\theta) {\rm Det}(\alpha_1)= {\rm Det}(\alpha_2){\rm Det}(\alpha_0)$

\bigskip 
\hspace{2mm} Therefore $3{\rm Det}(\theta){\rm Det}(1+a-ba)= -3 {\rm Det}(1+a-ba)$.  As ${\rm Det}(1+a-ba)$ is non-zero, we can conclude that 

\bigskip
\hspace{2mm} ${\rm Det}(\theta)=-1$.  

\bigskip
\hspace{2mm}Hence $\theta$ is an isomorphism.

\hfill $\Box$ $\,\,\,$

\bigskip
\hspace{2mm}  Recall the definition of the Swan map (Definition 1.2.7).

\cor{If $3 \in (\Z_{4n})^*$ then $3$ is in the image of the Swan Map ${\rm Aut(J)} \to (\Z_{4n})^*$ }

\bigskip
\hspace{2mm}  Let us now consider dihedral groups of order $2^m$ for $m \geq 2$.  Clearly $2^m$ is divisible by $4$
and coprime to $3$.  Hence we know that $3$ is in the image of the Swan Map.     

\bigskip
\lem{$2^m$ divides $3^{2^{m-3}}-1+2^{m-1}$ for $m \geq 4$.}

\hspace{2mm} Proof:  We proceed by induction.  $3^{2^{4-3}}-1+2^{4-1} = 16$.  So the proposition holds for $m=4$. 
Now suppose it holds for some $m$.  Then $2^mz =3^{2^{m-3}}-1+2^{m-1}$ for some $z$.  Rearranging gives

$$
3^{2^{m+1-3}}= (3^{2^{m-3}})^2 = (2^mz+1-2^{m-1})^2
$$

So

$$
3^{2^{m+1-3}}-1+2^{m+1-1} = (2^mz+1-2^{m-1})^2 -1 + 2^m
$$

$$
= 2^{2m}z^2+2^{2m-2} +2^{m+1}z-2^{2m}z = 2^{m+1}(2^{m-1}(z^2-z) + 2^{m-3} + z)
$$

So the proposition holds for $m+1$.  Hence by induction it holds for all $m$.
\hfill $\Box$

\bigskip
\bigskip
\prop{The elements $3$, $-3$ generate $(\Z/{2^m})^*$ for $m \geq 2$.}

\hspace{2mm} Proof:  The order of $(\Z/{2^m})^*$ is $2^{m-1}$.  $(\Z/{4})^* = \{1,3\}$ and $(\Z/{8})^* = \{1, -1,3,
-3\}$, so only the case $m \geq 4$ remains.  We know that the order of $3$ in $(\Z/{2^m})^*$ is a power of $2$.  The
previous lemma shows us that for $m \geq 4$ it is at least $2^{m-2}$, as 

\bigskip
$3^{2^{m-3}} \equiv 1+2^{m-1}$ Mod $2^m$.  

\bigskip
\hspace{2mm}It remains to show that $-1$ is not a power of $3$, as
then the $\pm 3^k$ give us all $2^{m-1}$ elements of $(\Z/{2^m})^*$.

\hspace{2mm} Suppose $3^k= -1 \quad {\rm Mod} \,\, 2^m$ for some $m \geq 4$.  Then $3^k = -1 \quad {\rm Mod} \,\, 8$ which is impossible as $3^k$ only
takes the values $1$ and $3$ modulo $8$.

\hfill $\Box$

\hspace{2mm}Combining this result with corollary 2.2.7 we obtain

\bigskip
\cor{The Swan Map $Aut(J) \to (\Z_{2^m})^*$ is surjective for all $m \geq 2$.}

\hspace{2mm}From theorem 1.2.8 we may conclude

\thm{An algebraic $2$- complex, $X$, over $\Z[D_{2^m}]$, with $\pi_2(X) \cong J \oplus \ZG^k$, is uniquely determined up to chain homotopy equivalence.}

\bigskip
\bigskip
\bigskip
\sec{The D(2) property for $\Zd$} 

\hspace{2mm}  Let $F_2$ denote the two element module over $\ZD$, on which the
action of $\ZD$ is trivial.

\prop{(i)$H^0(D_{4n},F_2)=F_2$
\newline (ii)$H^1(D_{4n},F_2)=F_2^2$
\newline (iii)$H^2(D_{4n},F_2)=F_2^3$}

\proof{We have the following resolution for $\Z$ over $\ZD$: \newline
\xymatrix{
\cdots\ar[r]&C_3 \ar[r]^{\partial_3} 
&C_2 \ar[r]^{\partial_2} 
&C_1 \ar[r]^{\partial_1} 
&C_0 \ar[r]&\Z\ar[r]&0}  \bigskip\newline
where $C_0$ is the free module 
generated by $*$;  $\,C_1$ is the free module generated by $e_1, e_2$;   
$\, C_2$ is the free module generated by $E_1, E_2, E_3$, and 
$C_3$ is the free module generated by $D_1, D_2, D_3, D_4$. 
\newline
\hspace{2mm} The maps $\partial_1,\partial_2,\partial_3$ are given by
\newline
$\partial_1 e_1=*(a-1)$ \newline
$\partial_1 e_2=*(b-1)$ \newline
\newline
$\partial_2 E_1=e_1\Sigma$ \newline
$\partial_2 E_2=e_2(b+1)$ \newline
$\partial_2 E_3=e_1(1+ba)+e_2(a-1)$ \newline
\newline
$\partial_3 D_1=E_1(a-1)$ \newline
$\partial_3 D_2=E_2(b-1)$ \newline
$\partial_3 D_3=E_1(b+1)-E_3\Sigma$ \newline
$\partial_3 D_4=E_2(a-1)-E_3(1-ba)$ \newline
\newline
${}$\hspace{2mm}  As all the coefficients above have even augmentation, applying
${\rm Hom}_{Z[D_{4n}]} (\bullet,F_2)$ to this resolution gives:
\newline
\xymatrix{
\cdots&F_2^4 \ar[l]
&F_2^3 \ar[l]_0 
&F_2^2 \ar[l]_0 
&F_2 \ar[l]_0} 
\newline
from which we deduce (i), (ii) and (iii) immediately.
} 
\bigskip

\hspace{2mm}Recall the sequence(1), from \S2.2.

\prop{$J$ has minimal $Z$- rank in its stable class.}

\hspace{2mm} Proof:  Given any finite algebraic 2-complex, consider the cochain obtained by applying ${\rm Hom}_{Z[D_{4n}]} (\bullet
,F_2)$:

$$
F_2^{d_2} \stackrel{v_2}{\leftarrow} F_2^{d_1} \stackrel{v_1}{\leftarrow} F_2^{d_0} 
$$

where ${d_0}$, ${d_1}$, ${d_2}$, are the $\ZD$ ranks of the zeroth, first and second modules
in the complex.  As $H^0(D_{4n},F_2)=F_2$, the kernel
of $v_1$ has $F_2$- rank $1$.  Consequently, the image of  $v_1$ has $F_2$- rank
$d_0-1$.  $H^1(D_{4n},F_2)=F_2^2$ so $v_2$ has kernel of $F_2$- 
rank $2+d_0-1=d_0+1$.  The image of $v_2$ is then seen to have rank $d_1-d_0-1$.  
$H^2(D_{4n},F_2)=F_2^3$ so we know that $d_2 \geq 3+ d_1- d_0 -1$. 
Rearranging gives $d_2-d_1+d_0 \geq 2$.

\hspace{2mm} Exactness implies that the $\Z$- rank of the algebraic $\pi_2$ of the algebraic complex must be $4n(d_2-d_1+d_0)-1$. Hence our inequality implies that this is at 
least $8n-1$, which is the $\Z$- rank of $J$.  

\hfill $\Box$

\hspace{2mm}We now restrict to the case $n=2$.

\bigskip
\prop{The only elements in the stable class of $J$ are modules
of the form $J \oplus \Zd^k$.}

\hspace{2mm}  Proof:  We refer to [11], Theorem 6.1.  This states that over $\Zd$, $A \oplus C =B \oplus C$ implies $A=B$ for torsion
free, finitely generated modules $A$, $B$, $C$.  

\hspace{2mm} If a module $M$ is in the stable class of $J$ then $M \oplus \Zd^r=J \oplus \Zd^s$.  From 
proposition 2.3.2 we have $s \geq r$.  From the theorem, we deduce that $M= J \oplus \Zd^{s-r}$.

\hfill $\Box$

\bigskip
\hspace{2mm} The only modules that can turn up as the second homology group of an algebraic 2-complex over $\Zd$ are ones of the form $J \oplus
\Zd^k$ for some $k \geq 0$.  From corollary 2.2.10 we know that for these modules, the Swan map is surjective. Hence theorem 1.2.8 implies that for 
each $k$, up to chain homotopy, there is a unique algebraic 2-complex with second homology group $J \oplus
\Zd^k$.  Given any $k$, the
homotopy class of this algebraic 2-complex is realized by the Cayley complex of the presentation

$$
<a,b \,\, |a^{2n}=b^2=e,\,\,aba =b, \,\, r_1=e,\, r_2=e, \dots r_k=e>
$$

where $r_i = e$ for $i=1, \dots ,k$.  We conclude

\bigskip
\thm{The group $D_8$ satisfies the D(2) property.}

\bigskip
\bigskip
\hspace{2mm}  In the next section we examine minimal elements of the stable class
$\Omega_3(\Z)$ for general groups $D_{4n}$, where we cannot immediately apply the 
torsion free cancellation which we used to prove proposition 2.3.3.

\sec{Minimal elements of $\Omega_3(\Z)$}

\hspace{2mm}  As before, we work over the group ring $\ZD$.  Consider the
resolution (1), in \S2.2:

$$
J\hookrightarrow \ZD^3 \stackrel{\partial_2}{\longrightarrow} \ZD^2 \stackrel{\partial_1}{\longrightarrow} \ZD \stackrel{\epsilon}{\twoheadrightarrow} \Z
\eqno (1)
$$

\hspace{2mm}  From the proposition 3.3.2 we know that $J$ has minimal $\Z$-
rank in its stable class.  From the Swan - Jacobinski Theorem, we know that the
only non-minimal modules in the stable class of $J$ are ones of the form $J
\oplus \ZD^k$.  It remains to investigate those modules, stably equivalent to
$J$ and of the same $\Z$- rank.  

\hspace{2mm}In this section, we take an arbitrary minimal element of $\Omega_3(\Z)$
and show that it must be isomorphic to one of a finite set of modules, which we parametrize 
by the group ${\Z_{2n}}^7 \oplus {\Z_2}^2$.

\hspace{2mm}Let $W_2$ denote the image of $\partial_2$.  We have a short exact sequence:

\bigskip
$J \to \ZD^3 \to W_2$

\bigskip
\hspace{2mm} Suppose $K$ is stably equivalent to $J$ and of the same $\Z$- rank.  Then we
have a short exact sequence:

\bigskip
$K\oplus F \to \ZD^3 \oplus F' \to W_2$

\bigskip
where $F$ and $F'$ are free module of the same rank. $(\ZD^3 \oplus F')/F$ is
stably free and of the same $\Z$- rank as $\ZD^3$.  In fact it is $\ZD^3$, as
dihedral groups satisfy the Eichler condition, hence stably free modules over
them are free.

\bigskip
\hspace{2mm}  Hence we have the following short exact sequence:

\bigskip
$K \stackrel{i}{\longrightarrow} \ZD^3  \stackrel{j}{\longrightarrow} W_2$

\bigskip
\hspace{2mm}  As $W_2$ is the image of $\partial_2$, it is the submodule of
$\ZD^2$, generated by 

$
e_1 \Sigma \newline
e_2 (1+b) \newline
e_1(1+ba)+e_2(a-1)
$

\bigskip
where, as before, $e_1, e_2$ are a basis of $\ZD^2$.

\bigskip
\hspace{2mm}  Define $p:W_2 \to \ZD$ by projection onto the factor generated by
$e_2$.  Let $W_1^T$ denote the image of the map.  Clearly, $W_1^T$ is generated
by $e_2(1+b)$ and $e_2(a-1)$.

\lem{The kernel of $p$ is generated by $e_1 \Sigma$.}
 
\proof{Clearly $e_1 \Sigma$ is in the kernel of $p$.  Its span has $\Z$- rank 2.
$W_2$ has $\Z$- rank $4n+1$ by exactness of (1).  Also $W_1^T$ has $\Z$- rank
$4n-1$, as it is generated by $e_2(1+b)$ and $e_2(a-1)$.  Hence the  $\Z$- rank
of the kernel of $p$ is 2.  
\newline
${}$\hspace{2mm} Therefore, given any element, $\alpha$ in the kernel of $p$, there 
exists some integer, r, such that $\alpha r$ is in the span of $e_1 \Sigma$. 
Hence $\alpha$ itself must be in the span of $e_1 \Sigma$.
}
 
\hspace{2mm}  The kernel of $p$ will be denoted by by $ZC_2$.

\hspace{2mm}  Let $f_1,\,f_2$
be a basis for a module isomorphic to $\ZD^2$.  Let $W_2^T$ be the submodule
generated by 

$f_1 \Sigma \newline
f_2 (1-b) \newline
f_1 (1-ba) -f_2(a-1)$

\bigskip
Let $s$ denote the natural inclusion of $W_2^T$ in $\ZD^2$, and let $t:\ZD^2
\to W_1^T$ be the map which sends $f_1$ to $e_2(a-1)$ and $f_2$ to $e_2 (1+b)$. 

\lem{We have a short exact sequence: \newline $$W_2^T
\stackrel{s}{\longrightarrow} \ZD^2 \stackrel{t}{\longrightarrow} W_1^T$$.}

\bigskip
\hspace{2mm} Let $M$ denote the kernel of $p \circ j$.

\lem{$M$ is isomorphic to $W_2^T \oplus \ZD$}.

\proof{We have short exact sequences:
\newline
$
M \to \ZD^3 \stackrel{p \circ j}{\longrightarrow} W_1^T
$
\newline
and 
\newline
$
W_2^T \oplus \ZD \stackrel{s \oplus 1}{\longrightarrow} \ZD^3 \stackrel{t \oplus 0}{\longrightarrow} W_1^T
$
\newline
${}$\hspace{2mm} Hence $M$ is stably equivalent to $W_2^T \oplus \ZD$, by
Schanuel's Lemma.  Dihedral groups satisfy the Eichler condition and clearly 
$W_2^T \oplus \ZD$ contains a free copy of the group ring as a summand, so by
the Swan-Jacobinski theorem, $M$ must be isomorphic to $W_2^T \oplus \ZD$.}

\hspace{2mm}  $K$ includes in $M$ and the cokernel of this inclusion is the
kernel of $p$, $ZC_2$.  As $K$ was chosen arbitrarily, we may conclude that any
module, stably equivalent to $J$ and of the same $\Z$- rank, occurs as the
kernel of some surjection $M \to ZC_2$.

\bigskip
\hspace{2mm}  As the action of $a$ on $ZC_2$ is trivial, the kernel of any map 
$M \to ZC_2$, must contain $M(a-1)$.  

\bigskip
\hspace{2mm} Let $\Z^T$ denote the $\ZD$ module whose underlying Abelian group
is isomorphic to the integers, and on which $a$ acts trivially and $b$ acts as
multiplication by $-1$. 

\lem{$W_2^T/W_2^T(a-1)$ has $\Z$- rank 3.}

\proof{$(W_2^T/W_2^T(a-1)) \otimes \Q = (W_2 ^T \otimes \Q)/((W_2^T \otimes \Q)
(a-1)$.  It is sufficient to show that this has dimension 3 over $\Q$.
\bigskip
\newline
${}$\hspace{2mm}
We have an exact sequence:
$$0 \to W_2^T \stackrel{s}{\longrightarrow} \ZD^2 \stackrel{t}{\longrightarrow}
 \ZD \to \Z^T \to 0$$
 ${}$ \hspace{2mm}Tensoring this sequence with $\Q$ and 
 applying "Whitehead's trick", yields
 $W_2^T \otimes \Q \oplus \QD= \QD^2 \oplus \Q^T$, where $\Q^T$ is the $\Q$-
 rank - 1 module on which $a$ acts trivially and $b$ acts as multiplication by
 -1.  Hence $W_2^T \otimes \Q = \QD \oplus \Q^T$. 
\bigskip
 \newline
${}$ \hspace{2mm} So $(W_2 ^T \otimes \Q)/((W_2^T \otimes \Q)= \QD/ \QD(a-1) \oplus
\Q^T/Q^T(a-1) = ZC_2 \otimes \Q \oplus \Q^T$ which has $\Q$ rank 3.
 }

\lem{The cokernel of the natural inclusion 
\newline $M(a-1) \hookrightarrow M$ \newline is $\Z
\oplus \Z^T \oplus \Z^T \oplus ZC_2$}

\proof{M is isomorphic to $W_2^T \oplus \ZD$.  We first observe that 
\newline $\ZD / \ZD(a-1) \cong 
ZC_2$.
\newline
${}$\hspace{2mm}  $W_2^T$ is generated by \newline
$f_1 \Sigma \newline
f_2 (1-b) \newline
f_1 (1-ba) -f_2(a-1)$
\newline
${}$\hspace{2mm}  Let $w_1$, $w_2$, $w_3$ denote the images of $f_1 \Sigma$,
$f_2(1-b)$, $f_1 (1-ba) -f_2(a-1)$ respectively, under quotienting by
$W_2^T(a-1)$.  $w_1$ generates a copy of $ZC_2$ and $w_2$ generates a copy of
$\Z^T$.
\bigskip
\newline
${}$\hspace{2mm}  Note that $w_3 2n=w_3 \Sigma = w_1(1-b)$.  Hence $
w_1b=(w_1 - w_3n) - w_3n$.  Also $w_1 = (w_1-w_3n) + w_3n$.  So
$W_2^T/W_2^T(a-1)$ is generated over $\Z$ by $w_2$, $w_3$, $w_3b$ 
and $w_1- w_3n$. 
\newline
${}$\hspace{2mm}
We will show that $w_3b$=$-w_3$ and hence that $W_2^T/W_2^T(a-1)$ is generated 
over $\Z$ by $w_2$, $w_3$,  
and $w_1- w_3n$:
\newline
$w_3(1+b)$ is equal to the image of $f_1(1-ba)(1+b)-f_2(a-1)(1+b)$ under
quotienting by $W_2^T(a-1)$.  But
\newline
$f_1(1-ba)(1+b)-f_2(a-1)(1+b)=(f_2(b-1) -(f_1(1-ba)-f_2(a-1))b)(a-1)       $
\newline
which is in $W_2^T(a-1)$.  Hence $w_3(1+b)=0$.
\newline
${}$\hspace{2mm}
The elements $w_2$, $w_3$,  
and $w_1- w_3n$ must be
$\Z$- linearly independent, in order for their span to have $\Z$- rank 3.  
\bigskip
\newline
${}$\hspace{2mm} Hence we know that $W_2^T/W_2^T(a-1)$ is 
torsion free.  We know that $(w_3)(1+b)2n=0$ and $(w_1-w_3n)(1-b)2=0$, and 
are able to conclude that $w_3b=-w_3$ and $(w_1-w_3n)b=(w_1-w_3n)$.    
\bigskip
\newline
${}$\hspace{2mm}
Hence $w_1- w_3n$, $w_2$ and $w_3$ 
generate $\Z\oplus \Z^T \oplus \Z^T$ and $M/M(a-1) \cong \Z
\oplus \Z^T \oplus \Z^T \oplus ZC_2$.
}

\hspace{2mm}  Any surjection $M \to ZC_2$ must therefore factor through 
$\Z \oplus \Z^T \oplus \Z^T \oplus ZC_2$. Hence we have a surjection

$\phi:\Z \oplus \Z^T \oplus \Z^T \oplus ZC_2 \to ZC_2$

and the cokernel of the natural inclusion of $M(a-1)$ in $K$, is naturally
identified with the kernel of $\phi$. 

\lem{The kernel of $\phi$ is isomorphic to $\Z \oplus \Z^T \oplus \Z^T$.}

\proof{Let $u_1$, $u_2$, $u_3$, $u_4$ generate copies of $\Z$, $\Z^T$, $\Z^T$
and $ZC_2$, respectively in the direct sum 
 $\Z \oplus \Z^T \oplus \Z^T \oplus ZC_2$.  $v$ will generate the image of
 $\phi$.
\bigskip
\newline
${}$\hspace{2mm} $\phi(u_1)$ is an integer multiple of $v(1+b)$ and $\phi(u_2)$,
 $\phi(u_3)$ are
integer multiples of $v(1-b)$.  Consequently $\phi(u_4)=v(x+yb)$, where $x$ and
$y$ are integers with odd sum, as otherwise $\phi$ would not be surjective.
\bigskip
\newline    
${}$\hspace{2mm}  If some linear combination of the $u_i$ is in the kernel of
$\phi$, the coefficient on $u_4$ must have even augmentation, as the augmentation
multiplied by $x+y$
must be even.  Hence the kernel
of $\phi$ is contained in the $\Z$- linear span of $u_1$, $u_2$, $u_3$,
$u_4(1+b)$, $u_4(1-b$.   Also the image of this span under $\phi$ must be the
whole of the span of $v(1+b)$, $v(1-b)$, as $x+y$ is odd.
\bigskip
\newline    
${}$\hspace{2mm} $\phi$ restricts to a surjection from the span of $u_1$, $u_4(1+b)$ to
the span of $v(1+b)$.  The kernel of a surjection $\Z^2 \to Z$ must be
isomorphic to $\Z$.
\bigskip
\newline    
${}$\hspace{2mm} $\phi$ restricts to a surjection from the span of $u_2$, $u_3$,
$u_4(1-b)$ to
the span of $v(1-b)$.  The kernel of a surjection $(\Z^T)^3 \to Z^T$ must be
isomorphic to $\Z^T \oplus \Z^T$.
\bigskip
\newline
${}$\hspace{2mm}  Hence the kernel of $\phi$ is $\Z \oplus \Z^T \oplus \Z^T$
}

\hspace{2mm} Hence we know that $K$ occurs in a short exact sequence of the
form:

$$
M(a-1)\to K \to \Z \oplus \Z^T \oplus \Z^T
$$

\bigskip
\hspace{2mm}  Hence all candidates for minimal
elements of $\Omega_3(\Z)$, are parametrized by the group 

${\rm Ext}^1_{\ZD}( \Z \oplus \Z^T \oplus \Z^T, \,M(a-1))$

\bigskip
\hspace{2mm}  In order to calculate this extension group, we will require a
better description of $W_2^T(a-1)$.  We know that $W_2^T$ is generated by

\bigskip
$f_1 \Sigma \newline
f_2 (1-b) \newline
f_1 (1-ba) -f_2(a-1)$

\lem{If $f_1x \in W_2^T$ then $\Sigma$ divides $x$.}

\proof{Suppose $f_1x = f_1 \Sigma \mu_1 + f_2(1-b)\mu_2+(f_1 (1-ba)
-f_2(a-1))\mu_3$.  We could replace$\mu_2$ with a polynomial in $a$, which must
be divisible by (a-1).  We will denote the polynomial $(a-1)p$ where $p$ is a
polynomial in $a$. So
\newline
$f_1x = f_1 \Sigma \mu_1 + f_2(1-b)(a-1)p+(f_1 (1-ba)
-f_2(a-1))\mu_3=
\newline
= f_1 \Sigma \mu_1 + f_2(a-1)(1+ba)p+(f_1 (1-ba)
-f_2(a-1))\mu_3
\newline
= f_1 \Sigma \mu_1  +(f_1 (1-ba)
-f_2(a-1))(\mu_3-(1+ba)p)
$
\newline
As $(a-1))(\mu_3-(1+ba)p)=0$, we know that $(\mu_3-(1+ba)p)=\Sigma y$, for some
element $y$ of $\ZD$.  So $x= \Sigma \mu_1 + \Sigma(1-b)y$.
}

\hspace{2mm} Let $W_1$ denote the augmentation ideal of
$\ZD$.

\lem{$W_2^T(a-1)\cong W_1(a-1)$}

\proof{We have a surjective homomorphism $W_2^T(a-1)\twoheadrightarrow W_1(a-1)$
given by projection onto the $f_2$ component.  It is sufficient to show that it
has $0$ kernel.
\newline
${}$\hspace{2mm}  Suppose $f_1x$ is in $W_2^T(a-1)$, for some $x \in \ZD$. 
Clearly $(a-1)$ divides $x$.  By the previous lemma, $\Sigma$ also divides $x$. 
Hence $x=0$.}

\hspace{2mm}  Recall we defined $ZC_2$ to be the quotient $\ZD/\ZD(a-1)$.  Let
$L$ denote the quotient $\ZD/ID(a-1)$.

\bigskip

${\rm Ext}^1_{\ZD}( \Z \oplus \Z^T \oplus \Z^T, \,M(a-1))$

\bigskip
$\cong{\rm Hom}_{\rm Der}(W_1 \oplus W_1^T \oplus W_1^T, \, W_2^T(a-1) \oplus \ZD(a-1))$ 

\hspace{2mm}  By dimension shifting, we have that this is equal to 

\bigskip

${\rm Hom}_{\rm Der}(\Z \oplus \Z^T \oplus \Z^T, \, L \oplus ZC_2)$.  

\bigskip

$\cong{\rm Hom}_{\rm Der}(\Z, ZC_2) \oplus 
{\rm Hom}_{\rm Der}(\Z^T,\,ZC_2 )^2
\oplus {\rm Hom}_{\rm Der} (\Z,\, L)  \oplus 
{\rm Hom}_{\rm Der} (\Z^T ,\,L)^2$

\bigskip
\hspace{2mm}  The module $ZC_2$ is generated freely 
over $\Z$ by the images of $1$ and
$b$.  $1b=b$ and $bb=1$.  The action of $a$ is trivial.  A map from $\Z$ to
$ZC_2$ is determined by where $1\in\Z$ is sent to.  It must go to an element
which is fixed by the action of $b$.  Hence it must go to some multiple of the
image of $(1+b)$.  

\bigskip
\hspace{2mm}  Any map from $\Z$ to $ZC_2$ which factors through a projective
will send $1 \in \Z$ to $x\Sigma(1+b)$ for some element $x \in ZC_2$.  Thus $1
\in \Z$ will go to some multiple of $(1+b)2n$.  So 
${\rm Hom}_{\rm Der}(\Z, ZC_2) \cong \Z_{2n}$.

\bigskip
\hspace{2mm} Similarly, a map from $\Z^T$ to
$ZC_2$ is determined by where $1\in\Z^T$ is sent to.  It must go to an element
which is sent to its negative by the action of $b$.  Hence it must 
go to some multiple of the
image of $(1-b)$.  

\bigskip
\hspace{2mm}  Any map from $\Z^T$ to $ZC_2$ which factors through a projective
will send $1 \in \Z^T$ to $x\Sigma(1-b)$ for some element $x \in ZC_2$.  Thus $1
\in \Z^T$ will go to some multiple of $(1-b)2n$.  So 
${\rm Hom}_{\rm Der}(\Z^T, ZC_2) \cong \Z_{2n}$.

\bigskip
\hspace{2mm} Given any element of $\ZD$, by subtracting an appropriate integral
combination of $1$ and $b$ we have a multiple of $(a-1)$.  By further 
subtracting an appropriate integral multiple of $(a-1)$, we have an element of
$W^1(a-1)$.  Hence $L$ is generated over $\Z$ by the images of $1$, $b$, 
and $(a-1)$,
denoted hereafter by $u_1$, $u_2$ and $u_3$.   Note $u_3 2n$ is equal to the
image of  $(a-1)\Sigma + (2n-\Sigma)(a-1)$ which equals the image of
$(2n-\Sigma)(a-1)$  which is $0$.  So $u_32n=0$.

\bigskip 
\hspace{2mm}  Now suppose $u_1r+u_2s+u_3t=0$ for integers $r,s,t$. 
Clearly $r=0$ and $s=0$.  Also $(a-1)t=x(a-1)$ for some $x \in W^1$.  Hence $t$
differs from something of augmentation $0$ by some multiple of $\Sigma$.  Hence
$2n$ divides $t$.  So as Abelian group, $L$ is isomorphic to $\Z^2 \oplus
\Z_{2n}$.

\bigskip 
\hspace{2mm} 
Next we calculate the group action on $L$. 
 
$a=1 +(a-1)$

$ba =b +b(a-1)=b+(a-1)+(b-1)(a-1)$

$(a-1)a = (a-1) + (a-1)(a-1)$

So

$(u_1r+u_2s+u_3t)a = u_1r+u_2s+u_3(r+s+t)$

Also

$1b=b$

$bb=1$

$(a-1)b= -ba^{-1}(a-1) = (1-ba^{-1})(a-1) -(a-1)$

So

$(u_1r+u_2s+u_3t)b = u_1s+u_2r-u_3t$

\bigskip
\hspace{2mm}  A map from $\Z$ to $L$ is determined by where $1 \in \Z$ goes.  It
must go to an element which is fixed by the action of $b$.  Hence it goes to an
element of the form $(u_1r+u_2r+u_3nt)$.  This element must also be fixed by
the action of $a$.  Hence it has the form $(u_1ns+u_2ns+u_3nt)$.

\bigskip
\hspace{2mm}  Any such map which factors through a 
projective must send $1 \in \Z$ to
an element of the form

$(u_1r+u_2s+u_3t)(1+b)\Sigma=(u_1(r+s)+u_2(r+s))\Sigma=
\newline
u_1(r+s)2n+u_2(r+s)2n +
u_3(2(r+s)(2n-1)n= u_1(r+s)2n+u_2(r+s)2n$

\bigskip
\hspace{2mm}  Hence ${\rm Hom}_{\rm Der}(\Z, L) \cong {\Z_2}^2$.

\hspace{2mm}  A map from $\Z^T$ to $L$ is determined by 
where $1 \in \Z^T$ goes.  It
must go to an element which is sent to its negative by fixed by the action of 
$b$.  Hence it goes to an
element of the form $(u_1s - u_2s+u_3t)$.  This element will automatically 
be fixed by
the action of $a$.  

\bigskip
\hspace{2mm}  Any such map which factors through a 
projective must send $1 \in \Z^T$ to
an element of the form

$(u_1r+u_2s+u_3t)(1-b)\Sigma=(u_1(r-s)+u_2(s-r)+u_32t)\Sigma=
\newline
u_1(r-s)2n+u_2(s-r)2n +
u_3 4tn= u_1(r-s)2n+u_2(s-r)2n$

\bigskip
\hspace{2mm}  Hence ${\rm Hom}_{\rm Der}(\Z^T, L) \cong {\Z_{2n}}^2$.

\bigskip
\hspace{2mm} Hence the extension group is

${\rm Hom}_{\rm Der}(\Z, ZC_2) \oplus 
{\rm Hom}_{\rm Der}(\Z^T,\,ZC_2 )^2
\oplus {\rm Hom}_{\rm Der} (\Z,\, L)  \oplus 
{\rm Hom}_{\rm Der} (\Z^T ,\,L)^2$

\bigskip 
$\cong \Z_{2n} \oplus {\Z_{2n}}^2 \oplus {\Z_2}^2 \oplus {\Z_{2n}}^4$ 
$\cong {\Z_{2n}}^7 \oplus {\Z_2}^2$

\bigskip
\hspace{2mm}  In conclusion, we have found an exhaustive list 
of candidates for minimal modules in the stable class $\Omega_3(\Z)$.  These candidates 
are parametrized by the elements of the finite group, ${\Z_{2n}}^7 \oplus {\Z_2}^2$.  
This gives us an upper bound of $512n^7$ for the number of minimal modules in $\Omega_3(\Z)$. 

\bigskip
\hspace{2mm}
Not
all of the candidates will be in $\Omega_3(\Z)$.  If one could show that those
which are in $\Omega_3(\Z)$, are all isomorphic, then one would have proved that
cancellation of free modules holds in $\Omega_3(\Z)$ for dihedral groups of
order $4n$.

\sec{$D_{8n+4}$}

\hspace{2mm} The dihedral groups $D_{4n+2}$ have balanced presentations and period 4 resolutions over $\Z$.  This makes them easier to work with in 
certain respects than the dihedral groups $D_{4n}$, which have neither.  We show in this section that the groups $D_{8n+4}$ have period 4
resolutions if one works over the ground ring $\Z[x]/\!\!<2x-1\!\!>$, which we denote $\Zh$.  This raises the possibility that the methods used
to prove the D(2) property for the groups $D_{4n+2}$ (see \cite{John} and \cite{John1}) may be generalized to the groups $D_{8n+4}$.

\bigskip
\lem{\it $D_{8n+4} \cong D_{4n+2} \times C_2$}

\hspace{2mm} Proof: Take $\D = <a,b| \,aba=a, \, b^2=a^{4n+2}=e>$.  Then $a^{2n+1}$ is central and generates a copy of $C_2$.  The 
elements $a^2$, $b$ generate a copy of $D_{4n+2}$. As $2n+1$ is odd, any element can be uniquely written as a product of an element in the 
copy of $C_2$ and an element in the copy of $D_{4n+2}$.

\hfill $\Box$

\bigskip
\hspace{2mm} We work over the 
ring, ${\mathbb {Z}[\frac {1}{2}][D_{8n+4}]}$, which 
we will denote $\ZhD$. Let 
$$\nabla_1 = (1+a^{2n+1})/2, \qquad \nabla_2=(1-a^{2n+1})/2$$.  

\hspace{2mm}Then $\nabla_1, \, \nabla_2$ are central idempotents, 
which sum to $1$.  Consequently, we have 
$$
\ZhD \cong \ZhD \nabla_1 \times \ZhD \nabla_2
$$ as rings.   Let $T$ denote the ring $R\nabla_1$ and $S$ denote the ring $R\nabla_2$.

\bigskip
\hspace{2mm} Regarding $a^2$ and $b$ as generators of the subgroup $D_{4n+2}$, as before, we have an isomorphism 
$T\cong \mathbb {Z}[\frac {1}{2}][D_{4n+2}]$.   From \cite{John1}, \$41(ii) 
we obtain a free period 4 resolution of $\Z$ over $\Z[D_{4n+2}]$. This naturally yields a period 4 resolution of $\Zh$ over $T$.  We write 
this resolution in our notation, regarding $T$ as a subring of $R$:

$$
0 \rightarrow \Zh \stackrel{\epsilon^*}{\longrightarrow} \ZhDf \stackrel{\partial_3}{\longrightarrow} {\ZhDf}^2 \stackrel{\partial_2}{\longrightarrow} {\ZhDf}^2 \stackrel{\partial_1}{\longrightarrow} \ZhDf \stackrel{\epsilon}{ \longrightarrow} \Zh \rightarrow 0
\eqno(1)
$$

\bigskip
\hspace{2mm}Going from left to right, let the basis elements of the free modules in this sequence be $\{c\}, \, \,\{f_1, f_2 \}, \, \, \{e_1, e_2\}$ and $\{v\}$.
$\epsilon$ is the augmentation map, which takes $v \in \ZhDf$ to $1 \in \Zh$ and $\epsilon^*$ is its dual.  

\bigskip
\hspace{2mm} 
For any integer $k$ let $\Sigma_k = \sum_{i=0}^{k-1}a^{2i}$.  With these conventions we set:

\begin {eqnarray*}
\partial_3 c = f_1 \nabla_1(1+a^2-a^{2(n+1)} - b) + f_2\nabla_1(a^2-a^{2n}b)\\ \\
\partial_2 f_1 = e_1 \nabla_1\Sigma_{2n+1}-e_2\nabla_1(1+b) \,\,\,\,\,\,\,\,\, \quad \qquad \quad \quad\quad\quad\\
\partial_2 f_2 = e_1 \nabla_1(\Sigma_n(b-1)+a^{2n}b)+e_2\nabla_1(1-a^{2n})\quad\quad\,\\ \\
\partial_1 e_1 = v \nabla_1(a^2-1) \qquad \qquad \qquad \qquad\quad \quad \, \,\quad\quad\,\,\,\quad\,\,\,\,\\
\partial_1 e_2 = v \nabla_1(b-1) \qquad \qquad \qquad \qquad\quad \quad \,\,\,\,\quad\quad\quad\quad\,\,
\end {eqnarray*}

\hspace{2mm}We also have an exact sequence over $S$:

$$
0 \longrightarrow 0 \longrightarrow S \stackrel{p_3}{\longrightarrow} S^2 \stackrel{p_2}
{\longrightarrow}
 S^2 \stackrel{p_1
 }{\longrightarrow} S {\longrightarrow} 0 \longrightarrow 0
$$

\bigskip
\hspace{2mm}Again, going from left to right, let the basis elements of the free modules in this sequence be $\{c'\}, \, \,\{f_1', f_2' \},
 \, \, \{e_1', e_2'\}$ and $\{v'\}$.

\hspace{2mm} With respect to this basis we set $$
p_3={\tiny
\left( \begin{array}{c}  0\\ 1\end{array} \right)}, \quad
p_2={\tiny
\left( \begin{array}{cc}  0&0\\ -1&0 \end{array} \right)}, \quad
p_1={\tiny \left( \begin{array}{cc}  1&0 \end{array} \right)}
$$

\bigskip
\hspace{2mm}We can take the direct product of these two resolutions, to get a resolution over the product ring, $\ZhD$: 

$$
0 \rightarrow \Zh \stackrel{\epsilon^* \times 0}{\longrightarrow} \ZhD \stackrel{\partial_3 \times p_3}{\longrightarrow} {\ZhD}^2 \stackrel{\partial_2 \times p_2}{\longrightarrow} {\ZhD}^2 \stackrel{\partial_1 \times p_1}{\longrightarrow} \ZhD \stackrel{\epsilon \times 0}{ \longrightarrow} \Zh \rightarrow 0
$$                   

\hspace{2mm}Let $\delta$ denote the augmentation map $R \to \Zh$ sending $1 \in R$ to $1 \in \Zh$.  Then  
$\delta(\nabla_1)=1$ and $\delta(\nabla_2)=0$.  Hence $\delta = \epsilon \times 0$.  Note also that 
$\delta^*(1)=2\nabla_1\epsilon(1)$.

\hspace{2mm} Hence, if we let $d_i = \partial_i \times p_i$  we get an exact sequence:
  
$$
0 \rightarrow \Zh \stackrel{\epsilon^*}{\longrightarrow} \ZhD \stackrel{d_3}{\longrightarrow} {\ZhD}^2 \stackrel{d_2}{\longrightarrow} {\ZhD}^2 \stackrel{d_1}{\longrightarrow} \ZhD \stackrel{\epsilon}{\longrightarrow} \Zh \rightarrow 0
\eqno(2)
$$     

\hspace{2mm} Let $c''=c + c',\,\, f_1''=f_1+f_1',\,\,f_2''= f_2+f_2',\,\,
e_1''=e_1+e_1', \,\,e_2''=e_2+e_2',\,\, v''=v+v'$. We have

\begin {eqnarray*}
d_3 c'' = f_1'' \nabla_1(1+a^2-a^{2(n+1)} - b) +f_2'' (\nabla_1(a^2-a^{2n}b)+ \nabla_2)\\\\ 
d_2 f_1'' = e_1'' \nabla_1 \Sigma_{2n+1}-e_2''(\nabla_1(1+b)+\nabla_2) \qquad \,\,\,\,\,\,\,\,\,\,\quad \quad \qquad \,\,\quad\\ 
d_2 f_2'' = e_1'' \nabla_1 (\Sigma_n(b-1)+a^{2n}b)+e_2''\nabla_1(1-a^{2n}) \,\,\, \,\,\,\quad \qquad\quad\quad\\ \\
d_1 e_1'' = v'' (\nabla_1(a^2-1)+\nabla_2) \,\,\,\,\quad \qquad \qquad \qquad \qquad\quad \,\,\quad\quad\quad \,\,\,\\
d_1 e_2'' = v'' \nabla_1(b-1)  \quad \qquad \qquad \qquad \qquad \qquad \qquad \qquad \quad\quad\quad\,\,
\end {eqnarray*}

\hspace{2mm}We have demonstrated a period 4 resolution over ${\mathbb {Z}[\frac {1}{2}][D_{8n+4}]}$.  We now consider how this relates to the
D(2) problem for $D_{8n+4}$.

\bigskip      
\hspace{2mm} Let $J$ denote the image of $\partial_3$.  Let $\iota$ denote its inclusion into ${\ZhDf}^2$.  The image 
of $d_3$ is equal to $J \times S$.  Let $\iota_S$ denote the the inclusion of $S$ into $S^2$, induced by $p_3$.

\hspace{2mm}We have a chain map:

\begin{eqnarray*}
S \stackrel{\iota_S}{\longrightarrow} S^2 \stackrel{p_2}{\longrightarrow} S^2 \stackrel{p_1}{\longrightarrow} S {\longrightarrow} 0 \,\,\,\,\\
\downarrow 1 \quad \,\, \downarrow 1 \quad \quad \downarrow 1 \,\,\, \quad\downarrow 1 \,\,\, \downarrow 0\\
S \stackrel{\iota_S}{\longrightarrow} S^2 \stackrel{p_2}{\longrightarrow} S^2 \stackrel{p_1}{\longrightarrow} S {\longrightarrow} 0\,\,\,\,\\
\end{eqnarray*} 

\hspace{2mm} So given $k \in \Zh$ if there exists an automorphism $\alpha$, and a chain map

\begin{eqnarray*}
J \stackrel{\iota}{\longrightarrow} {\ZhDf}^2 \stackrel{\partial_2}{\longrightarrow} {\ZhDf}^2 \stackrel{\partial_1}{\longrightarrow} \ZhDf \stackrel{\epsilon}{ \longrightarrow} \Zh \quad
\\
\downarrow \alpha  \, \quad \downarrow m_2 \quad  \downarrow m_1   \quad\downarrow m_0 \quad \,\, \downarrow \times k
\\
J \stackrel{\iota}{\longrightarrow} {\ZhDf}^2 \stackrel{\partial_2}{\longrightarrow} {\ZhDf}^2 \stackrel{\partial_1}{\longrightarrow} \ZhDf \stackrel{\epsilon}{ \longrightarrow} \Zh \quad 
\end {eqnarray*}

we have the product chain map

\begin{eqnarray*}
J \times S \quad\stackrel{\iota\times \iota_S}{\longrightarrow}\quad {\ZhD}^2 
\quad\stackrel{d_2}{\longrightarrow} \quad{\ZhD}^2 \quad\stackrel{d_1}{\longrightarrow} \quad 
\ZhD \quad\stackrel{\epsilon}{ \longrightarrow}\quad \Zh \quad
\\
{\downarrow \alpha \times 1 \quad \quad \downarrow  m_2 \times 1  \,\,\, \quad \downarrow {m_1 \times 1} \,\,\, \quad 
\downarrow m_0 \times 1 \,\, \quad \downarrow \times k}
\\
J \times S \quad\stackrel{\iota\times \iota_S}{\longrightarrow}\quad {\ZhD}^2 
\quad\stackrel{d_2}{\longrightarrow} \quad{\ZhD}^2 \quad\stackrel{d_1}{\longrightarrow} \quad 
\ZhD \quad\stackrel{\epsilon}{ \longrightarrow}\quad \Zh \quad
\end{eqnarray*} 

\hspace{2mm} However, even if the $m_i$ are chosen to be representable by
matrices with coefficients all in $\Z[D_{4n+2}] \subset T$, the maps $m_i \times 1$ may not be. Consequently, we have not completely reduced
the problem of realizing $k$- invariants via automorphisms over $D_{8n+4}$, to the technically simpler problem over $D_{4n+2}$ 
(see \cite{John} and \cite{John1}).  However, the above construction does provide a link between the two problems.

\newpage
\sec{Summary of results concerning the D(2) problem for dihedral groups}

\hspace{2mm}  In \cite{John1} (62.3) it is shown that the dihedral groups $D_{4n+2}$ satisfy the D(2) property.  The smallest dihedral group
not covered by this is $D_8$.  We show that
the D(2) property does hold for $D_8$ (theorem 2.3.4).  

\hspace{2mm}More generally, for dihedral groups of order $4n$, we show that a minimal 
element of $\Omega_3(\Z)$ is realized as the $\pi_2$ of a presentation (proposition
2.3.2).  We parametrize all possible minimal elements of $\Omega_3(\Z)$ by the group ${\Z_{2n}}^7 \oplus {\Z_2}^2$ (\S2.4).

\hspace{2mm}In the case of dihedral groups of order $2^n$, $n \in \Z$, we are further able to show that up to chain homotopy equivalence, there is a unique algebraic $2$-
complex with "standard" algebraic $\pi_2$ (theorem 2.2.11).

\chapter{$\pi_3$ of geometric 2- complexes}\label{}
\stepcounter{ch}
\setcounter{sec}{0}

\hspace{2mm}  Let $X$ be a finite geometric $2$- complex with finite fundamental 
group $G$. Lemma 2.1.5 implies
that $\pi_2(X) = H_2(C_*(\tilde{X});\Z)$ as modules over $\ZG$.  Hence 
$\pi_2(X)$ is determined by $C_*(\tilde{X})$.  

\hspace{2mm}In fact, $C_*(\tilde{X})$ determines the
homotopy type of $X$, (see \cite{John1}, theorem 49.5) and hence it determines 
all the
homotopy groups, as modules over $\ZG$.  In this chapter we compute $\pi_3(X)$
from $C_*(\tilde{X})$.  We also show that if $G$ has even order, then $G$
determines the stable class of the module $\pi_3(X)$.

\hspace{2mm}   Let $G$ be a finite group and let $J$ be 
a finitely generated, torsion free,
module over it.  There is a $G$- action on $J \otimes_{\Z}
J$, given by $(a \otimes b)g = ag \otimes bg$, making it a $\ZG$ module.  Let
$t$ be the $\ZG$- linear automorphism of $J \otimes_{\Z} J$ defined by 
$t(a \otimes b) =b \otimes a$.  

\bigskip
{\bf Definition} $S^2(J)=\{x \in J \otimes J|tx=x\}$

\hspace{2mm} The main theorem of this chapter is now stated.

{\bf Theorem 3A} {\it Let $X$ be a finite geometric $3$- complex, with finite 
fundamental group $G$.  If $C_*(\tilde{X}) \sim  \mathcal{A}$, for some 
finite algebraic $2$- complex $\mathcal{A}$, then  $\pi_3(X)\cong S^2(J)$, where 

$J=H_2(C_*(\tilde{X});\,\ZG)$.}

\hspace{2mm} In particular,

{\bf Corollary} {\it Let $X$ be any finite geometric $2$- complex, with 
finite fundamental group $G$.  Let $J \cong\pi_2(C_*(\tilde{X}))$.  Then  
$\pi_3(X)\cong S^2(J)$.}

\bigskip
\hspace{2mm}The claim that $\pi_3(X)$ may be computed directly from $J$ is non -
trivial in the sense that $J$ does not by itself determine the homotopy 
type of $C_*(\tilde{X})$.

\hspace{2mm} Let $\mathcal{A}$ be a finite algebraic $2$- complex.

\bigskip
{\bf Definition}
{$\pi_3(\mathcal{A})= S^2(J)$, where $J \cong \pi_2(\mathcal{A})$.}

\bigskip
\hspace{2mm}From the theorem, we see that this definition is 
consistent, in that if $C_*(\tilde{X}) \sim \mathcal{A}$, for a finite geometric $3$-
complex $X$, then $\pi_3(\mathcal{A})\cong\pi_3(X)$.

\sec{Higher Covers}

\bigskip
\hspace{2mm}  We computed $\pi_2$ of a geometric $2$- complex by noting that 
passing to the universal cover preserved it, and then applying the Hurewicz
isomorphism theorem to compute it from the homology of the algebraic $2$-
complex.  Our first step in computing $\pi_3$ is to generalize the notion 
of universal cover.

\bigskip
\hspace{2mm}We say a $X$ space is $r$- connected if $\pi_i(X)=0$, for
$i=0,1,...,r$.

\define{} For $r \geq 1$, an $r$- cover of an $r-1$- connected space $X$, is a
fibre bundle map $f:Y \to X$, with fibre $k(\pi_r(X), r-1)$ and with  

$\delta_r:\pi_r(X) \to \pi_{r-1}(k(\pi_r(X), r-1))$ an isomorphism. Here
$\delta_r$ is the boundary operator in the long exact sequence associated 
to the fibre bundle.

\lem{We have:}

$\pi_i(Y)=0$, $\qquad i \leq r$

$\pi_i(Y)=\pi_i(X)$, $\qquad i >r$

\proof{$X$ is $r-1$- connected, and $k(\pi_r(X), r-1)$ is an Eilenberg- Mac lane
space.  Consequently, the long exact sequence associated to the fibre 
bundle has the form:
\newline
{\tiny
\xymatrix{&&\cdots\ar[r]^{\quad\quad0}&\\
0\ar[r]&\pi_{r+2}(Y)\ar[r] & \pi_{r+2}(X) \ar[r]^{\quad\quad 0}&\\
0\ar[r]&\pi_{r+1}(Y)\ar[r] & \pi_{r+1}(X) \ar[r]^{\quad\quad 0}&\\
0\ar[r]&\pi_{r}(Y)\ar[r]^0 & \pi_{r}(X) \ar[r]^{\quad\quad\delta_{r}}&\\
\pi_r(X)\ar[r]^0&\pi_{r-1}(Y)\ar[r] & 0 \ar[r]^{\quad\quad0}&\\
0 \ar[r]^0&\pi_{r-2}(Y)\ar[r] & 0 \ar[r]^{\quad\quad0}&\\
\cdots&&&
}}
\newline
\newline
${}$\hspace{2mm}So $Y$ is $r$- connected and we have isomorphisms 
$\pi_i(Y) \to \pi_i(X)$ for $i>r$.
}

\hspace{2mm}  Note also that the universal cover of a connected space, together
with the associated covering map, are a $1$- cover.  In fact all  $1$- covers
have that form.

\bigskip
\hspace{2mm}  Suppose we have a connected space $X$ and a tower of maps

\bigskip
\xymatrix{\cdots\ar[r]^{f_3} &X_2 \ar[r]^{f_2} &X_1 \ar[r]^{f_1} &X}

\bigskip
where for each $i$, $f_i$ is an $i$- cover.  Then each $X_i$ is $i$-
connected and by induction on $r$, $\pi_s(X_r) \cong \pi_s(X)$ as long as $s>r$.

\bigskip
\hspace{2mm}  In particular, $\pi_{r+1}(X_r) \cong \pi_{r+1}(X)$, and $X_r$ is
$r$- connected, for $r \geq 1$, so by the Hurewicz isomorphism theorem, we 
may conclude: 

\lem{$\pi_{r+1}(X) \cong
H_{r+1}(X_r; \Z)$, $r \geq 1$}.  

\hspace{2mm}We will use this identity in the next section to
calculate $\pi_3$ of a geometric $2$- complex.

\bigskip
\hspace{2mm}  Note that the fibre of the map $f_r$ is 
$k(\pi_r(X_{r-1}),r-1)$, which is equal to 

$k(\pi_r(X),r-1)$.  The inclusion of
this fibre in $X_r$ induces a map on homology: 

$H_{r+1}(k(\pi_r(X),r-1);\Z) \to 
H_{r+1}(X_r; \Z)$.  

\hspace{2mm}We already have $\pi_{r+1}(X) \cong
H_{r+1}(X_r; \Z)$, so a map is induced
$$
H_{r+1}(k(\pi_r(X),r-1);\Z) \to \pi_{r+1}(X)
$$

for each $r \geq 1$.  Both the domain and codomain of this map are 
invariants of $X$.



\sec{The Hopf fibration}

\hspace{2mm}  Our first example of a $2$- cover will be the Hopf fibration. In this
section we define it and look at the map of sets it induces: $\pi_2(X) \to \pi_3(X)$. 
Crucially, for our main theorem, will show that this map respects the 
action of $\pi_1 (X)$.

\hspace{2mm}We
may regard $S^3$ as the set of pairs $(z,w) \in \C ^2$, satisfying $z \bar{z}+w
\bar{w}=1$. Define a relation $\sim$, by setting $(z,w) \sim (z \lambda,
w\lambda)$ whenever $\lambda \in \C$ satisfies $\lambda \bar{\lambda} =1$.

\define{Hopf fibration} The Hopf fibration, $h:S^3 \to S^2$ is the natural map 
$S^3 \to S^3/\sim$ composed with the topological identifications $S^3 /\sim \cong 
\C P^1 \cong S^2$.

\hspace{2mm}  This map is a fibre bundle map with fibre $S^1$.  If for some $x \in
S^2$ we pick a 
point $y$ in $h^{-1}(x)$, we may parametrize the elements of $h^{-1}(x)$ as $\lambda
y$ for $\lambda
\in \C$, satisfying $\lambda \bar{\lambda} =1$.

\hspace{2mm}  Let ${\tiny +}= h(0,1) \in S^2$.  We take ${\tiny +}$ as base 
point and identify the
complement of ${\tiny +}$ with the interior of a closed disk $D^2$.  We have a
map $I:D^2 \to S^2$, which restricts to the identification on the interior, and maps
the boundary to ${\tiny +}$.  Then $I:(D_2, \partial D_2) \to (S^2,+)$ 
represents a generator
of $\pi_2(S^2)$.

\hspace{2mm} $I$ lifts to a map sending $D^2$ to $\{(\surd(1-w \bar{w}),w)|w
\bar{w} \leq 1 \}$.  This restricts to a homeomorphism $\partial D^2 \to
h^{-1}({\tiny +})$.

\hspace{2mm}  Hence the map $\pi_2(S^2) \to \pi_1(S^1)$, associated to the
fibre bundle, sends a generator of $\pi_2(S^2)$ to a generator of $\pi_1(S^1)$, and
is an isomorphism.  Consequently,
this map is a $2$- cover, so by lemma 3.1.3, we have $\pi_3(S^2) \cong \pi_3(S^3)
\cong H_3(S^3; \Z)$.  

\lem{$H_3(S^3; \Z) \cong \Z$}

\proof{$S^3$ may be constructed from a 
single point and a single $3$- cell, with associated algebraic complex:
$$\Z \to 0 \to 0 \to \Z$$}

\hspace{2mm}  Let $X$ be a topological space with base point $*$.  We consider 
$(0,1)$ to be the base point of $S^3$.  

\define{}
$h^*:\pi_2(X) 
\to \pi_3(X)$ 
is the map which sends an 
element of $\pi_2(X)$,
represented by a map $\alpha:(S^2,+) \to (X,*)$ to the element of $\pi_3(X)$
represented by $\alpha \circ h$.

\hspace{2mm}We will now define a construction, which takes an arbitrary 
map 

$\alpha:(D^2, \partial D_2) \to (X,*)$ to a map $\alpha^T:(S^3,(0,1)) \to (X,*)$.

 \hspace{2mm} We may identify $D^2$ with the set of complex numbers $\{w|w\bar{w} \leq
\frac{1}{2}\}$.  

\hspace{2mm} Consider the solid torus in $S^3$, given by, 
$\{(z,w) \in S^3| w\bar{w} \leq \frac{1}{2}\}$.  We denote this torus $N$.  Points in 
$N$ may be parametrized by the argument of $z$ and the element of $D^2$, $w$.  Let 
$R_{\theta}$ denote a rotation of $D^2$, about $\theta$ radians anti clockwise.

\hspace{2mm} Given an arbitrary map $\alpha:(D^2, \partial D_2) \to (X,*)$, the map 
$\alpha^T:(S^3,(0,1)) \to (X,*)$ is defined as follows:

$\alpha^T$ sends the
complement of $N$ to the base point.  

If $(\theta, w) $ is an element of $N$
(with respect to the above parametrization), then 

$$
\alpha^T(\theta,w)= \alpha (R_{\theta}w)
$$

\hspace{2mm}Pictorially then, $\alpha^T$ maps the solid torus $N$ to $X$ by rotating the map $\alpha$ as one goes round $N$.

\hspace{2mm}  If $\alpha_t$, $t \in [0,1]$ is a homotopy from $\alpha_0$ to
$\alpha_1$, then $\alpha_t^T$, $t \in [0,1]$ is a homotopy from $\alpha_0^T$ to
$\alpha_1^T$.  Hence  
the map sending $\alpha$ to $\alpha^T$ may regarded as a map of sets $\pi_2(X)
\to \pi_3(X)$.  In fact, we will show that $h^*(\alpha)=\alpha^T$.  (We abuse notation
by denoting elements of homotopy groups by maps representing them.)

\hspace{2mm}Let $1$ denote the generator of $\pi_2(S^2)$ represented by the identity
map 

${\rm Id}:(D^2,\partial D^2) \to (D^2,\partial D^2)$, composed with the natural 
collapse

$c:(D^2,\partial D^2) \to (S^2,+)$.  We may regard $h$ as a map $(S^3,(0,1)) \to (S^2,+)$.

\lem{As elements of $\pi_3(S^2)$, we have $h=1^T$}.

\proof{We may deform $h$ to $h'$ by thickening the preimage of the base
point and pushing the other fibers out accordingly.  Consider the preimage, 
under $h'$, of a point other
than the base point.  Following it round a
circle in the $z$ - plane,
centered on the origin of the $z$- plane, takes one round a circle in the $w$-
plane, centered on the origin of the $w$- plane.  Hence $h'$ is homotopic to $1^T$.}

\hspace{2mm}Note that a map, $\alpha:(D^2, \partial D_2) \to (X,*)$, represents an
element of $\pi_2(X)$, so $h^*(\alpha)$ is a well defined element 
of $\pi_3(X)$.  Specifically, $h^*(\alpha)=\alpha' \circ h$, where 

$\alpha':(S^2,+) \to (X,*)$
satisfies $\alpha' \circ c = \alpha$.

\lem{For any map $\alpha:(D^2, \partial D_2) \to (X,*)$, we have 
$h^*(\alpha)=\alpha^T$ as elements of $\pi_3(S^2)$.} 

\proof{$h^*(\alpha)$ is represented by $\alpha' \circ h$, where 
$\alpha' \circ c = \alpha$.  By the previous lemma
$\alpha' \circ h$ is homotopic to $\alpha' \circ 1^T$. Clearly $\alpha' \circ 1^T$
maps the complement of $N$ to $*$.  Given $(\theta,w) \in N$ we have
$$
\alpha' \circ 1^T(\theta,w)= \alpha'c(R_\theta w)=\alpha(R_\theta w)=\alpha^T
(\theta,w)
$$
${}$\hspace{2mm}So $\alpha' \circ 1^T= \alpha^T$.
}

\hspace{2mm}  The element $0 \in \pi_2(X)$ is represented by the map which sends
$D^2$ to $*$.  Therefore the map $0^T$ sends every point of $S^3$ to the 
$*$, and we have $0^T=0$.  Hence $h^*(0)=0$.

\hspace{2mm} Given $\alpha:(D^2,\partial D_2) \to (X,*)$, we define 
$-\alpha:(D^2,\partial D_2) \to (X,*)$, by setting $-\alpha(w)= \alpha(\bar{w})$ for
all $w \in D^2$.  Clearly as elements of $\pi_2(X)$, we have $\alpha + (-\alpha)=0$,
so there is no ambiguity in the minus sign.

\hspace{2mm}  If we let $z=z_1 +iz_2$, $w=w_1+iw_2$ for $(z_1,z_2,w_1,w_2) \in
\R^4$, then we may regard $S^3$ as sitting in $\R^4$.  The base point $(0,1)$ is 
$(0,0,1,0)$, with respect to this
parametrization. The following matrix represents a rotation of $S^3$, fixing
the base point:

$$
\left[\begin{array}{cccc}  
1&0&0&0 \\
0&{\rm cos}(t)&0&{\rm sin}(t)  \\
0&0&1&0 \\
0&-{\rm sin}(t) &0&{\rm cos}(t) \\
\end{array} \right]
$$

\hspace{2mm}  For $t \in [0, \pi]$, we denote this rotation $L_t$.  We now have
a homotopy, $\alpha^T L_t$, from $\alpha^T$ to $\alpha^T L_\pi$, which keeps the
image of the base point fixed.

\lem{Given $\alpha \in \pi_2(X)$ we have $h^*(-\alpha) = h^*(\alpha)$.}

\proof{We will show that as maps, $(-\alpha)^T =\alpha^T L_\pi$. The 
map $L_\pi$ preserves the solid 
torus, $\{(z,w) \in S^3| w\bar{w} \leq \frac{1}{2}\}$, sending $(\theta,w)$ 
to $(-\theta,\bar{w})$.  Hence 
$$
\alpha^T L_\pi(\theta,w)=\alpha^T(-\theta,\bar{w})=\alpha(R_{-\theta}\bar{w})
= \alpha(\bar{R_\theta w}) = (-\alpha)^T(\theta,w)
$$
\hspace{2mm}So $(-\alpha)^T=\alpha^T L_\pi$, which we know is homotopic to 
$\alpha^T$.  Hence $h^*(-\alpha) = h^*(\alpha)$.
}

\newpage
\lem{Let $g$ be an element of $\pi_1(X)$ and let $\alpha$ be an 
element of $\pi_2(X)$.  Then 
$h^*(\alpha)g =h^*(\alpha g)$.}

\hspace{2mm}Proof:  We must show that $\alpha^T g$ is homotopic 
to  $(\alpha g)^T$.  We regard $\alpha^T$ as a map from the ball $\{p \in
\R^3|\,\,|p| \leq 1\} \to X$, which maps the boundary to $*$.  We interpret $g$
as a map $I \to X$, mapping endpoints to $*$.  We 
regard $\alpha ^Tg$ as a map from the ball $\{p \in
\R^3|\,\,|p| \leq 3\} \to X$ defined as follows:

$\alpha ^Tg(p) = \alpha^T(p)$, if $|P| \leq 1$,

$\alpha ^Tg(p) = g(2-p)$, if $1 \leq |P| \leq 2$,

$\alpha ^Tg(p) = *$, if $|P| \geq 2$.

\hspace{2mm}   Consider the the cylinder
$\{(a,b,c) \in \R^3 | a^2 +b^2  \leq 0.3,\,\,\, -2.5  \leq c \leq 2.5\}$, which
we will denote $C$.  
Without loss of generality, we may assume that $C$ passes through the hole of
the solid torus corresponding to $N$.  

\hspace{2mm} By deforming $\alpha ^Tg$ slightly, to a map $f$, we get its 
restriction to $C$ to be given by 

$$
f(a,b,c)= m(c)
$$

where $m$ is the path $e \cdot g^{-1} \cdot e \cdot g \cdot e$.  Let $m_t$ be a
homotopy, fixing end points, with $m_0(c)=m(c)$ and $m_1(c)=*$.

\hspace{2mm}  We define a homotopy $f_t$ as follows:

$f_t(p)=f(p)$, for $p \notin C$,

$f_t((a,b,c))=m_t(c)$, for $a^2+b^2 \leq 0.2$,

$f_t((a,b,c))=m_{10t(0.3-a^2-b^2 )}(c)$, for $0.3 \leq a^2+b^2 \leq 0.2$,

\hspace{2mm}  Now the solid torus corresponding to $N$ is contained in a larger solid 
torus, $N'$, with the smaller cylinder $C'=\{(a,b,c) \in \R^3 | a^2 +b^2  
\leq 0.2,\,\,\, -2.5  \leq c \leq 2.5\}$, in its hole.  The map $f_1$ sends 
the complement of $N'$ to $*$.  The map $f_1$, restricted to
the cross section of $N'$, corresponding to an argument $\theta$, is homotopic 
to $\alpha g \circ R_\theta$.  

\hspace{2mm} Hence as elements of $\pi_3(X)$, we have $h^*(\alpha)g=\alpha^T g= 
f = f_1= (\alpha g)^T =h^*(\alpha g)$.

\hfill $\Box \,\,$

\newpage
\sec{$\pi_3$ of a $2$- complex}

\hspace{2mm}  In this section we will calculate the Abelian group structure of 
$\pi_3(X)$ where $X$ is a finite geometric complex
whose universal cover is homotopy equivalent to a finite wedge of spheres. This 
was first done by Milnor and Hilton (see \cite{Miln1}).

\bigskip
\hspace{2mm}  Let $X$ be a finite geometric $3$- complex, satisfying 
$C_*(\tilde{X}) \sim \mathcal{A}$, for some algebraic $2$- complex 
$\mathcal{A}$.  From the proof of lemma 2.1.7, for some integer $k$, we have a 
homotopy equivalence $ \tilde{X} \sim W$ where $W= \bigvee_{i=1}^k W_i$, with each
$W_i \cong S^2$.   Let $\{w\}= \bigcap_{i=1}^k W_i$.  As Abelian groups, 

$\pi_n(X) \cong \pi_n(W)$ for $n>1$.

\bigskip
\hspace{2mm}We proceed to construct a $2$- cover of $W$.  The fibre of a $2$- cover 
of $W$,will be $k(\pi_2(W),1)$.  We have $\pi_2(W)=H_2(W; \Z)=\Z^k$, as the 
associated algebraic complex of $W$ is:
$$
\Z^k \to 0 \to \Z$$

\hspace{2mm}  Let $T^k$ denote the product of $k$ copies of $S^1$.

\lem{$k(\pi_2 (X),1)=T^k$}
   
\proof{$\tilde{T^k}= \R^k$, which is contractible. Hence $\pi_n(T^k)=0$ for
$n>1$.  Using the identity $\pi_1(X \times Y)= \pi_1(X) \times \pi_1(Y)$, we
have $\pi_1(T^k)= \Z^k$.}

\bigskip
\hspace{2mm}
Let $S_i$, $i= 1,\cdots,k$ denote $3$- spheres.  For each $i$, let $h_i:S_i \to W_i$ 
denote the Hopf fibration and let $C_i$ denote the fibre over $w$ in each
fibration.

\hspace{2mm} There is a natural inclusion, $\iota:{\lor}_{i=1}^k W_i \hookrightarrow
\prod_{i=1}^k W_i$ which satisfies: 

$\iota(x)_i= x$ if $x \in W_i$ 

$\iota(x)_i= w$ if $x \notin W_i$

for each $i \in 1,\cdots,k$.

\newpage
\hspace{2mm}Hence we may consider the pullback of the fibration $\prod_{i=1}^k h_i$
along $\iota$:

$$
\xymatrix{Y \ar[d]_f \ar[r]&\prod_{i=1}^k S_i \ar[d]^{\prod_{i=1}^k h_i}\\
{\lor}_{i=1}^k W_i \ar[r]^\iota& \prod_{i=1}^k W_i}
$$

\bigskip
\bigskip
\hspace{2mm}
The fibre of $\prod_{i=1}^k h_i$ (and consequently $f$) is $\prod_{i=1}^k C_i 
\cong T^k$

\prop{$f$ is a $2$-cover.}

\proof{The boundary operator $\pi_2(X) \to \pi_1(T^k)$ sends a generator of
$\pi_2(X_i)$ to a generator of $\pi_1(C_i)$.  Hence it may be represented by the
identity matrix and is certainly an isomorphism.}

\hspace{2mm}
From the previous section we know that $\pi_3(X) \cong H_3(Y; \Z)$.  We proceed to
calculate this homology group.

\hspace{2mm}As each $S_i$ is a 3- sphere, it may be constructed
as a CW- complex from the following cells; 

the point $w$, 

a $1$- cell, $e_i$ for $C_i$, 

a pair of $2$- cells $U^+_i, U^-_i$ 

a pair of $3$- cells $V^+_i, V^-_i$

\bigskip
\hspace{2mm}  For each $i$, the boundary maps on these cells are: 

\bigskip

$\partial e_i =0$ 

$\partial U_i^+ =e_i$ 

$\partial U_i^- =-e_i$ 

$\partial V_i^+ =U_i^+ +U_i^-$ 

$\partial V_i^- =-U_i^+-U_i^-$

\bigskip
\hspace{2mm}  The product space $\prod_{i=1}^k S_i$ naturally inherits the structure
of a CW- complex with cells given by tensor products of the cells constituting the
$S_i$.  The boundary map of a tensor product is calculated using the identity
$$\partial(C \otimes D)=\partial C \otimes D + (-1)^{{\rm deg}C} \otimes 
\partial D$$  (See  \cite{Vick}, p120).




\bigskip
\hspace{2mm}$Y$ is a subcomplex of $\prod_{i=1}^k S_i$.  It is the preimage
under $\prod_{i=1}^k h_i$ of $W \subset \prod_{i=1}^k S_i$.  Hence $(x_1, \cdots x_k)
\in \prod_{i=1}^k S_i$ is an element of $Y$ precisely when $\#\{i|\,\,x_i \notin C_i\}
\leq 1$.  We enumerate the 3- cells 
which constitute $Y$ and use the identity above to calculate the boundary of each.

\hspace{2mm}The generators of $C_3(Y)$ are:

$e_i \otimes e_j \otimes e_l$, for $i <j<l$,

$V_i^+$ and $V_i^-$,

$e_i \otimes U_j^+$ for $i<j$ and $U_j^+ \otimes e_i$ for $j<i$,

$e_i \otimes U_j^-$ for $i<j$ and $U_j^- \otimes e_i$ for $j<i$,

\bigskip
\hspace{2mm}
Applying the boundary operator to each of these gives:

$\partial (e_i \otimes e_j \otimes e_l)=0$,

$\partial V_i^+=U_i^++U_i^-$,  $\quad \partial V_i^-=-U_i^+-U_i^-$,

$\partial (e_i \otimes U_j^+)=-e_i \otimes e_j$,  
$\quad \partial U_j^+ \otimes e_i=e_i \otimes e_j $,

$\partial (e_i \otimes U_j^-)=e_i \otimes e_j$,  
$\quad \partial U_j^- \otimes e_i=-e_i \otimes e_j $.

\bigskip
\hspace{2mm}
Hence we see that the kernel of the boundary operator is generated by

\bigskip
$e_i \otimes e_j \otimes e_l$, for $i <j<l$,

$V_i^+ + V_i^-$,

$e_i \otimes U_j^+  + U_i^+ \otimes e_j$ for $i<j$,

$e_i \otimes U_j^+ + e_i \otimes U_j^-$ for $i<j$,

$U_i^+ \otimes e_j  + U_i^- \otimes e_j$ for $i<j$.

\bigskip
\hspace{2mm}  We now enumerate the 4- cells of $Y$.  The generators of $C_4(Y)$ are:

\bigskip
$e_i \otimes e_j \otimes e_l \otimes e_m$ for $i<j<l<m$,

$e_i \otimes e_j \otimes U^+_l$ for $i<j<l$,

$e_i \otimes U^+_j \otimes e_l$ for $i<j<l$,

$U^+_i \otimes e_j \otimes e_l$ for $i<j<l$,

$e_i \otimes e_j \otimes U^-_l$ for $i<j<l$,

$e_i \otimes U^-_j \otimes e_l$ for $i<j<l$,

$U^-_i \otimes e_j \otimes e_l$ for $i<j<l$,

$e_i \otimes V_j^+$ for $i<j$,

$V_i^+ \otimes e_j$ for $i<j$,

$e_i \otimes V_j^-$ for $i<j$,

$V_i^- \otimes e_j$ for $i<j$.

\bigskip
\hspace{2mm}  Applying the boundary operator gives 

$\partial (e_i \otimes e_j \otimes e_l \otimes e_m)=0$

$\partial (e_i \otimes e_j \otimes U^+_l)=e_i \otimes e_j \otimes e_k$ 

$\partial (e_i \otimes U^+_j \otimes e_l)=-e_i \otimes e_j \otimes e_k$ 

$\partial (U^+_i \otimes e_j \otimes e_l)=e_i \otimes e_j \otimes e_k$ 

$\partial (e_i \otimes e_j \otimes U^-_l)=-e_i \otimes e_j \otimes e_k$ 

$\partial (e_i \otimes U^-_j \otimes e_l)=e_i \otimes e_j \otimes e_k$ 

$\partial (U^-_i \otimes e_j \otimes e_l)=-e_i \otimes e_j \otimes e_k$ 

$\partial (e_i \otimes V_j^+)=-e_i\otimes U_j^+-e_i\otimes U_j^-$ 

$\partial (V_i^+ \otimes e_j)=U_i^+ \otimes e_j + U_i^- \otimes e_j $ 

$\partial (e_i \otimes V_j^-)=e_i\otimes U_j^++e_i\otimes U_j^-$ 

$\partial (V_i^- \otimes e_j)=-U_i^+ \otimes e_j - U_i^- \otimes e_j$ 

\bigskip
\hspace{2mm}Hence the image of $C_4(Y)$ under the boundary
operator is generated by:

\bigskip
$e_i \otimes e_j \otimes e_l$, for $i <j<l$,

$e_i \otimes U_j^+ + e_i \otimes U_j^-$ for $i<j$,

$U_i^+ \otimes e_j  + U_i^- \otimes e_j$ for $i<j$.

\bigskip
\hspace{2mm}  This leaves: 

$V_i^+ + V_i^-$,

$e_i \otimes U_j^+  + U_i^+ \otimes e_j$ for $i<j$,

\bigskip
to generate $H_3(Y; \Z)$.

\bigskip
\hspace{2mm}
Let $J$ denote the algebraic $\pi_2$ of
$\mathcal{A}$.  As a $\ZG$ - module, $J$ is isomorphic to $\pi_2(X)$.  This is
isomorphic to  $\pi_2(W)$ as an Abelian group, which in turn is equal to 

$$
\bigoplus_{i=1}^k \pi_2(W_i)
$$

\hspace{2mm} We regard $w$ as the base point of each $W_i$.  For each $i$, let $w_i$ 
denote the generator of $\pi_2(W_i)$, induced from the identity map 
$W_i \to S^2$.  The isomorphism 

$\pi_2(W) \to
\pi_2(X)$ is induced by a homotopy equivalence $W \to \tilde{X}$ composed with
the covering map $\tilde{X} \to X$.  Let $q:W \to X$ denote this composition.

\hspace{2mm} For each $i$, let $\alpha_i:S^2 \to X$ denote $q \circ w_i$ and 
 let $*$ denote $qw$.  Taking $*$ once and for all as
the 
base point of $X$, we may
regard $\alpha_i$ as representing an element of $\pi_2(X)$.  In particular,
$\alpha_i$ represents $q_*w_i \in \pi_2(X)$.

\hspace{2mm}
The $\alpha_i$, $i\,\, \in\, 1, \cdots ,k$ are a basis for $J$ over $\Z$.  Therefore 
the $\alpha_i \otimes \alpha_j$,
$i,j\,\, \in\, 1, \cdots ,k$ are a basis for $J \otimes_{\Z}J$ over $\Z$ .  Hence, 
any element of 
$S^2(J)$ may be written uniquely as a $\Z$- linear combination 
of $\alpha_i \otimes \alpha_j$.  
However, in any such linear combination, the coefficient on 
$\alpha_i \otimes \alpha_j$ would have to equal the
coefficient on $\alpha_j \otimes \alpha_i$.  Consequently, $S^2(J)$
is generated, over $\Z$, by the $\alpha_i \otimes \alpha_i$ and the 
$\alpha_i \otimes\alpha_j+\alpha_j \otimes \alpha_i$, $i \neq j$. 

\newpage
\hspace{2mm}  We may define a $\Z$- linear 
isomorphism $\phi:\pi_3(X) \to S^2(J)$,
by sending:

$V_i^+ + V_i^- \mapsto \alpha_i \otimes \alpha_i$,

$e_i \otimes U_j^+  + U_i^+ \otimes e_j \mapsto \alpha_i \otimes 
\alpha_j + \alpha_j \otimes
\alpha_i$. 

\bigskip
\hspace{2mm}  The purpose of the next section is to show that $\phi$ is $\ZG$- 
linear.

\sec{Realizing elements of $\pi_3 (X)$.}

\hspace{2mm}We have the following maps of sets:

\xymatrix{J\ar[r]^{h^*} \ar[d]_q& \pi_3{X} \ar[dl]_\phi\\
S^2(J)}

\bigskip
where $q: \alpha \mapsto \alpha \otimes \alpha$, for all $\alpha \in J$.  

\hspace{2mm}From lemma 3.2.7, we know that $h^*$ respects the action 
of $G$. Also $q$ respects the action of $G$ by construction.  In order to prove 
that $\phi$ respects the action of $G$ (and consequentially is a $\ZG$- linear isomorphism), we
will show the following:

\bigskip
i) $h^*(J)$ generates $\pi_3(X)$ over $\Z$.  We obtain this in corollary 3.4.8.  

ii)The diagram above commutes.  In other 
words, $\phi(h^*(\alpha))= \alpha \otimes \alpha$ for all $\alpha \in J$.  We obtain this
in lemma 3.4.13.

\bigskip

\hspace{2mm}  In the previous section $\pi_3 (X)$ was computed, as an 
Abelian group.  We now realize the elements of this 
Abelian group as actual maps $S^3 \to X$.

\hspace{2mm}  Consider the identity map, $S^3 \to S_i$, taking $w$ as base 
point.  This map sends the generator of $H_3(S^3; \Z)$ to the
element of $H_3(Y; \Z)$ represented by $V_i^+ + V_i^-$.  Hence the Hurewicz
isomorphism takes $V_i^+ + V_i^-$ to the element of $\pi_3(Y)$ represented by the
identity map $S^3 \to S_i$.

\hspace{2mm}  The $2$- cover, $f$, restricts to the Hopf 
fibration $S_i \to W_i$.  So  $V_i^+ + V_i^-$ is sent to the element of
$\pi_3(W)$, represented by the Hopf fibration $S^3 \to W_i$.  This equals $h^*(w_i)$. 

\hspace{2mm}  
Composing $h^*(w_i)$ with $q$
gives $q \circ  h^*(w_i) = q \circ w_i \circ h=\alpha_i \circ h =h^*(\alpha_i)$.  
Hence the isomorphism $H_3(Y; \Z) \to \pi_3(X)$, due to lemma 3.1.3, takes  
$V_i^+ + V_i^-$ to 

$h^*(\alpha_i) \in \pi_3(X)$.  

\lem{For each $i\in1,\cdots,k$, we have $\phi(\alpha_i^T)=\alpha_i \otimes \alpha_i$}.

\proof{From the construction of $\phi$ and the discussion above, we have
\newline
$\phi(h^*(\alpha_i))=\alpha_i \otimes \alpha_i$.$\quad$  So by lemma 
3.2.5 $\phi(\alpha_i^T)=\phi(h^*(\alpha_i))=\alpha_i \otimes \alpha_i$}

\hspace{2mm}Once again, we may view $S^3$ as $\{(z,s) \in \C ^2|z \bar{z}+s
\bar{s}=1\}$.  It decomposes into two solid tori, given by $s\bar{s} \leq 
\frac{1}{2}$ and $s\bar{s} \geq \frac{1}{2}$.  Given $i<j$, we may identify the
disk $U_j^+$ with the disk $s\bar{s} \leq \frac{1}{2}$ in the complex plane.  We
may also identify the set of arguments of $z$, on the complex plane, with $C_i$.

\hspace{2mm} Similarly we may identify the
disk $U_i^+$ with the disk $z\bar{z} \leq \frac{1}{2}$ in the complex plane.  We
may also identify the set of arguments of $s$, on the complex plane, with $C_j$.

\hspace{2mm}  Hence we have a natural identification map $I:S^3 \to 
C_i \times U_j^+  \cup U_i^+ \times C_j$.  This takes a generator of
$H_3(S^3;\Z)$ to $e_i \otimes U_j^+  + U_i^+ \otimes e_j$.  Hence we have

\lem{The Hurewicz
isomorphism, $H_3(Y;\Z) \to \pi_3(Y)$, takes $e_i \otimes U_j^+  + U_i^+ 
\otimes e_j$ to the element of $\pi_3(Y)$ represented by $I$.}  

\hspace{2mm}  For any $a \in C_i$ and $b \in C_j$, the map $f$ restricts to 
$w_j:(a,U_j^+) \to W_j$ and to $w_i:(U_i^+,b) \to W_i$, where $w_i$, $w_j$ are now
regarded as maps $D^2 \to W$, which send the boundary to $w$.

\hspace{2mm}Fix two simply linked solid tori, $A, B$, each parametrized 
$S^1 \times D^2 $, in $S^3$.  Let $Q$ be a topological space with base 
point $\%$.  Given $\alpha, \beta\in\pi_2(Q)$, we construct an 
element of $\pi_3(Q)$:

\define{} $\alpha \vee \beta:S^3 \to Q$ is defined by

$\alpha \vee \beta(p)=\%$ for $p \notin A,B$

$\alpha \vee \beta(\theta,d)=\alpha(d)$ for $(\theta, d) \in A$

$\alpha \vee \beta(\theta, d)=\beta(d)$ for $(\theta, d) \in B$

\hspace{2mm}  By thinning the solid tori, $s\bar{s} \leq \frac{1}{2}$ and 
$s\bar{s} \geq \frac{1}{2}$, we have that $f_*I$ and $ w_i \vee w_j$ represent 
the same element of $\pi_3(W)$.  

\lem{$q \circ w_i \vee w_j = \alpha_i \vee \alpha_j$.}

\proof{
\newline
$
q \circ w_i \vee w_j (\theta,d)= q \circ w_i (d)= \alpha_i(d)$ 
for $(\theta,d) \in A$.
\newline
$
q \circ w_i \vee w_j (\theta,d)= q \circ w_j (d)= \alpha_j(d)$ 
for $(\theta,d) \in B$.
\newline
${}$\hspace{2mm}Hence we have 
$q \circ w_i \vee w_j = \alpha_i \vee \alpha_j$.}

\hspace{2mm}So to recap: the Hurewicz isomorphism takes $e_i \otimes U_j^+  + U_i^+ 
\otimes e_j$ to $I \in \pi_3(Y)$.  Composition with $f$ takes this to $w_i \vee w_j$. 
Composition with $q$ takes this to $\alpha_i \vee \alpha_j$.  

\hspace{2mm} We may conclude:

\lem{As an Abelian group, $\pi_3(X)$ is generated freely by 
$\alpha_i^T$,$i\in\{1,\cdots, k \}$ and $\alpha_i \vee \alpha_j$, $i<j \in 
\{1,\cdots,k\}$.  Also,
\newline
${}\qquad\phi(\alpha_i^T)= \alpha_i \otimes
\alpha_i$, $\quad i\in\{1,\cdots, k \}$
\newline
${}\qquad\phi(\alpha_i \vee \alpha_j)= \alpha_i \otimes \alpha_j+
\alpha_j \otimes \alpha_i$, $i<j \in \{1,\cdots,k\}$.}

\bigskip
\hspace{2mm}  We now describe two homotopies which will be used in the
lemmas which follow.
\newline

i)\hspace{2mm}  Let $\gamma \in \pi_3(X)$ agree with $\alpha^T$ on the solid torus
$\{(z,s) \in S^3| s\bar{s} \leq \frac{1}{2}\}$, for some $\alpha \in \pi_2(X)$. 
So given a point $(\theta, s)$ in the solid torus, we have
$\gamma(\theta,s)=\alpha(R_\theta(s))$.

\hspace{2mm}  For $t \in [0,\pi]$, let

$\gamma_t(\theta, s)= \alpha(s), \quad$ for $\theta \leq t$.

$\gamma_t(\theta, s)= \alpha (R_{\small \frac{2\pi(\theta -t)}{2\pi-t}}(
s)), \quad$ for $\theta \geq t$.

$\gamma_t(p)= \gamma(p), \quad$ for $p$ not in the solid torus.

\hspace{2mm}  Intuitively, $\gamma_\pi$ is $\gamma$ with the "twist" moved round
to one side of the torus.
\newline

ii)\hspace{2mm}  This homotopy will be referred to as "pinching" the solid torus.

\hspace{2mm}  Consider $\gamma_\pi$ as before, and select a cylinder,
in $S^3$
which intersects the solid torus in two disconnected places, on the "untwisted"
side, and which does not intersect any other points which do not map to $*$,
under $\gamma_\pi$.  Parametrize this cylinder $A \times I$, where $A$ is a
disk of radius $3$ about the origin in $\C$, and $I$ is the interval $[-1,1]$.  
Let $A_1$ denote the disk of radius $\frac{1}{\surd2}$ about 
$i$, in $A$ and $A_2$ denote the 
disk of radius $\frac{1}{\surd2}$ about $-i$, in $A$.

\hspace{2mm}By deforming $\gamma_\pi$ to $\gamma'$, fixing the base point, we
may have $\gamma'$ restricted to the cylinder giving:

$\gamma'(a,s)= *, \quad$ if $a \notin A_1,A_2$,

$\gamma'(a,s)= \alpha(a-i), \quad$ if $a \in A_1$,

$\gamma'(a,s)= \alpha(-a-i), \quad$ if $a \in A_2$.

\hspace{2mm}For each $s \in I$, we have the map $\gamma'(\_,s):A \to X$ 
representing
$\alpha-\alpha=0 \in \pi_2$. Hence, fixing some $s$, we have a 
homotopy 

$h_t:A \to X$, for $t \in
[0,1]$, with $h_0(a)=\gamma'(a,s) $, and $h_1(a)=*$.

\hspace{2mm}  Let $\gamma'_t$, $\quad t\in[0,1]$ be defined as follows:

$\gamma'_t(p) =\gamma'(p), \quad$ if $p \notin A \times I$,

$\gamma'_t(a,s)=h_t(a), \quad$ for $|s| \leq \frac{1}{2}$,

$\gamma'_t(a,s)=h_{(2-2|s|)t}(a), \quad$ for $|s| \geq \frac{1}{2}$,

\hspace{2mm}  So $\gamma'_1$ represents the same element of $\pi_2$ as
$\gamma_\pi$, but it has the "twisted" part of the solid torus "pinched off"
from the untwisted part.

\lem{Let $\alpha,\beta \in \pi_2(X)$.  Then $\alpha \vee \beta = (\alpha +
\beta)^T - \alpha^T - \beta^T$.}

\proof{Consider the map $(\alpha +
\beta)^T:S^3 \to X$.  We may partition the solid torus 
$\{(z,w) \in S^3| w\bar{w} \leq \frac{1}{2}\}$, into two linked solid tori, such
that altering $(\alpha +
\beta)^T$ to map one solid tori to $*$ would leave $\alpha^T$, and altering 
$(\alpha +
\beta)^T$ to map the other solid tori to $*$, would leave $\beta^T$.
\newline
${}$\hspace{2mm} Hence "pinching off" the "twists" in the two solid tori gives
$$
(\alpha +\beta)^T =  \alpha^T + \beta^T +\alpha \vee \beta  
$$}
\hspace{2mm}  Note the left hand side of the equality in lemma 3.4.6 is
symmetric in $\alpha$ and $\beta$.  

\hspace{2mm}Hence we have:

\cor {$\alpha \vee \beta= \beta \vee \alpha$}

\hspace{2mm}Also, as each $\alpha \vee \beta$ may now be written in terms of
$(\alpha + \beta)^T$, $\alpha^T$ and $ \beta^T$, we have: 

\cor{Elements of the form $\alpha^T$, for $\alpha \in \pi_2(X)$, span $\pi_3(X)$
over $\Z$.}

\lem{Let $\alpha,\beta, \gamma \in \pi_2(X)$.  Then
$$
(\alpha+\beta+\gamma)^T=(\alpha + \beta)^T +(\beta +\gamma  )^T +
(\gamma + \alpha)^T - \alpha^T - \beta^T - \gamma^T
$$}

\proof{Consider the map $(\alpha +
\beta + \gamma)^T:S^3 \to X$.  We may partition the solid torus 
$\{(z,w) \in S^3| w\bar{w} \leq \frac{1}{2}\}$ into three mutually linked
solid tori, $A,B,C$ so that if the complement of $A$, $B$ or $C$ was
mapped to $*$, the map would be $\alpha^T$, $\beta^T$, $\gamma^T$ respectively.
\newline
${}$\hspace{2mm} Pinching off the twist in each solid torus gives
$$
(\alpha +\beta + \gamma)^T = \alpha^T+\beta^T+\gamma^T + \Delta
$$
where $\Delta$ is represented by a map $S^3 \to X$, sending the complement of
three linked tori, $A'$, $B'$, $C'$
 to $*$, and projecting each solid tori ($=S^1 \times D^2$)
onto $D^2$ and mapping via $\alpha$, $\beta$, $\gamma$ respectively.
\newline
${}$\hspace{2mm}We may pinch each solid torus between its links with the other
two, to get
$$
\Delta = \alpha \vee \beta + \beta \vee \gamma + \gamma \vee \alpha
$$
${}$\hspace{2mm} Hence
$$
(\alpha +\beta + \gamma)^T = \alpha^T+\beta^T+\gamma^T + 
\alpha \vee \beta + \beta \vee \gamma + \gamma \vee \alpha
$$
$$
=(\alpha + \beta)^T +(\beta +\gamma  )^T +
(\gamma + \alpha)^T - \alpha^T - \beta^T - \gamma^T
$$}

\hspace{2mm}  We now apply this lemma to get:

\lem{Suppose $\alpha, \beta, \gamma \in \pi_2(X)$ satisfy the following:
\newline
$\phi(\alpha^T)= \alpha\otimes \alpha$, $\quad$
$\phi(\beta^T)= \beta\otimes \beta$, $\quad$
$\phi(\gamma^T)= \gamma\otimes \gamma$, $\quad$
\newline
$\phi((\alpha+ \beta)^T) = (\alpha+ \beta) \otimes (\alpha+ \beta)$,
\newline
$\phi((\alpha+ \gamma)^T) = (\alpha+ \gamma) \otimes (\alpha+ \gamma)$,
\newline
$\phi((\gamma+ \beta)^T) = (\gamma+ \beta) \otimes (\gamma+ \beta)$,
\newline
Then $\phi((\alpha+\beta + \gamma)^T )= (\alpha+\beta + \gamma) \otimes
(\alpha+\beta + \gamma)$.
}

\hspace{2mm} Proof:  
$$\phi((\alpha+\beta + \gamma)^T) \quad =  \quad \quad
\phi((\alpha + \beta)^T +(\beta +\gamma  )^T +
(\gamma + \alpha)^T - \alpha^T - \beta^T - \gamma^T)
$$\hfill (by 3.4.9)
{\small
$$
=(\alpha + \beta)\otimes (\alpha + \beta) +(\beta +\gamma  ) \otimes 
(\beta +\gamma)  +
(\gamma + \alpha)\otimes (\gamma + \alpha) - \alpha \otimes \alpha - 
\beta \otimes \beta - \gamma \otimes \gamma
$$\hfill} (by hypothesis)
$$
=\alpha \otimes \alpha + \beta \otimes \beta  + \gamma \otimes \gamma + 
\alpha \otimes \beta + \beta \otimes \alpha +
\beta \otimes \gamma  + \gamma \otimes \beta +
\alpha \otimes \gamma + \gamma \otimes \alpha 
$$
$$
=(\alpha+\beta + \gamma) \otimes (\alpha+\beta + \gamma)
$$

\hfill $\Box \,\,$

\hspace{2mm}  As we have said, the $\alpha_i$ form a basis for $J$ over $\Z$.  

\define{Norm}  $\quad$ Given 
$ \alpha = \sum_r \alpha_i \lambda_i$, $\lambda_i \in \Z$,
we define the norm of $\alpha$, denoted $|\alpha|$, by
$$
|\alpha| = \sum_r |\lambda_i|
$$

\lem{Let $\alpha \in J$ satisfy $|\alpha| \leq 2$.  Then 
$\phi(\alpha^T)= \alpha \otimes \alpha$}

\hspace{2mm}Proof:  The only element of $J$ with norm equal to $0$ is $0$.  We
have previously noted that $0^T=0$, so we have $\phi(0^T)=\phi(0)=0= 0 \otimes
0$.

\hspace{2mm}The only elements of $J$ with norm equal to $1$ are ones of the 
form $\alpha_i$ or $-\alpha_i$.  By lemma 3.4.5, $\phi(\alpha_i^T)=\alpha_i
\otimes \alpha_i$.  By lemma 3.2.6, we have
 
$\phi(-\alpha_i^T)=\phi(\alpha_i^T)=
\alpha_i
\otimes \alpha_i=  
-\alpha_i
\otimes -\alpha_i$.

\hspace{2mm}Elements of $J$ with norm equal to $2$ are of the form $\alpha_i +
\alpha_j  $, $\quad-\alpha_i -
\alpha_j $, $\quad\alpha_i -
\alpha_j $, $\quad\alpha_i2 $,  \hspace{2mm} or $-\alpha_i2$, where 
$i \neq j$.  By lemma 3.2.6, it is
sufficient to consider $\alpha_i +\alpha_j  $, $\quad \alpha_i -\alpha_j$,
\hspace{2mm} and $\alpha_i2$.

\hspace{2mm}By lemma 3.4.6, for $i \neq j$, we 
have $$\phi((\alpha_i +\alpha_j)^T ) = \quad 
\phi((\alpha_i +\alpha_j)^T - \alpha_i^T - \alpha_j^T + \alpha_i^T +
\alpha_j^T)$$ 
$$=\phi(\alpha_i \vee \alpha_j) +\phi(\alpha_i^T) +\phi( \alpha_j^T)
= \quad \alpha_i \otimes \alpha_j +\alpha_j \otimes \alpha_i +\alpha_i \otimes 
\alpha_i +\alpha_j \otimes \alpha_j
$$
$$
=(\alpha_i+\alpha_j) \otimes (\alpha_i+\alpha_j)
$$

\hspace{2mm}  By lemma 3.4.9, we have 
$$
\alpha_i^T=(\alpha_i+\alpha_j-\alpha_j)^T= (\alpha_i+\alpha_j)^T +
(\alpha_i-\alpha_j)^T - \alpha_i^T-\alpha_j^T2
$$
\hspace{2mm} Hence for $i \neq j$, we have
$$
\phi((\alpha_i-\alpha_j)^T)= \phi(\alpha_i^T)2
+\phi(\alpha_j^T)2 -\phi((\alpha_i+\alpha_j)^T)
$$
$$
=\alpha_i \otimes \alpha_i +\alpha_j \otimes \alpha_j  -
\alpha_i \otimes \alpha_j
-\alpha_j \otimes \alpha_i
$$
$$
=(\alpha_i-\alpha_j) \otimes (\alpha_i-\alpha_j)
$$

\hspace{2mm}  Again, by lemma 3.4.9, we have 
$$
\alpha_i^T=(\alpha_i+\alpha_i-\alpha_i)^T= (\alpha_i+\alpha_i)^T 
 -\alpha_i^T3
$$

\hspace{2mm}  So $\phi((\alpha_i2)^T)=\phi(\alpha_i^T)4= (\alpha_i \otimes
\alpha_i) 4 = \alpha_i2 \otimes \alpha_i2$.

\hspace{2mm} Hence we have dealt with each element of $J$, with norm
less than or equal to $2$.
 $ \Box$

\lem{Let $\alpha \in J$.  Then $\phi(\alpha^T)= \alpha \otimes \alpha$.}

\hspace{2mm}Proof:  Suppose not.  Then there exists some minimal number $n$ such
that there exists $\alpha \in J$, with $|\alpha| =n $ 
and $\phi(\alpha^T) \neq \alpha \otimes \alpha$.  By lemma 3.4.12, $n\geq 3$, so we
have $\alpha = \alpha' \pm \alpha_i \pm \alpha_j$, for some $i,j \in
\{1,\cdots,k\}$ and $|\alpha'|=n-2$. 

\hspace{2mm}  The hypothesis' for lemma 3.4.10 are fulfilled, as $1,2,n-1,n-2 <
n$, so we have 
$$
\phi(\alpha^T)= \phi((\alpha' \pm \alpha_i \pm \alpha_j)^T)
$$
$$
=(\alpha' \pm \alpha_i \pm \alpha_j) \otimes (\alpha' \pm \alpha_i \pm \alpha_j)
$$
$$
= \alpha \otimes \alpha
$$
which contradicts $\phi(\alpha^T) \neq \alpha \otimes \alpha$.

\hfill $\Box \,\,$

\hspace{2mm}Hence by lemma 3.2.5, we have $\phi(h^*(\alpha))= \alpha \otimes \alpha$.

\lem{The isomorphism $\phi:\pi_3(X) \to S^2(J)$ is $\ZG$- linear.}

\hspace{2mm}Proof:  Recall lemma 3.2.7, which states that given $g \in G$, and
$\alpha \in J$, we have $h^*(\alpha g)=h^*(\alpha)g$.  From lemma 3.4.13 we have 

$$
\phi(h^*(\alpha) g)= \phi(h^*(\alpha g))=\alpha g \otimes \alpha g = 
(\alpha \otimes \alpha)g= \phi(h^*(\alpha))g
$$ 

\hspace{2mm}  By corollary 3.4.8, given an arbitrary element $\beta \in
\pi_3(X)$, we may write 

$$
\beta= \sum_r h^*(\beta_r) \lambda_r
$$

with $\beta_r \in J$, $\lambda_r \in \Z$.

\hspace{2mm} So
$$
\phi(\beta g)=\sum_r\phi(h^*(\beta_r) g)\lambda_r = 
\sum_r\phi(h^*(\beta_r)) g\lambda_r = \phi(\beta)g
$$

\hfill $\Box \,\,$

\hspace{2mm}  From this lemma we may conclude:

{\bf Theorem 3A} {\it If $C_*(\tilde{X}) \sim \mathcal{A}$, for some finite algebraic
$2$- complex $\mathcal{A}$, then 
$\pi_3(X) \cong S^2(J)$ as $\ZG$- modules.}

\hspace{2mm} In particular:

\cor{If $X$ is a geometric $2$- complex, and $\pi_2(X)=J$, then $\pi_3(X) \cong
S^2(J)$ as $\ZG$- modules.}

\sec{The effect of stabilizing $\pi_2$.}

\hspace{2mm}  Let $X$ be a finite geometric $2$- complex with finite fundamental
group $G$.  Let $X'$ be another finite geometric $2$- complex with 
fundamental group $G$.  From corollary 1.1.2 and lemma 2.1.5 we know that there
exist integers $a,b$ such that 

$$
\pi_2(X) \oplus \ZG^a \cong \pi_2(X') \oplus \ZG^b
$$

\hspace{2mm}  In this section we investigate the corresponding relationship
between $\pi_3(X)$ and $\pi_3(X')$.

\hspace{2mm}  If we let $J= \pi_2(X)$ and $J'= \pi_2(X')$, then from the last
section we have 

$\pi_3(X)= S^2(J)$ and $\pi_3(X')= S^2(J')$.  We know that 

$$
S^2(J\oplus \ZG^a) \cong S^2(J' \oplus \ZG^b)
$$

\hspace{2mm}  The next few lemmas give us an expansion of this.

\lem{Let $A_i$, $i=1, \cdots, n$ be modules over $\ZG$, with finitely generated, 
free underlying Abelian groups.  Then
$$
S^2(\bigoplus_{i=1}^n A_i)= \bigoplus_{i=1}^n S^2(A_i) \oplus \bigoplus_{i <j} A_i
\otimes_\Z A_j
$$
}

\proof{For each $i$, let the $e_{i,r}$ be a basis over $\Z$ for $A_i$.  Then 
$S^2(\bigoplus_{i=1}^n A_i)$ is freely generated over $\Z$, by elements of the
form: 
\newline
$e_{i,r} \otimes e_{i,r}$, 
\newline
$e_{i,r} \otimes e_{i,s}+e_{i,s} \otimes e_{i,r},\quad$ $r<s$,
\newline
$e_{i,r} \otimes e_{j,s}+e_{j,s} \otimes e_{i,r},\quad$ $i<j$,
\newline
${}$\hspace{2mm}  For each $i$, the $\Z$- linear span of the 
$e_{i,r} \otimes e_{i,r}$, for all $r$, and the 
$e_{i,r} \otimes e_{i,s}+e_{i,s} \otimes e_{i,r},\quad$ $r<s$, is closed under
the group action and is isomorphic over $\ZG$ to $S^2(A_i)$.
\newline
${}$\hspace{2mm}Similarly, for each pair $i,j$, with $i <j$, the 
$\Z$- linear span of the 
\newline
$e_{i,r} \otimes e_{j,s}+e_{j,s} \otimes e_{i,r}$, for
all $r$ and $s$, is closed under the group action. Furthermore we have an 
isomorphism from it to $A_i \otimes_\Z A_j$, which maps 
$$e_{i,r} \otimes e_{j,s}+e_{j,s} \otimes e_{i,r} \mapsto e_{i,r} \otimes
e_{j,s}$$\hspace{2mm}So over $\ZG$, the module $S^2(\bigoplus_{i=1}^n A_i)$ 
decomposes into
the sum
$$
\bigoplus_{i=1}^n S^2(A_i) \oplus \bigoplus_{i <j} A_i
\otimes_\Z A_j
$$}

\hspace{2mm} Note that for $\ZG$- modules $A$, $B$, we have a $\ZG$- linear 
isomorphism 

$A \otimes_\Z B \to B \otimes_\Z A$, given by sending $a
\otimes b \mapsto b \otimes a$.

\lem{Let $A$ be a $\ZG$ module whose underlying Abelian group is isomorphic to
$\Z^k$.  Then $A \otimes_\Z \ZG \cong \ZG^k$.}

\hspace{2mm}Proof:  The action of an element $g \in G$ on $A$ is a $\Z$- linear
isomorphism.  Hence, if $e_1, \cdots, e_k$ is a $\Z$- linear basis for $A$, then
so is $e_1g, \cdots, e_kg$.  Hence $A \otimes_\Z \ZG$ is freely generated 
over $\Z$ by elements of the form $e_ig \otimes g$, for $g\in G$ and $i \in
\{1,\cdots,k\}$.

\hspace{2mm}  For each $i$, the $\Z$ linear span 
of the $e_ig \otimes g$, $g\in G$, is closed under the action of $G$. 
Furthermore, we have a $\ZG$- linear isomorphism from it to $\ZG$, which maps 
$e_ig \otimes g \mapsto g$.

\hspace{2mm}Hence,
$$
A \otimes_\Z \ZG \cong \bigoplus_{i=1}^k\ZG \cong \ZG^k
$$

\hfill $\Box \,\,$

\hspace{2mm}Let $p$ denote the number of pairs, $\{g, g^{-1}\} \in G$, with $g
\neq g^{-1}$.  

\bigskip
\hspace{2mm}Also, let $T = \{g \in G | g^2=e\}$.

\define{} For a finite group $G$, we set 

$$V_G = \bigoplus_{t \in T} (1+t) \ZG
$$

\hspace{2mm}Let $n$ denote the order of $G$.

\lem{$S^2(\ZG) \cong \ZG^{1+p} \oplus V_G$}

\hspace{2mm}Proof:  Consider the $\Z$- linear basis for $S^2(\ZG)$ given by $g
\otimes g$, $g\in G$, and 

$g \otimes h + h\otimes g$, $\{g,h\} 
\subset G$.  

\hspace{2mm}  We have $<e \otimes e> \ZG=<g \otimes g|\,\,g\in G>\Z$. 

\hspace{2mm}  Suppose $g \neq g^{-1}$.  Then $<g\otimes
e + e \otimes g>\ZG =$

$ <h \otimes l +l
\otimes h|\,\,h,l \in G,\, hl^{-1}=g>\Z$.  Further, given an element 
$h \otimes l + l\otimes h$, satisfying $hl^{-1}=g$, we may  write $h \otimes l + 
l\otimes h =(g\otimes e + e \otimes g)l$.  This is the unique way of writing 
$h \otimes l + l\otimes h$ as a $\ZG$- linear multiple of $
g\otimes e + e \otimes g$, because $h \otimes l + 
l\otimes h \neq (g\otimes e + e \otimes g)h$, as $gh=l$ would imply
$g=lh^{-1}=(hl^{-1})^{-1}=g^{-1}$.

\hspace{2mm}So the $\Z$- linear span of the $g \otimes g$ and the $h \otimes l 
+ l\otimes h$, $lh^{-1} \neq hl^{-1}$, is isomorphic to $\ZG^{1+p}$.  The
remaining elements of the basis are of the form $h \otimes l + l\otimes h$, with
$hl^{-1}=lh^{-1}$.  Let $t=hl^{-1}$.  then $h \otimes l 
+ l\otimes h = (e \otimes t + t \otimes e)l$.  The $\Z$-linear span of the 
$h \otimes l + l\otimes h$, $hl^{-1}=lh^{-1}$, is therefore equal to 
$$
\bigoplus_{t \in T}(e \otimes t + t \otimes e)\ZG
$$
\hspace{2mm}Regarding $S^2(\ZG)$ as a submodule of $\ZG \otimes_\Z \ZG$, we have
$e \otimes t + t \otimes e=(e\otimes t)(1+t)$.  Therefore 
$$
\bigoplus_{t \in T}(e \otimes t + t \otimes e)\ZG = V_G
$$
and
$$
S^2(\ZG) \cong \ZG^{1+p} \oplus V_G$$
\hfill $\Box \,\,$

\hspace{2mm}Returning to our geometric complexes $X$ and $X'$, we had 

$$
S^2(J\oplus \ZG^a) \cong S^2(J' \oplus \ZG^b)
$$

\hspace{2mm}By lemma 3.5.1, this expands to:

$$
S^2(J) \oplus S^2(\ZG)^a \oplus (J \otimes_\Z \ZG)^a \oplus 
(\ZG \otimes_\Z \ZG)^{a(a-1)/2}
$$
$$
=S^2(J') \oplus S^2(\ZG)^b \oplus (J' \otimes_\Z \ZG)^b \oplus 
(\ZG \otimes_\Z \ZG)^{b(b-1)/2}
$$
\hspace{2mm}  Let $k$ denote the $\Z$- rank of $J$, and let $k'$ denote the
$\Z$- rank of $J'$.  By theorem 3A we have $S^2(J)= \pi_3(X)$ 
and $S^2(J')= \pi_3(X')$.  Applying lemmas 3.5.2 and 3.5.4 to the equation above
gives:
$$
\pi_3(X) \oplus (\ZG^{1+p} \oplus V_G)^a \oplus \ZG^{ka} \oplus 
ZG^{na(a-1)/2}
$$
$$
=\pi_3(X') \oplus (\ZG^{1+p} \oplus V_G)^b \oplus \ZG^{k'b} \oplus 
ZG^{nb(b-1)/2}
$$

\hspace{2mm} Hence 
$$
\pi_3(X) \oplus \ZG^{(1+p+k+n(a-1)/2)a} \oplus {V_G}^a
=\pi_3(X') \oplus \ZG^{(1+p+k'+n(b-1)/2)b} \oplus {V_G}^b
$$

and we have the following theorem:

\thm{Let $X$ and $X'$ be finite geometric $2$- complexes, with finite 
fundamental group $G$. Then there exist integers $p,q,r,s$ such that:
$$
\pi_3(X) \oplus \ZG^p \oplus {V_G}^q
=\pi_3(X') \oplus \ZG^r \oplus {V_G}^s
$$
}

\hspace{2mm}Note that if the order of $G$ is odd, it does not contain any
elements of order $2$.  Hence $V_G=0$ and we have.  

\cor{Let $X$ and $X'$ be finite geometric $2$- complexes, with odd finite 
fundamental group $G$. Then $\pi_3(X)$ and $\pi_3(X')$ are stably equivalent.
}






\newpage
\sec{The case $\pi_2=IG^*$} 

\hspace{2mm}  In this section we consider the case of a finite geometric $2$-
complex, $X$, with finite fundamental group $G$ and $\pi_2(X)=IG^*$.  Recall
that we have a short exact sequence

$$
0 \to \Z \to \ZG \stackrel{p}{\to} IG^*\to0
$$

\hspace{2mm}  We now consider 
the $\ZG$- linear surjection $p':S^2(\ZG) \to S^2(IG^*)$, which
sends 

$g \otimes g \mapsto p(g) \otimes p(g)$ 

and $g \otimes h +h \otimes g
\mapsto p(g) \otimes p(h) + p(h) \otimes p(g)$.

\hspace{2mm}  Let $\Sigma$ denote the sum of the elements of $G$.

\lem{The kernel of $p'$ is generated over $\ZG$, by $\Sigma \otimes e + e
\otimes \Sigma$ and $\Sigma \otimes \Sigma$.}

\proof{We may take a $\Z$- linear basis for $\ZG$, given by $\{g \in G |g \neq
e\}$ together with $\Sigma$.  We have a $\Z$- linear basis for $IG^*$ given by
\newline
$\{p(g) |g \in G, g \neq e \}$.  Hence the kernel of $p'$ is generated over $\Z$
by $\Sigma \otimes \Sigma$, and $g \otimes \Sigma + \Sigma \otimes g$.
Clearly this is contained in the, $\ZG$ linear span of  $\Sigma \otimes e + e
\otimes \Sigma$ and $\Sigma \otimes \Sigma$.}

\hspace{2mm} Let $S$ be a subset of $G$, containing precisely one of $g$ or
$g^{-1}$, for each pair $g, g^{-1}$, with $g \neq g^{-1}$.  Let $e_g = g \otimes
e + e \otimes g$ for each $g \in S$.  Let $e_t=t \otimes
e + e \otimes t$ for each $t \in T$.  Let $e_0= e \otimes e$.  Then from the
proof of lemma 3.5.4 we have:
$$
S^2(\ZG)= e_0\ZG \oplus \bigoplus_{g \in S} e_g\ZG \oplus 
\bigoplus_{t \in T} e_t\ZG
$$
\hspace{2mm}Also, 
$$
e\otimes \Sigma + \Sigma \otimes e = 
e_0 2+\sum_{g \in S}e_g(1+g^{-1}) + \sum_{t\in T}e_t
$$
and
$$
\Sigma \otimes \Sigma  = 
(e_0 2+\sum_{g \in S}e_g(1+g^{-1}) + \sum_{t\in T}e_t) \frac{\Sigma}{2}
$$

Let 

$$
u=(e_0 2+\sum_{g \in S}e_g(1+g^{-1}) + \sum_{t\in T}e_t)
$$ 

We have 

$$
S^2(IG^*)= \frac {S^2(\ZG) } { \quad <u, u {\tiny \frac{\Sigma}{2}}>}
$$

\hspace{2mm}The relation $u=0$ is equivalent to 
$e_0 2=-(\sum_{g \in S}e_g(1+g^{-1}) + \sum_{t\in T}e_t)$.  Let 
$$M=\bigoplus_{g \in S} e_g\ZG \oplus 
\bigoplus_{t \in T} e_t\ZG
$$

and let $u_M \in M$ be equal to 
$-(\sum_{g \in S}e_g(1+g^{-1}) + \sum_{t\in T}e_t)$.




\thm{If $\pi_2{X}=IG^*$ then $\pi_3(X)=M[\frac{u_M}{2}]$.}

\proof{We need to check that $\frac{u_M}{2}$ satisfies the relations 
which $e_0$ was subject to:
\newline
$\frac{u_M}{2}2+\sum_{g \in S}e_g(1+g^{-1}) + \sum_{t\in T}e_t=0$
\newline
$(\frac{u_M}{2}2+\sum_{g \in S}e_g(1+g^{-1}) + \sum_{t\in T}e_t)\Sigma/2=0$
}

\eg{$G=C_3=<x|x^3=e>$}

\hspace{2mm}  If $\pi_2(X)=IG^*$ then 
$$\pi_3(X)=S^2(IG^*)=e_{x^2}\ZG[(1+x)/2] \cong \ZG[\frac{(1+x)}{2}]
$$
\eg{$G=Q_8=<x,y|x^2=y^2,xyx=y>$}

\hspace{2mm}There is only one element of order two in $Q_8$ so if
$\pi_2(X)=IG^*$ we have
$$
\pi_3(X)=S^2(IG^*) = (e_x \ZG \oplus e_y \ZG \oplus e_{xy}  \ZG \oplus e_{y^2}
\ZG)[\frac{u_M}{2}]
$$
where $u_m=-(e_x(1+x^3)+e_y(1+y^3)+e_{xy}(1+xy^3)+e_{x^2})$.  

\hspace{2mm}As $e_{y^2}\ZG \cong (1+y^2)\ZG$, (see lemma 3.5.4), we have a
single relation $e_{y^2}((1-y^2)=0$.

\bigskip
\hspace{2mm}  We remark that rationally, for some $a \in \Z$, 
$\pi_2(X) \otimes \Q \cong IG^* \otimes
\Q \oplus \QG^a$, for any finite geometric $2$- complex $X$, 
with finite fundamental group, $G$.  Consequentially, we have 
$$
S^2(J) \otimes \Q = S^2(J \otimes \Q)= S^2(IG^* \oplus \ZG^a) \otimes \Q
$$

\hspace{2mm} From lemma 3.6.1, the kernel of $p'$ is generated by $u$ and
$u\Sigma/2$.  Hence the kernel of $p' \otimes \Q:S^2(\QG) \to S^2(IG^*) \otimes
\Q$, is generated by $u$, as $u\Sigma/2$ is in the $\QG$- linear span of $u$.

\hspace{2mm}We therefore have a short exact sequence:
$$
0 \to \QG \to S^2(\QG) \to S^2(IG^*) \otimes \Q \to 0
$$

\hspace{2mm}As surjections over $\QG$ split and cancellation of finitely
generated modules holds over $\QG$, we may write
$$
 S^2(IG^*) \otimes \Q =  S^2(\QG) / \QG
$$

\hspace{2mm}  So from lemmas 3.5.1, 3.5.2 and 3.5.4 we may conclude 

\thm{There exist integers
$a$, $b$, such that
$$
\pi_3(X) \otimes \Q = \QG^a \oplus (V_G \otimes \Q)^b
$$
}

\hspace{2mm}  Again, if the order of $G$ is odd, then $V_G=0$, so 

\cor{If the order of $G$ is odd, then $\pi_3(X)$ is rationally free.}

\hspace{2mm} Note that example 3.6.3 is a case in point.

\sec{Summary}

\hspace{2mm}In this chapter we have shown that given a geometric $2$- complex, $X$,
with finite fundamental group $G$, we
have $\pi_3(X) \cong S^2(J)$, where $J=\pi_2(X)$ (theorem 3A).  We have defined a 
module over
$\ZG$, $V_G$ and shown that $\pi_3(X)$ is determined by $G$, up to stabilization by 
copies of $\ZG$ and copies of $V_G$ (theorem 3.5.5).  Rationally, we have shown that 
$\pi_3(X) \otimes \Q \cong \QG^a \oplus (V_G \otimes \Q)^b$ for integers $a,b$
(theorem 3.6.5).

\hspace{2mm}In the case where $G$ is a group of odd order, we have 
$V_G \cong 0$.  Hence in this case, the stable class of $\pi_3(X)$ is determined and
$\pi_3(X)$ is rationally free (corollaries 3.5.6 and 3.6.6).

\chapter{Algebraic Poincare 5- complexes}\label{}
\stepcounter{ch}
\setcounter{sec}{0}

\hspace{2mm}Let $M$ be a closed, connected, orientable $5$- dimensional manifold,
with finite fundamental group $G$ (we assume manifolds to be without boundary).  In this chapter we consider algebraic
complexes $C_*(\tilde{M'})$, where $M'$ is a finite CW- complex, with $M \sim M'$.  $C_*(\tilde{M'})$ must satisfy 
Poincare duality.  We use this to show 
that up to chain homotopy equivalence, we may represent it by an algebraic $2$- complex, $\mathcal{A}$, connected to its dual via a 
$G$- invariant bilinear form, $\beta$,
on $(\pi_2(\mathcal{A}))^*$.  We denote the resulting algebraic $5$- complex $(\mathcal{A}, \beta)$.

\hspace{2mm}In \S4.3 we show that the algebraic $2$- complex $\mathcal{A}$, is not important in the sense that any algebraic $2$- complex
may be stabilized to one which, together with the appropriate bilinear form, represents the homotopy type of $C_*(\tilde{M'})$.In \S4.4, we 
describe chain homotopy
equivalences between these algebraic complexes.

\hspace{2mm}  We next consider the homotopy equivalence induced by Poincare Duality.  In particular we are interested in how similar it can
be made to the identity.  In \S4.5 we show that it can be taken as the identity on $4$ of the $6$ terms of the chain complex.  In \S4.6
however, we find a homological obstruction to this homotopy equivalence actually being the identity.  In particular, certain
manifolds described in \cite{Bard} do not satisfy the homological condition necessary, for being able to write the homotopy equivalence as
the identity.

\sec{The category ${\rm TOP}^5$}   

\hspace{2mm}   Fix a finite group $G$.   Let ${\rm TOP}^5$ denote the category of closed connected orientable five dimensional 
topological manifolds with base point, with respect to which the fundamental group is identified with $G$. The morphisms in this category are continuous 
maps which preserve the base point and induce the identity on $G$.

\hspace{2mm} Given an object of ${\rm TOP}^5$, we may find a finite $CW$- complex, which is homotopy equivalent to it 
(see \cite{Kirb}).  Let the following be the algebraic chain complex of the universal cover of the  
CW- complex:

$$
C_5 \stackrel{\partial_5}{\longrightarrow} C_4 \stackrel{\partial_4}{\longrightarrow}C_3 \stackrel{\partial_3}{\longrightarrow}
C_2 \stackrel{\partial_2}{\longrightarrow} C_1 \stackrel{\partial_1}{\longrightarrow}C_0 
$$

\hspace{2mm}  This is an algebraic complex of free $\ZG$- linear modules and $\ZG$- linear maps.  It is exact at $C_1$ and the
cokernel of $\partial_1$ is $\Z$.

\hspace{2mm}  We say that an algebraic complex satisfies Poincare Duality if it is chain homotopy equivalent to to its dual.

\prop{(Poincare Duality) The algebraic complex, $(C_*,\partial_*)$ satisfies Poincare Duality.}

\hspace{2mm} Note that Poincare duality tells us that the algebraic complex above is chain homotopy equivalent to to its dual.  Therefore 
it is exact at $C_4$ and the kernel of $\partial_5$ is $\Z$.  As the Euler characteristic
of this complex is minus that of its dual, it must be $0$.

\hspace{2mm} Let ${\rm ALG}^5$ denote the category of algebraic 5-complexes of finitely generated free $\ZG$ modules, 

$$
F_5 \stackrel{\partial_5}{\longrightarrow} F_4 \stackrel{\partial_4}{\longrightarrow}F_3 \stackrel{\partial_3}{\longrightarrow}
F_2 \stackrel{\partial_2}{\longrightarrow} F_1 \stackrel{\partial_1}{\longrightarrow}F_0 
$$
satisfying:

i) Poincare duality.

ii) Exactness at $F_4$ and $F_1$

iii) The cokernel of $\partial_1$ and the kernel of $\partial_5$, equaling $\Z$. 

\bigskip
\hspace{2mm}The morphisms in this category are homotopy equivalence classes of chain maps.

\hspace{2mm} We may define a functor $C: {\rm TOP}^5 \to {\rm ALG}^5$ by choosing a homotopy equivalence, $h_{M}:M \to M'$ with $M'$ a
finite CW- complex, for each $M \in{\rm TOP}^5$.  $C(M)$ is then defined to be $C_*(\tilde{M'})$.

\hspace{2mm}Given a continuous map, which is a morphism in ${\rm TOP}^5$, 
$f:M_1 \to M_2$, we may select a cellular map $f'$, which is homotopic to $(h_{M_2} \circ f \circ h_{M_1}^{-1}):M_1' \to M_2'$.  Define $C(f)$ to 
be the equivalence class of the 
chain map $f'_*:C_*(\tilde{M_1')} \to C_*(\tilde{M_2')}$.

\hspace{2mm}The isomorphism class of $C(M)$ in ${\rm ALG}^5$ is an invariant of $M$ as by construction, different choices of $M'$ must be
homotopy equivalent to $M$ and hence each other.

\sec{Dual Form}

\hspace{2mm} Given an algebraic two complex, over $\ZG$, $J^*\to F_2 \to F_1 \to F_0 \to \Z \,$ and a $G$- invariant bilinear form $\beta$ on $J$, we can
 associate an algebraic 5- complex:

\begin{eqnarray*}
F_0^* \to F_1^* \to F_2^* \longrightarrow F_2 \to F_1 \to F_0 \quad\, \\
\searrow \,\,\,\,\quad \nearrow \qquad \qquad\qquad\,\,\\
J \to J^* \qquad \qquad \qquad\,\,\,\,\,
\end{eqnarray*}

\hspace{2mm} Let ${\rm DUAL}^2$ denote the category whose objects consist of:

(i)  An algebraic 2- complex of finitely generated free modules over $\ZG$, 

$J^*\dashrightarrow F_2 \to F_1 \to F_0 \dashrightarrow \Z \,$, with exactness at $F_1$.

\bigskip
(ii) A $G$- invariant bilinear form $\beta$, on $J$, such that the associated algebraic 5- complex is an
element of ${\rm ALG}^5$.

\hspace{2mm}  As before we define the morphisms of ${\rm DUAL}^2$, to be homotopy equivalence classes of chain maps
between the associated algebraic 5- complexes of objects in ${\rm DUAL}^2$.  If an element of ${\rm ALG}^5$ is the associated
algebraic 5- complex of an element of ${\rm DUAL}^2$, we say it is in "dual form".

\hspace{2mm}  We have a functor $i:{\rm DUAL}^2 \to {\rm ALG}^5$ which sends an object to its associated
algebraic 5- complex, and sends a morphism to the class of chain map which it represents.  This functor is clearly full and faithful. 
We will show that every object in ${\rm ALG}^5$ is isomorphic, in the category ${\rm ALG}^5$, to an object in
 the image of $i$.

\thm{Every element of ${\rm ALG}^5$ is chain homotopy equivalent to an algebraic 5- complex in dual form.}

\bigskip
\hspace{2mm} Proof:  We start with an arbitrary element of ${\rm ALG}^5$:

$$
C_5 \stackrel{\partial_5}{\longrightarrow} C_4 \stackrel{\partial_4}{\longrightarrow}C_3 \stackrel{\partial_3}{\longrightarrow}
C_2 \stackrel{\partial_2}{\longrightarrow} C_1 \stackrel{\partial_1}{\longrightarrow}C_0 
\eqno(1)
$$

\hspace{2mm} We perform three simple homotopy equivalences.  Firstly, the complex (1) is chain homotopy equivalent to to 

$$
C_5 \oplus C_0^*\stackrel{\delta_5}{\longrightarrow} C_4 \oplus C_0^*\stackrel{\partial_4 \oplus 0}{\longrightarrow}C_3 
\stackrel{\partial_3}{\longrightarrow}
C_2 \stackrel{\partial_2}{\longrightarrow} C_1 \oplus C_5^* \stackrel{\delta_1}{\longrightarrow}C_0 \oplus C_5^* 
\eqno(2)
$$

where

$$
\delta_1 = \left( \begin{array}{cc} \partial_1&0\\ 0&1 \end{array} \right) \qquad\qquad
\delta_5 = \left( \begin{array}{cc} \partial_5&0\\ 0&1 \end{array} \right)
$$ 

\bigskip
\bigskip
\hspace{2mm} Let $R_0=C_0$, $R_5=C_5$, and $R_1 = C_1 \oplus C_5^*$, $R_4 = C_4 \oplus C_0^*$.  Then (2) can be written

$$
R_5 \oplus R_0^*\stackrel{\delta_5}{\longrightarrow} R_4 \stackrel{\partial_4 \oplus 0}{\longrightarrow}C_3 
\stackrel{\partial_3}{\longrightarrow}
C_2 \stackrel{\partial_2}{\longrightarrow} R_1 \stackrel{\delta_1}{\longrightarrow}R_0 \oplus R_5^* 
\eqno(3)
$$ 

\bigskip
\bigskip
\hspace{2mm} Again we perform a pair of simple homotopy equivalences.  The complex (3) is chain homotopy equivalent to to 

$$
R_5 \oplus R_0^*\stackrel{\delta_5}{\longrightarrow} R_4 \oplus R_1^*\stackrel{\delta_4}{\longrightarrow}C_3 \oplus R_1^* 
\stackrel{\partial_3 \oplus 0}{\longrightarrow}
C_2 \oplus R_4^* \stackrel{\delta_2}{\longrightarrow} R_1 \oplus R_4^*\stackrel{\delta_1 \oplus 0}{\longrightarrow}R_0 \oplus R_5^* 
\eqno(4)
$$ 

where 

$$
\delta_2 = \left( \begin{array}{cc} \partial_2&0\\ 0&1 \end{array} \right) \qquad\qquad
\delta_4 = \left( \begin{array}{cc} \partial_4 \oplus 0 &0\\ 0&1 \end{array} \right)
$$ 

\bigskip
\bigskip
\hspace{2mm} Let $R_2 = C_2 \oplus R_4^*$, $R_3 = C_3 \oplus R_1^*$.  Then (4) can be written

$$
R_5 \oplus R_0^*\stackrel{\delta_5}{\longrightarrow} R_4 \oplus R_1^*\stackrel{\delta_4}{\longrightarrow}R_3  
\stackrel{\partial_3 \oplus 0}{\longrightarrow}
R_2 \stackrel{\delta_2}{\longrightarrow} R_1 \oplus R_4^*\stackrel{\delta_1 \oplus 0}{\longrightarrow}R_0 \oplus R_5^* 
\eqno(5)
$$ 

\bigskip
\hspace{2mm}  As all the modules in this complex are free and the Euler characteristic is 0, we can assume the existence of some
isomorphism $\theta:R_2^* \to R_3^*$

\bigskip
We perform a final homotopy equivalence to get 

$$
R_5 \oplus R_0^*\stackrel{\delta_5}{\longrightarrow} R_4 \oplus R_1^*\stackrel{\delta_4}{\longrightarrow}R_3 \oplus R_2^* 
\stackrel{\delta_3}{\longrightarrow}
R_2 \oplus R_3^*\stackrel{\delta_2 \oplus 0}{\longrightarrow} R_1 \oplus R_4^*\stackrel{\delta_1 \oplus 0}{\longrightarrow}R_0 
\oplus R_5^* 
\eqno(6)
$$ 

where 

$$
\delta_3 = \left( \begin{array}{cc} \partial_3 \oplus 0&0\\ 0&\theta \end{array} \right) \qquad\qquad
$$ 

\bigskip
The algebraic complex (1) is therefore chain homotopy equivalent to to (6).  We will show that (6) is chain isomorphic to an algebraic 
5-complex in dual form.

\bigskip
\lem{There exist maps $h_0$, $k_0$, such that the following diagrams commute:}

\begin{eqnarray*}
R_0 \oplus R_5^* \stackrel{\epsilon \oplus 0}{\longrightarrow} M\,\,   
\qquad \qquad       R_0 \oplus R_5^* \stackrel{\epsilon \oplus 0}{\longrightarrow} M\,\,\\
\downarrow h_0 \,\,\quad\,\, \quad \downarrow 1         \quad \qquad \qquad       \uparrow k_0 \,\, \quad\,\, \quad 
\uparrow 1\\ 
R_5^* \oplus R_0 \stackrel{\epsilon'\oplus 0}{\longrightarrow} M\,\,  \qquad \qquad       R_5^* \oplus R_0 \stackrel{\epsilon' 
\oplus
0}{\longrightarrow} M\,\,\\   
\end{eqnarray*}

\bigskip
\hspace{2mm}  Proof:  As the $R_i$ are projective, we may pick $f_0$, $g_0$ so that the following diagrams commute:

\setcounter{equation}{8}
\begin{eqnarray*}
R_0 \stackrel{\epsilon}{\longrightarrow} \Z\,\,   \qquad \qquad       R_0 \stackrel{\epsilon}{\longrightarrow} \Z\,\, \,\,
\nonumber\\ \downarrow f_0 \,\, \quad \downarrow 1             \quad\,\,\,\, \qquad       \uparrow g_0  \,\, \quad \uparrow 1 
\nonumber\\ 
R_5^* \stackrel{\epsilon'}{\longrightarrow} \Z\,\,  \qquad \qquad       R_5^* \stackrel{\epsilon'}{\longrightarrow} \Z\,\,\,\,  \\
\end{eqnarray*}
$\hfill (9)$

Define $h_0:R_0 \oplus R_5^* \to R_5^* \oplus R_0$ and $k_0:R_5^* \oplus R_0 \to R_0 \oplus R_5^*$ by

$$
h_0= \left(\begin{array}{cc} f_0&1-f_0g_0\\ 1& -g_0 \end{array}\right) \qquad \qquad 
k_0= \left(\begin{array}{cc} g_0&1-g_0f_0\\ 1& -f_0 \end{array}\right)
$$

Direct calculation shows that $h_0k_0 = 1$ and $k_0h_0 = 1$.

\bigskip
Also from commutativity of (9), we deduce

$$
(\epsilon' \quad 0)\left(\begin{array}{cc} f_0&1-f_0g_0\\ 1& -g_0 \end{array}\right) = (\epsilon'f_0 \quad
\epsilon'(1-f_0g_0))=(\epsilon \quad 0) 
$$

and

$$
(\epsilon \quad 0)\left(\begin{array}{cc} g_0&1-g_0f_0\\ 1& -f_0 \end{array}\right) = (\epsilon g_0 \quad
\epsilon(1-g_0f_0))=(\epsilon' \quad 0) 
$$

\bigskip
Hence the following diagrams commute:

\begin{eqnarray*}
R_0 \oplus R_5^* \stackrel{\epsilon \oplus 0}{\longrightarrow} M\,\,   
\qquad \qquad       R_0 \oplus R_5^* \stackrel{\epsilon \oplus 0}{\longrightarrow} M\,\,\\
\downarrow h_0 \,\,\quad\,\, \quad \downarrow 1         \quad \qquad \qquad       \uparrow k_0 \,\, \quad\,\, \quad 
\uparrow 1\\ 
R_5^* \oplus R_0 \stackrel{\epsilon'\oplus 0}{\longrightarrow} M\,\,  \qquad \qquad       R_5^* \oplus R_0 \stackrel{\epsilon' 
\oplus
0}{\longrightarrow} M\,\,\\   
\end{eqnarray*}

\hfill $\Box$ \,\,

\bigskip
\lem{There exist a pair of inverse chain isomorphisms between the exact
sequences:}

$$
R_2 \oplus R_3^*\stackrel{\delta_2 \oplus 0}{\longrightarrow} R_1 \oplus R_4^*\stackrel{\delta_1 \oplus 0}{\longrightarrow}R_0 
\oplus R_5^* \stackrel{\epsilon \oplus 0}{\dashrightarrow} \Z \to 0
\eqno(7)
$$

{\it and}

$$
R_3^* \oplus R_2\stackrel{\delta_4^* \oplus 0}{\longrightarrow} R_4^* \oplus R_1\stackrel{\delta_5^* \oplus 0}{\longrightarrow}R_5^* 
\oplus R_0 \stackrel{\epsilon' \oplus 0}{\dashrightarrow} \Z \to 0
\eqno(8)
$$

\hspace{2mm} Proof: We will construct a pair of inverse chain isomorphisms, $h$, $k$, between $(7)$
and  $(8)$.

\bigskip
\hspace{2mm}  We have already defined $h_0$ and $k_0$.  Now suppose that for $i= 0 \,\,{\rm or}\,\, 1$, we have defined 
$h_j:R_j \oplus R_{5-j}^* \to R_{5-j}^* \oplus R_j$ and 
$k_j:R_{5-j}^*\oplus R_j 
\to R_j \oplus R_{5-j}^*$ for $j=0,\dots,i-1$, so that for each $j$, we have $h_jk_j=1$ and $k_jh_j=1$.  We proceed by induction.

\bigskip
As before, pick $f_i$, $g_i$ so that the following diagrams commute:

\begin{eqnarray*}
R_i \quad \stackrel{\delta_i}{\longrightarrow} \quad R_{i-1} \oplus R_{n-i+1}^*\,\,   
\qquad \qquad    R_i \quad \stackrel{\delta_i}{\longrightarrow} \quad R_{i-1} \oplus R_{n-i+1}^*\,\, \nonumber\\  
\downarrow f_i \quad \,\quad\quad \quad \, \quad \downarrow  h_{i-1} \quad \,\,\,\qquad \qquad \uparrow g_i \qquad \,\, \quad \quad \quad 
\uparrow k_{i-1} \,\,\,\,\quad \nonumber\\ 
R_{n-i}^*\stackrel{\delta_{n-i+1}^*}{\longrightarrow} R_{n-i+1}^* \oplus R_{i-1} \,\,  
\qquad \qquad   R_{n-i}^* \stackrel{\delta_{n-i+1}^*}{\longrightarrow} R_{n-i+1}^* \oplus R_{i-1} \,\, \quad 
\end{eqnarray*}
$\hfill (10)$

\bigskip
Define $h_i:R_i \oplus R_{n-i}^* \to R_{n-i}^* \oplus R_i$ and $k_i:R_{n-i}^* \oplus R_i \to R_i \oplus R_{n-i}^*$ by

$$
h_i= \left(\begin{array}{cc} f_i&1-f_ig_i\\ 1& -g_i \end{array}\right) \qquad \qquad 
k_i= \left(\begin{array}{cc} g_i&1-g_if_i\\ 1& -f_i \end{array}\right)
$$

Direct calculation shows that $h_ik_i = 1$ and $k_ih_i = 1$.

\bigskip
Recall $h_{i-1}k_{i-1}=1$ and $k_{i-1}h_{i-1}=1$.  From commutativity of (10) we deduce

$$
(\delta_{n-i+1}^* \quad 0)\left(\begin{array}{cc} f_i&1-f_ig_i\\ 1& -g_i \end{array}\right) = (\delta_{n-i+1}^* f_i \quad
\delta_{n-i+1}^*(1-f_ig_i))=h_{i-1}(\delta_i \quad 0) 
$$

and

$$
(\delta_i \quad 0)\left(\begin{array}{cc} g_i&1-g_if_i\\ 1& -f_i \end{array}\right) = (\delta_i g_i \quad
\delta_i(1-g_if_i))=k_{i-1}(\delta_{n-i+1}^* \quad 0) 
$$

Hence the following diagrams commute:

\begin{eqnarray*}
R_i \oplus R_{n-i}^* \quad \stackrel{\delta_i \oplus 0}{\longrightarrow}  R_{i-1} \oplus R_{n-i+1}^*\,\,   
\qquad \qquad    R_i \oplus R_{n-i}^* \quad \stackrel{\delta_i \oplus 0}{\longrightarrow} \quad R_{i-1} \oplus R_{n-i+1}^*\,\,\\  
\downarrow h_i \,\,\,\quad \quad \, \quad \quad \qquad \downarrow h_{i-1} \,\, \qquad\,\,\,\qquad \qquad \uparrow k_i 
\qquad\qquad \,\,  \quad \quad 
\quad 
\uparrow k_{i-1}\,\,\,\,\quad \\ 
R_{n-i}^*\oplus R_i \stackrel{\delta_{n-i+1}^* \oplus 0}{\longrightarrow} R_{n-i+1}^* \oplus R_{i-1} \,\,  
\qquad \qquad   R_{n-i}^* \oplus R_i\stackrel{\delta_{n-i+1}^* \oplus 0}{\longrightarrow} R_{n-i+1}^* \oplus R_{i-1} \,\,  \quad
\end{eqnarray*}

\bigskip
\hspace{2mm} Together with the identity on $\Z$, the $h_i$, $k_i$ are therefore a pair of mutually inverse chain isomorphisms, 
between $(7)$ and $(8)$.

\begin{eqnarray*}
R_2 \oplus R_3^*\stackrel{\delta_2 \oplus 0}{\longrightarrow} R_1 \oplus R_4^*\stackrel{\delta_1 \oplus 0}{\longrightarrow}R_0 
\oplus R_5^* \stackrel{\epsilon \oplus 0}{\dashrightarrow} \Z \to 0\\
\downarrow h_2 \,\,\, \quad\quad\quad \downarrow h_1 \quad \quad \quad \quad \downarrow  h_0 \,\quad\quad \downarrow 1 \,\, \quad \\
R_3^* \oplus R_2\stackrel{\delta_4^* \oplus 0}{\longrightarrow} R_4^* \oplus R_1\stackrel{\delta_5^* \oplus 0}{\longrightarrow}R_5^* 
\oplus R_0 \stackrel{\epsilon' \oplus 0}{\dashrightarrow} \Z \to 0
\end{eqnarray*}

\hfill $\Box$ \,\,

\bigskip
Let $S_0 = R_0 \oplus R_5^*$,  $S_1 = R_1 \oplus R_4^*$,  $S_2 = R_2 \oplus R_3^*$.  

\bigskip
Also let $d_1= \delta_1 \oplus 0$ and $d_2= \delta_2 \oplus 0$.  Let $d_3 = \delta_3 k_2^*$.

\bigskip
\hspace{2mm} We complete the proof of the theorem with the following lemma:

\bigskip
\lem{ The complex (6) is chain isomorphic to 
\newline
$$
S_0^* \stackrel{d_1^*}{\longrightarrow} S_1^* \stackrel{d_2^*}{\longrightarrow} S_2^* \stackrel{d_3}{\longrightarrow} 
S_2 \stackrel{d_2}{\longrightarrow} S_1 \stackrel{d_1}{\longrightarrow} S_0 
\eqno(11)
$$
}

\hspace{2mm} Proof:  We have the following chain isomorphism:

\begin{eqnarray*}
R_5 \oplus R_0^*\stackrel{\delta_5}{\longrightarrow} R_4 \oplus R_1^*\stackrel{\delta_4}{\longrightarrow}R_3 \oplus R_2^* 
\stackrel{\delta_3}{\longrightarrow}
R_2 \oplus R_3^*\stackrel{\delta_2 \oplus 0}{\longrightarrow} R_1 \oplus R_4^*\stackrel{\delta_1 \oplus 0}{\longrightarrow}R_0 
\oplus R_5^* \quad\\
\downarrow h_0^* \qquad \qquad \downarrow h_1^* \,\,\, \quad \qquad \downarrow h_2^* 
\qquad \qquad\downarrow 1 \qquad \qquad\downarrow 1 \qquad \qquad \downarrow 1 \qquad 
\\
\quad S_0^* \quad\stackrel{d_1^*}{\longrightarrow} \quad S_1^* \quad \stackrel{d_2^*}{\longrightarrow} \quad \quad S_2^*  \quad
\stackrel{d_3}{\longrightarrow} \quad\quad
S_2 \quad \stackrel{d_2}{\longrightarrow} \quad S_1 \quad \stackrel{d_1}{\longrightarrow} \quad S_0 \quad \quad\quad
\end{eqnarray*}

We need only verify that the central square commutes: $d_3 h_2^* = \delta_3 k_2^* h_2^* = \delta_3$.

\hfill $\Box$ \,\,

\bigskip
This completes the proof of theorem 4.2.1
\hfill $\Box$ \,\,

\hspace{2mm}We may conclude:

\thm{The functor $i:{\rm DUAL}^2 \to {\rm ALG}^5$ is full, faithful and surjective up to
isomorphism.  Hence $i$ is a natural equivalence.}

\hspace{2mm} This means that when using the functor $C$ to provide an invariant of an element of ${\rm TOP}^5$, up 
to isomorphism in ${\rm ALG}^5$, we may 
work in the category ${\rm DUAL}^2$.  In order to parametrize the values this invariant can take, we need only classify the forms $\beta$ on 
elements of $\Omega_{-3}^{\ZG}(\Z)$, which give rise to elements of ${\rm ALG}^5$.

\sec{Polarization}   

\bigskip
\hspace{2mm} We would prefer to work over a fixed algebraic 2- complex, and and have the form completely determine the resulting
5- complex.  To that end fix any algebraic 2- complex, of finitely generated free modules: 

$$
\mathcal{A} \quad =\quad  J^* \stackrel{}{\to} F_2 \stackrel{d_2}{\to} F_1 \stackrel{d_1}{\to} F_0 \to \Z \,
$$

with exactness at $F_1$.  Let 

$$
\mathcal{A}^n \quad =\quad  F_2 \oplus \ZG^n \stackrel{d_2 \oplus 0}{\to} F_1  \stackrel{d_1}{\to} F_0 
$$

\bigskip
\hspace{2mm}  Let $J_n$ denote $J \oplus \ZG^n$.  The following result opens the possibility of classifying algebraic
5- complexes, without in any way having to classify algebraic 2-complexes. 

\bigskip
\thm{Any element of ${\rm ALG}^5$ is chain homotopic to the associated algebraic 5- complex of the element of 
${\rm DUAL}^2$, represented by $(\mathcal{A}^n,\,\, \gamma)$, for some $n$ and bilinear form, $\gamma$ on $J_n$.}

\bigskip
\hspace{2mm} Proof:  Any element of ${\rm ALG}^5$ is chain homotopy equivalent to to an algebraic complex of the form 

$$
T_0^* \stackrel{\Delta_1^*}{\longrightarrow}  T_1^*  \stackrel{\Delta_2^*}{\longrightarrow} T_2^* \stackrel{\Delta_3}
{\longrightarrow} T_2  \stackrel{\Delta_2}{\longrightarrow} T_1  \stackrel{\Delta_1}{\longrightarrow} T_0 
\eqno(12)
$$

for finitely generated free modules $S_i$.

\bigskip
\lem{For some integer, $n$, and free module $T$, we may apply a pair of simple homotopy equivalences to $\mathcal{A}^n$ 
to get 
\newline
$$L_2 \stackrel{D_2}{\longrightarrow} L_1  \stackrel{D_1}{\longrightarrow} L_0$$
\newline
and a pair of simple homotopy equivalences to  
\newline
$$T_2 \oplus T \stackrel{\Delta_2\oplus 0}{\longrightarrow} T_1  \stackrel{\Delta_1}{\longrightarrow}
T_0$$
\newline
to get  
\newline
$$S_2 \stackrel{\partial_2 }{\longrightarrow}  S_1 \stackrel{\partial_1}{\longrightarrow}S_0$$
\newline
so that we have a chain isomorphism:
\newline
\begin{eqnarray*}
S_2 \stackrel{\partial_2 }{\longrightarrow}  S_1 \stackrel{\partial_1}{\longrightarrow}S_0\,\,\\ 
\downarrow \theta_2    \quad      \downarrow \theta_1 \quad   \downarrow \theta_0\\
L_2\stackrel{D_2}{\longrightarrow}  L_1 \stackrel{D_1}{\longrightarrow}L_0\,\,\\ 
\end{eqnarray*}
\newline
}

\bigskip
\hspace{2mm} Proof:  See theorem 1.1.1.

\hfill $\Box$ \,\,

\bigskip

\hspace{2mm}  Note that (12) is chain homotopy equivalent to to 

$$
T_0^* \stackrel{\Delta_1^*}{\longrightarrow}  T_1^*  \stackrel{\Delta_2^*}{\longrightarrow} T_2^*\oplus T \stackrel{\delta_3}
{\longrightarrow} T_2 \oplus T \stackrel{\Delta_2 \oplus 0}{\longrightarrow} T_1  \stackrel{\Delta_1}{\longrightarrow} T_0 
$$

where $\delta_3=\left(\begin{array}{cc} \Delta_3&0\\ 0& 1 \end{array}\right)$. This in turn is chain homotopy equivalent to to

$$
S_0^* \stackrel{\partial_1^*}{\longrightarrow}  S_1^*  \stackrel{\partial_2^*}{\longrightarrow} S_2^* 
\stackrel{\partial_3}{\longrightarrow} S_2   \stackrel{\partial_2}{\longrightarrow} S_1  \stackrel{\partial_1}
{\longrightarrow} S_0 
$$

where $\partial_3$ is induced from $\Delta_3$.  

\bigskip
\hspace{2mm}  Let $D_3 = \theta_2 \circ \partial_3 \circ \theta_2^*$.  Let $\gamma$ be the bilinear form induced on $J_n$ by $D_3$.

We have a chain isomorphism

\begin{eqnarray*}
S_0^* \stackrel{\partial_1^*}{\longrightarrow}  S_1^*  \stackrel{\partial_2^*}{\longrightarrow} S_2^* \, \,\,\,\,
\stackrel{\partial_3}{\longrightarrow} \,\,\quad S_2 \quad  \stackrel{\partial_2}{\longrightarrow} S_1  \stackrel{\partial_1}
{\longrightarrow} S_0 \,\,\\ 
\quad \downarrow \theta_0^{*-1}\,\,\, \downarrow \theta_1^{*-1}
\,\,\, \downarrow \theta_2^{*-1} \qquad \quad  \downarrow \theta_2    \quad \quad      \downarrow \theta_1 \quad   \downarrow \theta_0\\
L_0^* \stackrel{D_1^*}{\longrightarrow}  \,\,L_1^*  \stackrel{D_2^*}{\longrightarrow}\,\, L_2^* \quad 
\stackrel{D_3}{\longrightarrow} \quad L_2 \quad   \stackrel{D_2}{\longrightarrow} L_1  \stackrel{D_1}
{\longrightarrow} L_0 
\,\,\\ 
\end{eqnarray*}

Finally note that the algebraic complex

$$
L_0^* \stackrel{D_1^*}{\longrightarrow}  \,\,L_1^*  \stackrel{D_2^*}{\longrightarrow}\,\, L_2^* \quad 
\stackrel{D_3}{\longrightarrow} \quad L_2 \quad   \stackrel{D_2}{\longrightarrow} L_1  \stackrel{D_1}
{\longrightarrow} L_0 
$$

is obtained from the algebraic 5- complex associated to $(\mathcal{A}^n,\, \gamma)$, by performing four simple homotopy
equivalences.  Hence (12) is chain homotopy equivalent to to the algebraic 5- complex associated to $(\mathcal{A}^n,\, \gamma)$.  

\hfill $\Box$ \,\,

\bigskip

\hspace{2mm} It may be convenient to regard all elements of ${\rm ALG}^5$ as being parametrized, up to homotopy, by 
forms on the same algebraic 2- complex.  To this end, note that for $n \leq m$ the associated algebraic complexes of 
$(\mathcal{A}^n,\, \beta)$ and $(\mathcal{A}^m,\, \gamma)$ are chain homotopy equivalent, where $\beta$ is a bilinear form on $J_n$ and $\gamma$ is
the bilinear form on $J_m$ defined by

$$
\gamma = \left(\begin{array}{cc} \beta &0\\ 0& 1 \end{array}\right)
$$

\hspace{2mm} Hence by identifying $\beta$ with $\gamma$, we may regard $\beta$ as a form on $J_m$, with no ambiguity as to which
element (up to isomorphism in the category) of ${\rm ALG}^5$ is associated to it.  

\sec{Homotopy equivalence in $DUAL^2$}

\hspace{2mm}  We proceed to give necessary and sufficient conditions 
for elements of
$DUAL^2$ to be chain homotopy equivalent.  Let 
$$\mathcal{B} \quad =\quad  J^* \stackrel{}{\to} F_2 \stackrel{d_2}{\to} F_1 
\stackrel{d_1}{\to} F_0  \,$$
$$\mathcal{C} \quad =\quad  K^* \stackrel{}{\to} E_2 
\stackrel{\delta_2}{\to} E_1 
\stackrel{\delta_1}{\to} E_0  \,$$

and let $\beta$, $\gamma$ be $G$- invariant bilinear forms on $\mathcal{B}$,
$\mathcal{C}$ respectively.

\prop{$(\mathcal{B}, \beta)$ and $(\mathcal{C}, \gamma)$ are chain homotopy
equivalent if and only if there exist maps $\phi_1,\psi_2: J^* \to K^*$, $\phi_2,
\psi_1: K^* \to J^*$, which all augment to $\pm1$, and maps $I:J^* \to J$
and $L:K^* \to K$ which factor through projective modules, such that the
following diagram commutes:
\newline
\xymatrix{
J \ar[r]^\beta \ar[d]_{\phi_2^*}&J^* \ar[d]^{\phi_1}\\
K \ar[r]^\gamma \ar[d]_{\psi_2^*}&K^* \ar[d]^{\psi_1}\\
J \ar[r]^\beta &J^* \\
}
\newline
\newline
and 
\newline
\newline
$1-\psi_1 \phi_1=\beta I$, \newline$1-\psi_2^* \phi_2^*=I\beta$, 
\newline
$1-\phi_1 \psi_1=\gamma L$, \newline$1-\phi_2^* \psi_2^*=L\gamma$.
\newline
and $\psi_1$ and $\phi_1$ augment to the same value and $\psi_2$ and $\phi_2$
augment to the same value.
}

\bigskip
\hspace{2mm} Proof:  Firstly suppose we have a pair of 
inverse (up to homotopy) homotopy equivalences  $f_i, g_i$, $i=0,1,2,3,4,5$,
between $(\mathcal{B}, \beta)$ and $(\mathcal{C}, \gamma)$.  Then we may induce 
$\phi_1, \phi_2$ from $f_2, f_3$ and induce $\psi_1,\psi_2$ from $g_1, g_2$. These induced maps all augment to
$\pm1$ and make the diagram commute.  Also and $\psi_1$ and $\phi_1$ will 
augment to the same value and $\psi_2$ and $\phi_2$
will augment to the same value.

\bigskip
\hspace{2mm} We also have maps $I_i, L_i$ such that $1-g_i f_i =I_i\partial+
\partial' I_{i-1}$, and $1-f_i g_i =L_i\partial+
\partial' L_{i-1}$ where $\partial,\, \partial'$ are the relevant boundary maps in
each case and $I_i$ and $L_i$ are taken to be $0$, for $i=-1,5$.  Let 
$\iota,\kappa$ denote the inclusions of $J^*$ in $F_2$, $K^*$ in
$E_2$ respectively.  Then we may construct $I$, $L$, by setting
$I=\iota^*\circ I_2 \circ \iota$ and $L=\kappa^*\circ L_2\circ\kappa$.

\bigskip
\hspace{2mm}  Clearly $I$ and $K$ factor through projective modules.  We know

\bigskip
$\iota \beta I = (\iota \beta\iota^*)I_2 \iota
=(1-g_2f_2-I_1d_2)\iota = \iota-g_2f_2\iota=\iota(1-\psi_1 \phi_1)$

\bigskip
As $\iota$ is injective, we have $1-\psi_1 \phi_1=\beta I$.

\bigskip
\hspace{2mm}  Similarly, we have 

\bigskip
$\kappa \gamma L = (\kappa \gamma\kappa^*)L_2 \kappa
=(1-f_2g_2-L_1 \delta_2)\kappa = \kappa-f_2g_2\kappa=\kappa(1-\phi_1 \psi_1)$

and as $\kappa$ is injective we have $1-\phi_1 \psi_1=\gamma L$.

\bigskip
\hspace{2mm} Also

$I\beta\iota^*=\iota^* I_2
(\iota\beta\iota^*)=\iota^*(1-g_3f_3-d_2^*I_3)=\iota_*-\iota^*g_3f_3=
(1-\psi_2^*\phi_2^*)\iota^*$

\bigskip
As $\iota^*$ is surjective, we have $1-\psi_2^* \phi_2^*=I\beta$.

\bigskip
\hspace{2mm} Again,

$L\gamma\kappa^*=\kappa^* L_2
(\kappa\gamma\kappa^*)=\kappa^*(1-f_3g_3-\delta_2^*L_3)=\kappa_*-\kappa^*f_3g_3
=(1-\phi_2^*\psi_2^*)\kappa^*$

\bigskip
and as $\kappa^*$ is surjective, we have $1-\phi_2^* \psi_2^*=L\gamma$.

\bigskip

Hence all the required conditions are satisfied.

\bigskip
\hspace{2mm}
Conversely, suppose that we have $\phi_1,\phi_2, \psi_1, \psi_2, I, L$
satisfying the required conditions.  We must construct a homotopy equivalence.

\hspace{2mm} We may write $\mathcal{B}$ in the form

$$\quad  J^* \stackrel{}{\to} F_2 \stackrel{d_2}{\to}
S \oplus F 
\stackrel{d_1'\oplus 1_F}{\longrightarrow} \ZG \oplus F  \,$$

for some free module $F$, and stably free module $S$.

\bigskip
\hspace{2mm}  Simple homotopy equivalences
connect  $(\mathcal{B}, \beta)$ and $(\mathcal{B'}, \beta)$  where 

$$\mathcal{B'} \quad =\quad  J^* \stackrel{}{\to} F_2 \oplus F \stackrel
{d_2 \oplus 1_F}{\longrightarrow} F_1 
\stackrel{d_1'}{\to} \ZG  \,$$

\hspace{2mm}Similarly, we may write
$\mathcal{C}$ in the form

$$\quad  K^* \stackrel{}{\to} E_2 \stackrel{\delta_2}{\to}
T \oplus E 
\stackrel{\delta_1'\oplus 1_E}{\longrightarrow} \ZG \oplus E  \,$$

for some free module $E$, and stably free module $T$.

\bigskip
\hspace{2mm}
Simple homotopy equivalences
connect  $(\mathcal{C}, \gamma)$ and $(\mathcal{C'}, \gamma)$  where 

$$\mathcal{C'} \quad =\quad  K^* \stackrel{}{\to} E_2 \oplus E \stackrel
{\delta_2 \oplus 1_E}{\longrightarrow} E_1 
\stackrel{\delta_1'}{\to} \ZG  \,$$

\bigskip
\hspace{2mm}  From now on we will write $\mathcal{B'}$ as 

$$J^* \stackrel{}{\to} F_2' \stackrel{d_2'}{\to} F_1' 
\stackrel{d_1'}{\to} \ZG  \,$$

and $\mathcal{C'}$ as

$$K^* \stackrel{}{\to} E_2' 
\stackrel{\delta_2'}{\to} E_1' 
\stackrel{\delta_1'}{\to} \ZG  \,$$

\bigskip
\hspace{2mm}  We may extend $\psi_1, \psi_2^*, \phi_1, \phi_2^*$ so that the
following diagram commutes:

\bigskip

{\scriptsize
\xymatrix{\Z \ar[d]_{\pm\pm1}\ar[r]&\ZG \ar[r]^{d_1'^*} \ar[d]_{f_5}&
F_1'^* \ar[d]_{f_4} \ar[r]^{d_2'^*}
&F_2'^* \ar[d]_{f_3}\ar[r]^{\iota^*}& J\ar[d]_{\phi_2^*} \ar[r]^{\beta}
& J^* \ar[d]_{\phi_1}\ar[r]^{\iota}&F_2'\ar[d]_{f_2}\ar[r]^{d_2'} &
F_1' \ar[d]_{f_1}\ar[r]^{d_1'}& \ZG \ar[r]\ar[d]_{f_0}& \Z\ar[d]_{\pm1}\\
\Z \ar[d]_{\pm\pm1}\ar[r]&\ZG \ar[r]^{\delta_1'^*} \ar[d]_{g_5}&
E_1'^* \ar[d]_{g_4} \ar[r]^{\delta_2'^*}
&E_2'^* \ar[d]_{g_3}\ar[r]^{\kappa^*}& K\ar[d]_{\psi_2^*} \ar[r]^{\gamma}
& K^* \ar[d]_{\psi_1}\ar[r]^{\kappa}&E_2'\ar[d]_{g_2}\ar[r]^{\delta_2'} &
E_1' \ar[d]_{g_1}\ar[r]^{\delta_1'}& \ZG \ar[r]\ar[d]_{g_0}& \Z\ar[d]_{\pm1}\\
\Z \ar[r]&\ZG \ar[r]^{d_1'^*} &
F_1'^* \ar[r]^{d_2'^*}
&F_2'^* \ar[r]^{\iota^*}& J \ar[r]^{\beta}
& J^* \ar[r]^{\iota}&F_2'\ar[r]^{d_2'} &
F_1' \ar[r]^{d_1'}& \ZG \ar[r]& \Z
}
}

\bigskip
\bigskip

\hspace{2mm}  For some projective module $P$, there exist maps $a:J^* \to P$ and
$b:P \to J$ such that $I=ba$.  As $\iota$ is torsion free and $P$ is projective,
we have maps $a':F_2' \to P$ and $b':P \to F_2'^*$ such that $a=a' \iota$ and
$b=\iota^* b'$.  Let $I_2=b'a'$.  Then we have 

$\iota^* I_2 \iota =\iota^* b'a' \iota=ba=I$

\bigskip
\hspace{2mm}Let $d_3'=\iota\beta\iota^*$ and $\delta_3'=\kappa\gamma\kappa^*$. 
Then

\bigskip

$(1-g_2f_2-d_3'I_2)d_3'= (\iota - g_2f_2\iota
-\iota\beta\iota^*I_2\iota))\beta\iota^*=\iota(1-\psi_1\phi_1-\beta I )
\beta\iota^*=0$

\bigskip
\hspace{2mm}  As $F_2'$ is relatively injective, we may conclude 
that $(1-g_2f_2-d_3'I_2)$ factors through $d_2'$.  Hence we have some map
$I_1:F_1' \to F_2'$,
such that $1-g_2f_2=d_3'I_2+I_1d_2'$.

\bigskip
\hspace{2mm}We have 

$ (1-g_1f_1-  d_2'I_1 )d_2'\quad=\quad d_2'-g_1f_1d_2'-d_2'(1-g_2f_2-d_3'I_2)$

$=d_2'-g_1f_1d_2' -d_2'+g_1f_1d_2' =0$

\bigskip
\hspace{2mm}  Again $F_2'$ is relatively injective so 
$ (1-g_1f_1-  d_2'I_1 )$ factors through $d_1'$.  Hence we have some map$I_0:\ZG
\to F_1'$ such that $ 1-g_1f_1=  d_2'I_1+I_0d_1'$.

\bigskip
\hspace{2mm}  Let $\epsilon:\ZG \to \Z$ denote the augmentation map. From
commutativity of the above diagram, we know that  
$\epsilon(1-f_0g_0)=0$. Hence we know that the image of $(1-f_0g_0)$ lies in $IG
\subset \ZG$.  Clearly the image of $d_1'I_0$ also lies in $IG$.  Hence the
image of $(1-f_0g_0-d_1'I_0)$ lies in $IG$.

\bigskip
\hspace{2mm}As before we have,

\bigskip

$ (1-g_0f_0-  d_1'I_0 )d_1'\quad=\quad d_1'-g_0f_0d_1'-d_1'(1-g_1f_1-d_2'I_1)$

$=d_1'-g_0f_0d_1' -d_1'+g_0f_0d_1' =0$

\bigskip
\hspace{2mm}Hence $(1-g_0f_0-  d_1'I_0 )$ factors through $\epsilon$.  Let $W:\Z
\to \ZG$ satisfy 

\bigskip
$(1-g_0f_0-  d_1'I_0 )= W\epsilon$.  

\bigskip
\hspace{2mm}As $\epsilon$ is
surjective, the image of $W$ must lie in $IG$.  But any map $\Z \to IG$ is
necessarily equal to $0$.  Hence $W=0$ and $1-g_0f_0=d_1'I_0$.

\bigskip
\hspace{2mm}We proceed to construct $I_3$, $I_4$ in a dual fashion.  We have
already shown that $(1-g_2f_2-d_3'I_2)d_3'=0$.  By commutativity of the diagram,
we may conclude that $d_3'(1-g_3f_3-I_2d_3')=0$.  By projectivity of $F_2'^*$,
we have some map $I_3$ satisfying $1-g_3f_3=I_2d_3'+d_2'^*I_3$.

\bigskip
\hspace{2mm}  Again,

\bigskip
$d_2'^*(1-g_4f_4-  I_3 d_2'^*)=d_2'^*-g_3f_3d_2'^*-(1-g_3f_3-I_2d_3')d_2'^*=0$

and by projectivity of $F_1'^*$, we have a map $I_4:F_1'^* \to \ZG$ such that

$1-g_4f_4=I_3d_2'^*+d_1'^*I_4$.

\bigskip
\hspace{2mm}  Repeating the method once more, we have

\bigskip
$d_1'^*(1-g_5f_5-  I_4 d_1'^*)=d_1'^*-g_4f_4d_1'^*-(1-g_4f_4-I_3d_2'^*)d_1'^*=0$

so $(1-g_5f_5-  I_4 d_1'^*)$ factors through $\epsilon^*$, and we have some map
$W':\ZG \to \Z$ satisfying $1-g_5f_5=  I_4 d_1'^*+ \epsilon^*W'$.

\bigskip
\hspace{2mm}  We know that $\epsilon^*W'\epsilon^*=(1-g_5f_5-  I_4
d_1'^*)\epsilon^*=0$.  As $\epsilon^*$ is injective, we may conclude that 
$W'\epsilon^*=0$.  So $W'$ factors through $\ZG/\epsilon(\Z)\cong IG^*$. 
However, any map $IG^* \to \Z$ is necessarily equal to $0$, so $W'=0$ and 
$1-g_5f_5=  I_4 d_1'^*$.

\bigskip
\hspace{2mm}  Collating, we have: 

$1-g_5f_5=  I_4 d_1'^*$

$1-g_4f_4=I_3d_2'^*+d_1'^*I_4$

$1-g_3f_3=I_2d_3'+d_2'^*I_3$

$1-g_2f_2=d_3'I_2+I_1d_2'$

$1-g_1f_1=  d_2'I_1+I_0d_1'$

$1-g_0f_0=d_1'I_0$

Hence the $I_i$, for $i=0,1,2,3,4$, form a chain homotopy from the identity
 to $g_if_i$.

\hspace{2mm}Similarly we may construct maps $L_0, L_1, L_2, L_3,L_4$, which form a
chain homotopy between the identity and the $f_ig_i$, $i=0,1,2,3,4,5$.

\bigskip
\hspace{2mm} Hence $(\mathcal{B}, \beta)$ is chain homotopy equivalent to 
$(\mathcal{B'}, \beta)$, which is chain homotopy equivalent to 
$(\mathcal{C'}, \gamma)$, which is chain homotopy equivalent to 
$(\mathcal{C}, \gamma)$, as required.

\hfill $\Box$ \,\,

\bigskip
\hspace{2mm}Note that this proposition does not imply that the underlying
algebraic $2$- complex is not important in defining the isomorphism class of an
object in $DUAL \pm^2$, as whether or not a map $J^* \to K^*$ augments to
$\pm1$, is dependent on the precise inclusions of $J^*$, $K^*$, in $F_2$, $E_2$
respectively.

\lem{ Let $(\mathcal{B}, \beta)$, be as before.  Suppose we have some homotopy 
equivalence $f_*$, from $(\mathcal{B}, \beta)$ to its dual.  Then $\pm f_*$ is chain homotopic to a chain map of the form:
\newline
\begin{eqnarray*}
F_0^* \stackrel{d_1^*}{\longrightarrow} F_1^* \stackrel{d_2^*}{\longrightarrow} F_2^* \stackrel{d_3}{\longrightarrow} 
F_2 \stackrel{d_2}{\longrightarrow} F_1 \stackrel{d_1}{\longrightarrow} F_0\,\,\,\,\,\\
\downarrow \pm1  \quad\downarrow \pm1 \quad\downarrow \theta_2 \,\,\,\quad \downarrow \theta_1 \,\,\quad\downarrow 1
\qquad\downarrow 1\\
F_0^* \stackrel{d_1^*}{\longrightarrow} F_1^* \stackrel{d_2^*}{\longrightarrow} F_2^* \stackrel{d_3^*}{\longrightarrow} 
F_2 \stackrel{d_2}{\longrightarrow} F_1 \stackrel{d_1}{\longrightarrow} F_0\,\,\,\,\,\\ 
\end{eqnarray*}
\newline
for some maps $\theta_1:\,F_2\,\to\,F_2$, $\theta_2:\,F_2^*\,\to\,F_2^*$.  
}

\bigskip
\hspace{2mm}  Proof:  Consider the chain homotopy equivalence:

\begin{eqnarray*}
F_0^* \stackrel{d_1^*}{\longrightarrow} F_1^* \stackrel{d_2^*}{\longrightarrow} F_2^* \stackrel{d_3}{\longrightarrow} 
F_2 \stackrel{d_2}{\longrightarrow} F_1 \stackrel{d_1}{\longrightarrow} F_0\,\,\,\,\,\\
\downarrow \alpha_5  \,\,\quad\downarrow \alpha_4 \quad\downarrow \alpha_3 \,\quad \downarrow \alpha_2 \,\,\quad\downarrow \alpha_1
\quad\downarrow \alpha_0\\
F_0^* \stackrel{d_1^*}{\longrightarrow} F_1^* \stackrel{d_2^*}{\longrightarrow} F_2^* \stackrel{d_3^*}{\longrightarrow} 
F_2 \stackrel{d_2}{\longrightarrow} F_1 \stackrel{d_1}{\longrightarrow} F_0\,\,\,\,\,\\ 
\end{eqnarray*}

where the $\alpha_i$ are equal to the $f_i$ or the $-f_i$, depending on which is necessary to force the induced action on the cokernel of
$d_1$ to be the identity.

\hspace{2mm}  The
kernel of $d_1^*$ and the cokernel of $d_1$ are both $\Z$.  Hence $\alpha_0$ and $\alpha_5$ must induce multiplication by $\pm1$ on
$\Z$.  Our choice of sign forces $\alpha_0$ to induce multiplication by $1$.

\begin{eqnarray*}
J^*\stackrel{l}{\longrightarrow} 
F_2 \stackrel{d_2}{\longrightarrow} F_1 \stackrel{d_1}{\longrightarrow} F_0 \to \Z\,\,\,\,\,\\
\downarrow \phi_1 \,\quad \downarrow \alpha_2 \,\,\quad\downarrow \alpha_1
\quad\downarrow \alpha_0\,\,\, \downarrow 1\\
J^*\stackrel{l}{\longrightarrow} 
F_2 \stackrel{d_2}{\longrightarrow} F_1 \stackrel{d_1}{\longrightarrow} F_0\to \Z\,\,\,\,\,\\ 
\end{eqnarray*}

\begin{eqnarray*}
J^*\stackrel{l}{\longrightarrow} 
F_2 \stackrel{d_2}{\longrightarrow} F_1 \stackrel{d_1}{\longrightarrow} F_0 \to \Z\,\,\,\,\,\,\,\,\,\\
\downarrow \phi_2 \,\quad \downarrow \alpha_3^* \,\,\quad\downarrow \alpha_4^*
\quad\downarrow \alpha_5^*\,\,\, \downarrow \pm 1\\
J^*\stackrel{l}{\longrightarrow} 
F_2 \stackrel{d_2}{\longrightarrow} F_1 \stackrel{d_1}{\longrightarrow} F_0\to \Z\,\,\,\,\,\,\,\,\,\\
\end{eqnarray*}

\hspace{2mm}  Hence we have maps $I_0:F_0 \to F_1$ and $I_1:F_1 \to F_2$, such that $d_1 I_0 =1- \alpha_0$, and
$I_0 d_1+d_2 I_1 =1- \alpha_1$.

\hspace{2mm}  Similarly we have maps $I_4:F_1^* \to F_0^*$ and $I_3:F_2^* \to F_1^*$, such that $d_1 I_4^* =\pm1- \alpha_5^*$, and
$I_4^* d_1+d_2 I_3^* =\pm1- \alpha_4^*$.

\bigskip
\hspace{2mm} Next set $\theta_1= \alpha_2+I_1d_2$ and $\theta_2= \alpha_3+ d_2^*I_3$.  Then taking $I_2=0$, the $I_i$, $i=0,1,2,3,4$ form the
required chain homotopy.

\hfill $\Box$ \,\,

\sec{Poincare Duality}

\hspace{2mm}Let $X$ be a CW- complex with finite fundamental group $G$. Let $L_*(\tilde{X})$ denote $C_*(\tilde{X})$ with coefficients
restricted to $\Z$.  As $G$ is finite, we may naturally identify $H_p(\tilde{X};\ZG)$ with $H_p(L_*(\tilde{X});\Z)$ and we may 
identify $H^p(\tilde{X};\ZG)$ with $H^p(L_*(\tilde{X});\Z)$.  We make frequent use of these identifications, and hence assume natural maps 
$H^p(\tilde{X};\ZG) \times H_p(\tilde{X};\ZG) \to \Z$.  We now make a more specific statement of Poincare Duality:

\thm{(Poincare Duality) Let $M$ be an element of ${\rm TOP}^5$ and $M'$ be homotopy equivalent to $M$.  Then, given a generator of 
$H_5(\tilde{M'};\ZG)$, denoted $\eta$, there exists a chain homotopy equivalence, over $\ZG$, $\phi:C_*(\tilde{M'})^* \to C_*(\tilde{M'})$
satisfying the following:  Given $\alpha \in H^p(\tilde{M'};\ZG)$, we have $\phi_*(\alpha) = \eta \,\,\,\widehat{}\,\,\, \alpha$. Here $\phi_*$ denotes the
induced map $H^p(\tilde{M'};\ZG) \to H_{5-p}(\tilde{M'};\ZG)$.
}

\cor{Given $p\in \{ 0,1,2,3,4,5\} $ and $\alpha \in H^p(\tilde{M'};\ZG)$, $\beta \in H^{5-p}(\tilde{M'};\ZG)$ we have 
$\beta(\phi_*(\alpha)) = \alpha(\phi_*(\beta))$.
}

\proof{$$\beta(\phi_*(\alpha))=\beta(\eta \,\,\,\widehat{}\,\,\, \alpha)= (\alpha \,\,\,\breve{}\,\,\, \beta) \eta = 
{-1}^{p(5-p)} (\beta \,\,\,\breve{}\,\,\, \alpha)\eta =
{-1}^{p(5-p)}\alpha(\eta \,\,\,\widehat{}\,\,\, \beta)={-1}^{p(5-p)}\alpha(\phi_*(\beta))$$  
\hspace{2mm} As either $p$ or $5-p$ must be even, we have 
$\beta(\phi_*(\alpha)) = \alpha(\phi_*(\beta))$
}

\hspace{2mm}Suppose now, that we have  chain homotopy equivalence $f:C_*(\tilde{M'}) \to \mathcal{A}$, for some algebraic 5- complex
$\mathcal{A}$.  Then we have a homotopy equivalence $f \circ \phi \circ f^*:\mathcal{A}^* \to \mathcal{A}$.  Let 
$\phi'=f \circ \phi \circ f^*$.  If a chain homotopy equivalence $\mathcal{A}^* \to \mathcal{A}$ is chain homotopic to one constructed in this
way, starting with some generator of $H_5(\tilde{M'};\ZG)$, we say it is a duality equivalence.  Note that (up to sign) the maps from
cohomology to homology, induced by a duality equivalence, are determined by $M$ (with respect to the isomorphisms $f_*$ and $f^*$).

\lem{Given $p\in \{ 0,1,2,3,4,5\} $ and $\alpha \in H^p(\mathcal{A};\ZG)$, $\beta \in H^{5-p}(\mathcal{A};\ZG)$ we have 
$\beta(\phi'_*(\alpha)) = \alpha(\phi'_*(\beta))$.
}

\proof{$\beta(\phi'_*(\alpha)) = \beta(f_* \phi_* f^*(\alpha))=f^*(\beta)(\phi_*f^*(\alpha))=
f^*(\alpha)(\phi_*f^*(\beta))$
\newline
$=\alpha(f_* \phi_* f^*(\beta))
=\alpha(\phi'_*(\beta))$}

\hspace{2mm}By theorem 4.2.1 we may choose $\mathcal{A}$ to be of the form $(\mathcal{B},\beta)$ for some algebraic $2$- complex 
$\mathcal{B}= J^* \stackrel{\iota}\to F_2 \stackrel{d_2} {\to}F_1 \stackrel{d_1}{\to} F_0 \to Z$ and bilinear form $
\beta: J \times J \to \Z$.  Let $d_3=\iota^*\beta\iota$.  The algebraic complex $(\mathcal{B},\beta)$ is written:

$$
\Z \dashrightarrow F_0^* \stackrel{d_1^*} {\to}F_1^* \stackrel{d_2^*}{\to} F_2^* \stackrel{d_3}{\to}
 F_2 \stackrel{d_2} {\to}F_1 \stackrel{d_1}{\to} F_0 \dashrightarrow \Z
$$

\hspace{2mm}As $F_0^*$ is the dual of $F_0$, we may apply elements of $\Z$ (occurring on the left of this sequence), to $\Z$ 
(occurring on the right of this sequence).  We follow the convention that the choice of identification of the kernel of $d_1^*$ with $\Z$,
forces this application to be given by multiplication.

\hspace{2mm}Note that the copy of $\Z$ occurring on the left of the sequence may be identified with $H_5((\mathcal{B},\beta);\ZG)$. 
Similarly, the copy of $\Z$ occurring on the right of the sequence may be identified with $H_0((\mathcal{B},\beta);\ZG)$.

\hspace{2mm}The complex $(\mathcal{B},\beta)^*$ may be written:

$$
\Z \dashrightarrow F_0^* \stackrel{d_1^*} {\to}F_1^* \stackrel{d_2^*}{\to} F_2^* \stackrel{d_3^*}{\to}
 F_2 \stackrel{d_2} {\to}F_1 \stackrel{d_1}{\to} F_0 \dashrightarrow \Z
$$

\hspace{2mm}This sequence only differs from $(\mathcal{B},\beta)$ in the middle term.

\hspace{2mm}Again, note that the copy of $\Z$ occurring on the left of this sequence may be identified with $H^0((\mathcal{B},\beta);\ZG)$. 
Similarly, the copy of $\Z$ occurring on the right of this sequence may be identified with $H^5((\mathcal{B},\beta);\ZG)$.

\hspace{2mm}Let $\alpha \in H^0((\mathcal{B},\beta);\ZG)$ be represented by the integer $a$, and let 
$\gamma \in H_0((\mathcal{B},\beta);\ZG)$ be represented by the integer $c$. Then our conventions imply that application of elements 
of $H^0((\mathcal{B},\beta);\ZG)$ to $H_0((\mathcal{B},\beta);\ZG)$ is given by the application of $\Z$ to $\Z$, which was forced to be
multiplication.  Hence $\alpha ( \gamma) =ac$.

\hspace{2mm}Let $\alpha \in H^5((\mathcal{B},\beta);\ZG)$ be represented by the integer $a$, and let 
$\gamma \in H_5((\mathcal{B},\beta);\ZG)$ be represented by the integer $c$. Then our conventions imply that the application of elements 
of $H^5((\mathcal{B},\beta);\ZG)$ to $H_5((\mathcal{B},\beta);\ZG)$ is given by evaluation on the application of elements of $\Z$ to $\Z$, 
which was forced to be multiplication.  Hence $\alpha (\gamma) =ca=ac$.

\hspace{2mm}Let $x,y \in \Z$ be chosen so that the following diagram commutes:

\begin{eqnarray*}
\Z \dashrightarrow F_0^* \stackrel{d_1^*} {\to}F_1^* \stackrel{d_2^*}{\to} F_2^* \stackrel{d_3^*}{\to}
 F_2 \stackrel{d_2} {\to}F_1 \stackrel{d_1}{\to} F_0 \dashrightarrow \Z\,\,\,\\
y \downarrow \,\,\,\,\, \phi'_5\downarrow
\,\,\,\,\phi'_4\downarrow \,\,\,\,\phi'_3\downarrow\,\,\, \phi'_2\downarrow \,\,\,
\phi'_1\downarrow\,\,\, \phi'_0 \downarrow\quad\, \,\, x\downarrow\\
\Z \dashrightarrow F_0^* \stackrel{d_1^*} {\to}F_1^* \stackrel{d_2^*}{\to} F_2^* \stackrel{d_3}{\to}
 F_2 \stackrel{d_2} {\to}F_1 \stackrel{d_1}{\to} F_0 \dashrightarrow \Z\,\,\,
\end{eqnarray*}

\newpage

\lem{x=y}

\hspace{2mm} Proof:  {Let $\alpha \in H^5((\mathcal{B},\beta);\ZG)$ be represented by the integer $a$, and let 
$\gamma \in H^0((\mathcal{B},\beta);\ZG)$ be represented by the integer $c$.  By lemma 4.5.3 we know 
\newline$\alpha(\phi_* (\gamma))=
\gamma(\phi_* (\alpha))$.  Hence we have $ayc=\alpha(\phi_* (\gamma))=
\gamma(\phi_* (\alpha))=cxa$.
\newline
${}$\hspace{2mm}As $\alpha$ and $\gamma$ were picked arbitrarily, we must have $x=y$.
}

\hfill $\Box$

\hspace{2mm}As $\phi'$ is a homotopy equivalence, we must have $x=\pm1$.  If $x=-1$, we can replace $\eta$ with $-\eta$.  Hence without loss
of generality we have a duality equivalence:

\begin{eqnarray*}
\Z \dashrightarrow F_0^* \stackrel{d_1^*} {\to}F_1^* \stackrel{d_2^*}{\to} F_2^* \stackrel{d_3^*}{\to}
 F_2 \stackrel{d_2} {\to}F_1 \stackrel{d_1}{\to} F_0 \dashrightarrow \Z\,\,\,\\
1 \downarrow \,\,\,\,\, \phi'_5\downarrow
\,\,\,\,\phi'_4\downarrow \,\,\,\,\phi'_3\downarrow\,\,\, \phi'_2\downarrow \,\,\,
\phi'_1\downarrow\,\,\, \phi'_0 \downarrow\quad\, \,\, 1\downarrow\\
\Z \dashrightarrow F_0^* \stackrel{d_1^*} {\to}F_1^* \stackrel{d_2^*}{\to} F_2^* \stackrel{d_3}{\to}
 F_2 \stackrel{d_2} {\to}F_1 \stackrel{d_1}{\to} F_0 \dashrightarrow \Z\,\,\,
\end{eqnarray*}

\hspace{2mm}Hence by lemma 4.4.2, for some maps $\theta_1$, $\theta_2$ we have a duality equivalence:

\begin{eqnarray*}
F_0^* \stackrel{d_1^*}{\longrightarrow} F_1^* \stackrel{d_2^*}{\longrightarrow} F_2^* \stackrel{d_3^*}{\longrightarrow} 
F_2 \stackrel{d_2}{\longrightarrow} F_1 \stackrel{d_1}{\longrightarrow} F_0\,\,\,\,\,\\
\downarrow 1\,\,\,\,\,\,  \quad\downarrow 1 \,\,\,\,\quad\downarrow \theta_2 \,\,\,\quad \downarrow \theta_1 \,\,\quad\downarrow 1
\qquad\downarrow 1\\
F_0^* \stackrel{d_1^*}{\longrightarrow} F_1^* \stackrel{d_2^*}{\longrightarrow} F_2^* \stackrel{d_3}{\longrightarrow} 
F_2 \stackrel{d_2}{\longrightarrow} F_1 \stackrel{d_1}{\longrightarrow} F_0\,\,\,\,\,\\ 
\end{eqnarray*}

\hspace{2mm}A natural question to ask at this stage is, whether or not we can choose $(\mathcal{B},\beta)$ such that the identity map
is a duality equivalence.  If 
this is possible for a manifold, we say that it is self dual.  In the next 
section, we will show that not all manifolds in 
${\rm TOP}^5$ are self dual.

\hspace{2mm}Note that if the identity map is a chain map $(\mathcal{B},\beta)^* \to (\mathcal{B},\beta)$, then 

$d_3^*=1d_3^*=d_3 1=d_3$.  $d_3=\iota^*\beta\iota$.

\hspace{2mm}As $d_3=\iota^*\beta\iota$, with $\iota$ injective and $\iota^*$ surjective, we have $\beta^*=\beta$.  So $\beta$ must in this
case be symmetric.

\sec{Linking Number}  

\hspace{2mm}In this section we construct the linking number. We show that if a manifold is self dual, then its linking number
is symmetric.  We then use the fact that the linking number is always antisymmetric, to show that a self dual manifold, $M$, must satisfy
Tor$(H_2(\tilde{M};\ZG))={\Z_2}^k$, for some integer $k$.  Finally we refer to the existence of manifolds not satisfying this homology
condition, which are therefore not self dual.

\hspace{2mm}Let $(F_i, d_i)$, $i=0,1,2,3,4,5$, be a free, finite, algebraic $5$- complex over $\ZG$.  Then an element of Tor$(H^3(F_i, \ZG))$ 
may be represented by a map $F_3 \to \Z$, of the form $\frac{1}{a}wd_3$, for some integer $a$ and $w:F_2 \to \Z$.

\hspace{2mm}Similarly, an element of Tor$(H_2(F_i, \ZG))$ may be represented by an element of $F_2$ of the form $\frac{1}{b}d_3x$, for some
integer $b$ and $x \in F_3$.

\hspace{2mm}Hence we have a bilinear map Tor$(H^3(F_i, \ZG)) \times $Tor$(H_2(F_i, \ZG)) \to \Q / \Z$, given by $$
(\frac{1}{a}wd_3, \frac{1}{b}d_3x) \mapsto \frac{1}{ab}w(d_3x)
$$

\hspace{2mm}Let Z denote this bilinear map.

\lem{This bilinear map is well defined.}

\hspace{2mm}Proof:  Suppose we had made a different choice of representative to $\frac{1}{a}wd_3$.  Then the different choice 
would differ by a map of the form $vd_3$, for some map $v:F_2 \to \Z$.  We have $$
(\frac{1}{a}wd_3 +vd_3, \frac{1}{b}d_3x) \mapsto \frac{1}{ab}w(d_3x) + v(\frac{1}{b}d_3 x)
$$ with $v(\frac{1}{b}d_3 x) \in \Z$.
\newline
${}$\hspace{2mm}Again, suppose we had made a different choice of representative to $\frac{1}{b}d_3x$.  Then the different choice 
would differ by an element of the form $d_3y$, for some $y\in :F_3$.  We have $$
(\frac{1}{a}wd_3 , \frac{1}{b}d_3x+d_3y) \mapsto \frac{1}{ab}w(d_3x) + \frac{1}{a}wd_3 y
$$ with $\frac{1}{a}wd_3 y \in \Z$.
\newline
${}$\hspace{2mm}Either way, the value of the bilinear map, in $\Q / \Z$ is unchanged.

\hfill $\Box$.

\hspace{2mm}Suppose now we have another free, finite algebraic $5$- complex, $(E_i, \delta_i)$ and a homotopy equivalence $
f:(F_i, d_i) \to (E_i, \delta_i)$.  Let $\frac{1}{a}w\delta_3$ represent an element of Tor$(H^3(E_i, \ZG))$ and let $\frac{1}{b}d_3x$ represent an
element of Tor$(H_2(F_i, \ZG))$.  We have: 

\begin{eqnarray*}
f^*(\frac{1}{a}w\delta_3)=\frac{1}{a}f^*(w\delta_3)=\frac{1}{a}w\delta_3f_3  =\frac{1}{a}wf_2d_3 =\frac{1}{a}wf_2d_3
\end{eqnarray*}

and

\begin{eqnarray*}
f_*(\frac{1}{b}d_3x)=\frac{1}{b}f_2(d_3x) =\frac{1}{b}f_*(d_3x)=\frac{1}{b}\delta_3 f_3x
\end{eqnarray*}

\hspace{2mm}Hence:
\begin{eqnarray*}
Z(f^*(\frac{1}{a}w\delta_3), \frac{1}{b}d_3x)=  Z(\frac{1}{a}(wf_2 d_3, \frac{1}{b}d_3x)
=\frac{1}{ab}wf_2(d_3x) \\  \\=\frac{1}{ab}w\delta_3f_3x = Z(\frac{1}{a}w\delta_3, \frac{1}{b}\delta_3 f_3x)
=Z(\frac{1}{a}w\delta_3, f_*(\frac{1}{b}d_3x))
\end{eqnarray*}

\hspace{2mm}So the maps $f_*$ and $f^*$ are adjoint with respect to $Z$.

\hspace{2mm} Suppose now that $M'$ is homotopic to some element of ${\rm TOP}^5$, and $(\mathcal{B},\beta)$ is chain homotopy equivalent to 
$C_*(\tilde{M'})$, via a homotopy equivalence $f:C_*(\tilde{M'})\to (\mathcal{B},\beta)$.  Let $\phi$ denote a 
duality equivalence $(\mathcal{B},\beta)^* \to (\mathcal{B},\beta)$ and let $\psi$ be a homotopy inverse to it.  

\hspace{2mm}Let $\Theta$ denote ${\rm Tor}(H_2((\mathcal{B},\beta);\ZG))$.  We define the linking number on $(\mathcal{B},\beta)$ to be the 
bilinear map $\nu:\Theta \times \Theta \to \Q/\Z$, given by:
$$
\nu(x,y)=Z(\psi_*x,y)
$$

\hspace{2mm}Let $h$ denote a homotopy inverse to $f$.  Suppose we have a homotopy 
equivalence $g:C_*(\tilde{M'})\to (\mathcal{C},\gamma)$.  Let $k$ denote a homotopy inverse to it.  Then a 
duality equivalence  $(\mathcal{C},\gamma)^*\to (\mathcal{C},\gamma)$ is given by $gh\phi h^*g^*$.  A homotopy inverse to that is given by
$k^*f^*\psi fk$.  Let $m=fk$.

\hspace{2mm}Now let $\Theta'={\rm Tor}(H_2((\mathcal{C},\gamma);\ZG))$ and let $\nu'$ be the linking 
number 

$\Theta' \times \Theta' \to \Q/\Z$.  We have:

$$
\nu'(x,y)=Z(m^* \psi_* m_*x,y)=\nu'(x,y)=Z(\psi_* m_*x,m_*y)=\nu(m_*x,m_*y)
$$

by the adjointness property.  Hence we have that the linking number is well defined up to isomorphism.  In particular, such properties as
being symmetric or antisymmetric are independent of the choice of $(\mathcal{B},\beta)$.

\thm{(see \cite{Bard}, Lemma D(ii) ) $\nu$ is antisymmetric.}

\hspace{2mm}  Note that \cite{Bard} is concerned with simply connected, closed, orientable 5- manifolds.  However, the fact that $G$ is finite 
means that $\tilde{M}$ will be closed, simply connected and orientable.

\hspace{2mm}We now describe $\nu$ in terms of $\beta$. We know that $H_2((\mathcal{B},\beta);\ZG)$ is the cokernel of $\beta:J \to J^*$. 
Hence $\Theta={\rm Tor}(H_2((\mathcal{B},\beta);\ZG))$ consists of maps $J \to \Z$ of the form $\frac{1}{a}\beta(x,\_)$, with $a \in \Z$, $x \in J$.

\hspace{2mm}  Let $H= \{h \in J \otimes \Q | \quad\beta(h,x) \in \Z \quad\forall x \in J \}$, where we take the natural extension of 
$\beta$ to $J \otimes \Q$.  This extension restricts to a form on $H$.  Any element of $\Theta$ may then be written in the form $\beta(h,
\_)$, for some $h \in H$.

\hspace{2mm}  Let $K= \{k \in H | \quad\beta(k,x) = 0 \quad\forall x \in J \}$.  Then the set of maps of the form  
$\beta(h,\_)$, $h \in H$, are naturally identified with $H / K$. Two such maps are homologous precisely if they differ by a map of the form 
$\beta(x,\_)$, $x \in J$. Such maps are naturally identified with $J/K$.  Consequently, we have $\Theta = H / (J + K)$.

\hspace{2mm}Let $\lambda_1:J \to J$ be the map induced by $\psi_3$ and let  $\lambda_2:J \to J$ be the dual of the map induced by $\psi_2$. 
We will now compute $\nu(x/a,y/b)$, for $x,y \in J$, $a,b \in \Z$.

\bigskip
\hspace{2mm}$x/a$ represents the element of $\Theta$ corresponding to the map $\frac{1}{a}\beta(x,\_)$.  Hence 

$\psi_*(x/a)=\frac{1}{a}\beta(x,\lambda_2 \_)$.  Applying this to $y$ and dividing by $b$ gives: $$
\nu(x/a,y/b)=\frac{1}{ab}\beta(x,\lambda_2( y))
$$

\hspace{2mm}From the commutativity of the square 

\xymatrix{J \ar[d]_{\lambda_1}\ar[r]^\beta &J^*\ar[d]^{\lambda_2^*}\\
J\ar[r]^{\beta^*}&J^*}

\bigskip
we know that 
$\beta(x, \lambda_2 (y))= \beta(y, \lambda_1 (x))$, for all $x, y \in J$.  Hence 

$\beta(x, \lambda_2(y))+\beta(y, \lambda_2(x)) =
\beta(x,(\lambda_1 + \lambda_2)(y))$.

\hspace{2mm}The fact that $\nu$ is antisymmetric may therefore be stated as follows:  

\lem{If $a |\beta(x,\_)$ and $b |\beta(y,\_)$ then 
$ab|\beta(x,(\lambda_1 + \lambda_2)(y))$.}

\bigskip
\hspace{2mm}We now suppose that $M$ is self dual.  We may choose $(\mathcal{B}, \beta)$ such that we may take $\phi$ and $\psi$ to 
be identity chain maps.  In this case $\beta$ will be symmetric and $\lambda_1=\lambda_2=1$.  Hence from lemma 4.6.3 we have that 
if $a |\beta(x,\_)$ and $b |\beta(y,\_)$ then $ab|2\beta(x,y)$.

\hspace{2mm}  We will denote the kernel of $\beta$ by $H_3$.  The inclusion of
$H_3$ in $J$ splits over $\Z$. Let $V$ denote a complementary space.  By combining a basis of $H_3$, with a basis of $V$, we obtain a basis
 of $J$, with respect to which, we may represent $\beta$ by a matrix $B$.  Then $B$ will have the form 

$$
B=\left(\begin{array}{cc} 0&0\\ 0& C \end{array}\right)
$$
  
where $\rm{Det}(C) \neq 0$.

\hspace{2mm}  The condition $m$ divides
$\beta(x,y)$ for all $y \in J$, is equivalent to $m|Bx$.  Hence we have:

\lem{$m|Bx$ and $l|By$ imply that $ml|2x^TBy\,\,$ for all integers $m,\, l$ and $x,\,y \in J$.}

\lem{$m|Cx$ and $l|Cy$ imply that $ml|2x^TCy$ for all integers $m,\, l$ and $x,\,y \in V$.}

\bigskip
\hspace{2mm}  Proof:  If $m|Cx$ and $l|Cy$ then $m|Bx'$ and $l|By'$, where $x'=(0,x)$ and $y'=(0,y)$.  Consequently $ml|2x'^TBy'$. But
$2x'^TBy' = 2x^TCy$ so $ml|2x^TCy$. 

\hfill $\Box$ \,\,

\bigskip
\hspace{2mm}  Note that $C$ is invertible over $\Q$ and symmetric.  Note also, that $\Theta = {{\rm coker}(C)}$.

\bigskip
\lem{If $\Theta$ has a non-trivial element of order $k$, then there exists some vector $x \in V$, such 
that $k|Cx$ and k, x are coprime.}

\bigskip
\hspace{2mm}  Proof:  Some vector must have order $k$ modulo the columns of $C$.  Hence multiplying that vector by $k$, gives $Cx$ for some
$x$.  If $l$ were some non-trivial divisor of $x$ and $k$, then $Cx/l$ would be in the image of $C$ and our original vector would have order
less than or equal to $k/l$. 

\hfill $\Box$ \,\,

\bigskip
\hspace{2mm}  Suppose $\Theta$ has a non-trivial element of order $k$.  Let $x'$ denote $x$ factored out by the highest common factor of 
the components 
of $x$.  As $k$ is coprime to $x$, we still
have $k|Cx'$.  We may extend the vector ${x'}$ to a basis of $V$.  Let $D$ denote the matrix representing the bilinear form represented by
$C$, with respect to the new basis.

\bigskip
\hspace{2mm} We know that $D$ is a symmetric matrix with non-zero determinant.  Also we know that if $m|Dx$ and $l|Dy$ then $ml|2x^TDy$, for
any $x,y \in J$.  Finally, we know that the first row and the first column of $D$ are divisible by $k$.

\bigskip
\lem{$k=2$}

\bigskip
\hspace{2mm}  Proof: Let $e_1$ denote 
$\left(\begin{array}{c} 1\\ 0\\0\\.\\.\\.\\.\\0 \end{array}\right)$

\bigskip

\hspace{2mm} As $D$ has non-zero determinant, we may choose some vector $v$ such that $Dv=me_1$ for some positive integer $m$.  

\hspace{2mm}  We have $v^TDe_1 = m$.  As $m|Dv$ and $k|De_1$, we know that $km$ divides $2m$.  $m$ is positive, so $k=2$.

\hfill $\Box$ \,\,

\bigskip
\bigskip
\cor{$\Theta = {\Z_2}^s$, for some $s$.}

\bigskip
\hspace{2mm}  Proof: $\Theta$ is a finitely generated Abelian group, whose non-trivial elements all have order 2.

\hfill $\Box$ \,\,

\hspace{2mm}  We have shown that if $M$ is self dual, then $H_2(\tilde{M};ZG)=\Z^r \oplus \Z_2^s$ for some integers $r$ and $s$.  To show
that elements of ${\rm TOP}^5$ are not always self dual, we need only show a manifold which does not satisfy this homological condition.

\hspace{2mm}We refer again to \cite{Bard}.

\prop{(See \cite{Bard}, lemma 1.1(i)) For each integer $k >1$, there exists a simply connected, closed, manifold $M_k$, with 
$H_2(M_k;\Z)=\Z_k \oplus \Z_k$.}

\bigskip

\hspace{2mm}Lemma 4.5.4 tells us that for any manifold, we may choose $(\mathcal{B}, \beta)$ so that we may take both $\lambda_1$ and
$\lambda_2$ to augment to $1$.  If the manifold is simply connected (as above), then  we are working over $\Z$, so $J$ is free and 
the augmentation condition is vacuous.  However, it is possible that some condition on $G$ could force closed, orientable, simply
connected manifolds with fundamental group $G$, to be self dual.






\addcontentsline{toc}{part}{Bibliography}
\bibliographystyle{plain}
\bibliography{bibthesis}

\end{document}